\documentclass[a4paper,11pt]{article}

\usepackage{amsmath,amstext,amscd,amsthm,amsfonts,diagrams}
\usepackage{amssymb}
\newtheorem{theorem}{Theorem}[section]
\newtheorem{lemma}[theorem]{Lemma}
\newtheorem{proposition}[theorem]{Proposition}

\newtheorem{corollary}[theorem]{Corollary}
\theoremstyle{remark}
\newtheorem{rk}[theorem]{Remark}

\newcommand{\quot}{\ensuremath{/ \hspace{-1.2mm}/}}

\def\Ad{\mathop{\rm Ad}\nolimits}
\def\ad{\mathop{\rm ad}\nolimits}
\def\Int{\mathop{\rm Int}\nolimits}
\def\Aut{\mathop{\rm Aut}\nolimits}
\def\Char{\mathop{\rm char}\nolimits}
\def\Lie{\mathop{\rm Lie}\nolimits}
\def\GL{\mathop{\rm GL}\nolimits}
\def\SL{\mathop{\rm SL}\nolimits}

\def\SO{\mathop{\rm SO}\nolimits}
\def\O{\mathop{\rm O}\nolimits}
\def\Sp{\mathop{\rm Sp}\nolimits}

\def\charac{\mathop{\rm char}\nolimits}
\def\dim{\mathop{\rm dim}\nolimits}
\def\codim{\mathop{\rm codim}\nolimits}
\def\Id{\mathop{\rm Id}\nolimits}
\def\Hom{\mathop{\rm Hom}\nolimits}

\def\ker{\mathop{\rm ker}\nolimits}

\def\tr{\mathop{\rm tr}\nolimits}
\def\rank{\mathop{\rm rk}\nolimits}
\def\min{\mathop{\rm min}\nolimits}

\def\rank{\mathop{\rm rk}\nolimits}
\def\Spec{\mathop{\rm Spec}\nolimits}

\def\diag{\mathop{\rm diag}\nolimits}
\def\pr{\mathop{\rm pr}\nolimits}

\title{Involutions of reductive Lie algebras in positive characteristic}
\author{Paul Levy \\
plevy@imf.au.dk}

\begin{document}

\maketitle{}{}

\maketitle

%\newpage

\begin{abstract}
\noindent Let $G$ be a reductive group over a field $k$ of characteristic
$\neq 2$, let ${\mathfrak g}=\Lie(G)$, let $\theta$ be an involutive
automorphism of $G$ and let ${\mathfrak g}={\mathfrak
k}\oplus{\mathfrak p}$ be the associated symmetric space
decomposition.
For the case of a ground field of characteristic zero, the action of
the isotropy group $G^\theta$ on ${\mathfrak p}$ is well understood,
since the well-known paper of Kostant and Rallis \cite{kostrall}.
Such a theory in positive characteristic has proved more difficult to develop.
Here we use an approach based on some tools from geometric invariant
theory to establish corresponding results in (good) positive characteristic.

Among other results, we prove that the variety ${\cal N}$ of nilpotent
elements of ${\mathfrak p}$ has a dense open orbit, and that the same
is true for every fibre of the quotient map ${\mathfrak p}\rightarrow{\mathfrak p}\quot G^\theta$.
However, we show that the corresponding statement for $G$,
conjectured by Richardson, is not true.
We provide a new, (mostly) calculation-free proof of the number of
irreducible components of ${\cal N}$, extending a result of Sekiguchi
for $k={\mathbb C}$.
Finally, we apply a theorem of Skryabin to describe the infinitesimal
invariants $k[{\mathfrak p}]^{\mathfrak k}$.
\end{abstract}

%\newpage

\section{Introduction}
\label{intro}
Let $G$ be a reductive algebraic group over the algebraically closed field $k$ of characteristic $p\neq 2$.
Let $\theta$ be an involutive automorphism of $G$ and let $d\theta:{\mathfrak g}\longrightarrow{\mathfrak g}$ be the corresponding linear involution of ${\mathfrak g}=\Lie(G)$.
There is a direct sum decomposition ${\mathfrak g}={\mathfrak k}\oplus{\mathfrak p}$, where ${\mathfrak k}=\{ x\in{\mathfrak g}|d\theta(x)=x\}$ and ${\mathfrak p}=\{ x\in{\mathfrak g}|d\theta(x)=-x\}$.
Let $G^\theta=\{ g\in G|\theta(g)=g\}$ and let $K$ be the connected component of $G^\theta$ containing the identity element.
$K$ is reductive and normalises ${\mathfrak p}$, and ${\mathfrak k}=\Lie(K)$.

The idea of the representation $G^\theta\rightarrow\GL({\mathfrak p})$
as a `generalized version' of the adjoint representation goes back at
least as far as Cartan; but achieved a certain maturity in the well-known work \cite{kostrall}.
There Kostant and Rallis show that the action of $G^\theta$ on ${\mathfrak p}$ exhibits similar properties to the adjoint action of $G$ on ${\mathfrak g}$.
In the set-up of \cite{kostrall}, ${\mathfrak g}$ is a complex reductive Lie algebra, $G$ is the adjoint group of ${\mathfrak g}$ and $\theta$ is an involution of ${\mathfrak g}$ defined over a real form ${\mathfrak g}_{\mathbb R}$.
Many of the arguments in \cite{kostrall} use compactness properties and $\mathfrak{sl}(2)$-triples.
These arguments are not valid in positive characteristic.
On the other hand, Kostant-Rallis' results are generally assumed to be true over arbitrary (algebraically closed) fields of characteristic zero.

More recent work by Vust \cite{vust} and Richardson \cite{rich2} considers an analogous `symmetric space' decomposition in a reductive algebraic group $G$.
The object corresponding to ${\mathfrak p}$ is the closed set $P=\{ g\theta(g^{-1})\,|\,g\in G\}$: $G$ acts on $P$ by the {\it twisted} action $x*(g\theta(g^{-1}))=xg\theta(g^{-1})\theta(x^{-1})$.
If $x\in K$, this action is just ordinary conjugation.
(It was proved by Richardson that the twisted action induces a $G$-equivariant isomorphism of varieties $\sigma:G/G^\theta\rightarrow P$, where $G/G^\theta$ is the space of left cosets of $G$ modulo $G^\theta$.)

This paper will extend the analysis in the first two chapters of \cite{kostrall} to the case where $p$ is a good prime.
Our exposition proceeds along similar lines to \cite{kostrall}.
The main obstacles to be overcome are the construction of a $d\theta$-equivariant trace form on ${\mathfrak g}$ (Sect. \ref{sec:3}) and the replacement of the language of $\mathfrak{sl}(2)$-triples with that of associated cocharacters (Sect. \ref{sec:5}).
These adjustments allow us to generalise all of the relevant parts of
\cite{kostrall}.
In addition: in Sect. \ref{sec:5.5} and
Sect. \ref{sec:6.3} we give a new proof the number of irreducible
components of the variety ${\cal N}$ of nilpotent elements of
${\mathfrak p}$ (following Sekiguchi \cite{sek} in characteristic zero); we show in Sect. \ref{sec:6.4} that a conjecture of Richardson concerning the quotient morphism $\pi:P\rightarrow P\quot K$ is false; finally, we apply a theorem of Skryabin to describe the ring $k[{\mathfrak p}]^{K_i}$, where $K_i$ is the $i$-th Frobenius kernel of $K$.

A torus $A$ in $G$ is {\it $\theta$-split} if $\theta(a)=a^{-1}$ for all $a\in A$.
It was proved by Vust that the set of maximal $\theta$-split tori are $K$-conjugate.
Let ${\mathfrak a}$ be a toral algebra contained in ${\mathfrak p}$.
If ${\mathfrak a}$ is maximal such, then by abuse of terminology we say that ${\mathfrak a}$ is a {\it maximal torus} of ${\mathfrak p}$.

\begin{lemma}
Let ${\mathfrak a}$ be a maximal torus of ${\mathfrak p}$.
Then ${\mathfrak z}_{\mathfrak g}({\mathfrak a})\cap{\mathfrak p}={\mathfrak a}$, and there exists a unique maximal $\theta$-split torus $A$ of $G$ such that $\Lie(A)={\mathfrak a}$.
\end{lemma}

We reintroduce Kostant and Rallis' definition of a Cartan subspace, and check that it is valid in positive characteristic.
We provide a short proof of the following result from \cite{kostrall}.

\begin{theorem}
Any two Cartan subspaces of ${\mathfrak p}$ are conjugate by an element of $K$.
The Cartan subspaces of ${\mathfrak p}$ are just the maximal tori of ${\mathfrak p}$.
An element of ${\mathfrak p}$ is semisimple if and only if it is contained in a Cartan subspace.
\end{theorem}

The only assumption required for the above is (A) that $p$ is good for $G$.
We make the further assumptions from Sect. \ref{sec:3} onwards: (B) the derived subgroup of $G$ is simply-connected, and (C) there exists a non-degenerate $G$-equivariant symmetric bilinear form $\kappa:{\mathfrak g}\times{\mathfrak g}\longrightarrow k$.
The hypotheses (A),(B), and (C) are sometimes known as the standard hypotheses.

In order to make maximum use of the assumption (C), we would like the form $\kappa$ to be $d\theta$-equivariant.
This is straightforward in characteristic zero, but requires a more subtle argument if the characteristic is positive.
In order to construct the required $\kappa$, we develop a $\theta$-stable version of a reduction theorem of Gordon and Premet.
We then use this reduction theorem to prove our desired result.

\begin{lemma}
The trace form $\kappa$ in (C) may be chosen to be $d\theta$-equivariant.
\end{lemma}

The $d\theta$-equivariance of $\kappa$ allows us to proceed as in \cite{kostrall} in Sect. \ref{sec:4}.

\begin{lemma}\label{centdimintro}
Let $x\in{\mathfrak p}$.
Then $\dim{\mathfrak z}_{\mathfrak k}(x)-\dim{\mathfrak z}_{\mathfrak p}(x)=\dim{\mathfrak k}-\dim{\mathfrak p}$.
\end{lemma}

With Lemma \ref{centdimintro} we can define regularity: an element $x\in{\mathfrak p}$ is {\it regular} if $\dim{\mathfrak z}_{\mathfrak g}(x)\leq\dim{\mathfrak z}_{\mathfrak g}(y)$ for all $y\in{\mathfrak p}$.

\begin{lemma}
Let $x\in{\mathfrak p}$.
The following are equivalent:

(i) $x$ is regular,
(ii) $\dim{\mathfrak z}_{\mathfrak g}(x)=\dim{\mathfrak g}^A$,
(iii) $\dim{\mathfrak z}_{\mathfrak k}(x)=\dim{\mathfrak k}^A$,
(iv) $\dim{\mathfrak z}_{\mathfrak p}(x)=\dim A$.
\end{lemma}

Recall that, for a rational representation $\rho:H\longrightarrow\GL(V)$, an element $v\in V$ is {\it $H$-unstable} if $0$ is in the closure of $\rho(H)(v)$.

\begin{lemma}
Let $x\in{\mathfrak p}$.
Then $x$ is $K$-unstable if and only if $x$ is nilpotent.
\end{lemma}

It follows fairly quickly that:

\begin{lemma}
Let $x\in{\mathfrak p}$, and let $x=x_s+x_n$ be the Jordan-Chevalley decomposition of $x$.
The unique closed (resp. minimal) $K$-orbit in the closure of $\Ad K(x)$ is $\Ad K(x_s)$.
\end{lemma}

It is well-known from Mumford's Geometric Invariant Theory that the closed $K$-orbits in ${\mathfrak p}$ are in one-to-one correspondence with the $k$-rational points of the quotient ${\mathfrak p}\quot K = \Spec(k[{\mathfrak p}]^K)$.
We have a Chevalley Restriction Theorem for ${\mathfrak p}\quot K$.
The proof follows Richardson's proof of the corresponding result for the action of $K$ on $P=\{ g^{-1}\theta(g)\,|\,g\in G\}$.

\begin{theorem}\label{Chevintro}
Let $A$ be a maximal $\theta$-split torus of $G$, and let $W=N_G(A)/Z_G(A)$.
Let ${\mathfrak a}=\Lie(A)$.
Then the natural embedding $j:{\mathfrak a}\rightarrow{\mathfrak p}$ induces an isomorphism of affine varieties $j':{\mathfrak a}/W\rightarrow{\mathfrak p}\quot K$.
Hence $k[{\mathfrak p}]^K\cong k[{\mathfrak a}]^W$.
\end{theorem}

If ${\mathfrak g}$ is a complex reductive Lie algebra with adjoint group $G$, then by a well-known classical result $k[{\mathfrak g}]^G$ is a polynomial ring in $(\rank{\mathfrak g})$ indeterminates.
Here a straightforward application of Demazure's theorem on Weyl group invariants gives the analogous result:

\begin{lemma}
There are $r=\dim A$ algebraically independent homogeneous polynomials $f_1,f_2,\ldots,f_r$ such that $k[{\mathfrak a}]^W=k[f_1,f_2,\ldots,f_r]$.
Moreover,
$$\sum_{w\in W} t^{l(w)}=\prod_{i=1}^r{\frac{1-t^{\deg f_i}}{1-t}}$$
where $l$ is the length function on $W$ corresponding to a basis of simple roots in $\Phi_A$.
\end{lemma}

In Sect. \ref{sec:5} we consider in more detail the set of nilpotent elements of ${\mathfrak p}$, denoted ${\cal N}$.
In general ${\cal N}$ is not irreducible (and therefore not normal as 0 is in every irreducible component).
However, it is straightforward to prove (following \cite{kostrall}):

\begin{theorem}
Let ${\cal N}_1,{\cal N}_2,\ldots,{\cal N}_m$ be the irreducible components of ${\cal N}$.
The number of $K$-orbits in ${\cal N}$ is finite.
Each irreducible component ${\cal N}_i$ is normalized by $K$, contains a unique open $K$-orbit, and is of codimension $r=\rank A$ in ${\mathfrak p}$ (where $A$ is a maximal $\theta$-split torus).
An element of ${\cal N}_i$ is in the open $K$-orbit if and only if it is regular.
\end{theorem}

Let $K^*=\{ g\in G\,|\,g^{-1}\theta(g)\in Z(G)\}$.
In \cite{kostrall} it was proved that the irreducible components of ${\cal N}$ are permuted transitively by $K^*$.
For the proof, Kostant and Rallis showed that any regular nilpotent element of ${\mathfrak p}$ can be embedded as the nilpositive element in a principal normal $\mathfrak{sl}(2)$-triple, and that any two principal normal $\mathfrak{sl}(2)$-triples are conjugate by an element of $K^*$.
(A principal normal $\mathfrak{sl}(2)$-triple $\{ h,e,f\}$ is one such that $e,f\in{\mathfrak p}$ are regular and $h\in{\mathfrak k}$.)
Clearly, this argument cannot be applied if the characteristic is small.
To prove it in our case we replace the language of $\mathfrak{sl}(2)$-triples with that of (Pommerening's) associated cocharacters.
A reinterpretation of Kawanaka's theorem \cite{kawanaka} on nilpotent orbits in graded semisimple Lie algebras gives the following:

\begin{corollary}\label{cocharintro}
Let $e\in{\cal N}$.
Then there exists a cocharacter $\lambda:k^\times\longrightarrow K$ which is associated to $e$.
Any two such cocharacters are conjugate by an element of $Z_K(e)^\circ$.
\end{corollary}

The key step in proving that the set of regular nilpotent elements is a single $K^*$-conjugacy class is the following lemma.
For the proof, we reduce by a number of tricks to the case where $G$ is almost simple, $e$ is semiregular in ${\mathfrak g}$, and $\theta=\Ad\lambda(t_0)$, where $\lambda$ is an associated cocharacter for $e$ and $t_0=\sqrt{-1}$.
It is then fairly straightforward to prove the Lemma case-by-case (see Sect. \ref{sec:5.4}).

\begin{lemma}
Let $e\in{\cal N}$ and let $\lambda:k^\times\longrightarrow K$ be associated to $e$.
There exists $g\in G$ such that $(\Int g)\circ\lambda$ is $\theta$-split.
Equivalently $\Int n\circ\lambda=-\lambda$, where $n=g^{-1}\theta(g)$.
\end{lemma}

As a consequence, we have:

\begin{corollary}\label{omegaintro}
Let $A$ be a maximal $\theta$-split torus of $G$ and let $\Pi$ be a basis for $\Phi_A=\Phi(G,A)$.
Then $e$ is regular in ${\mathfrak p}$ if and only if $\lambda$ is $G$-conjugate to the cocharacter $\omega:k^\times\longrightarrow A\cap G^{(1)}$ satisfying $\langle\alpha,\omega\rangle=2$ for all $\alpha\in\Pi$.
\end{corollary}

The above Corollary shows that any regular nilpotent element of ${\mathfrak p}$ is even (see Rk. \ref{eiseven} for details).
It is now a fairly straightforward task to deduce that:

\begin{theorem}\label{denseorbitintro}
The regular elements ${\cal N}_{reg}\subset{\cal N}$ are a single $K^*$-orbit.
Hence $K^*$ permutes the components of ${\cal N}$ transitively.
\end{theorem}

It is easy to give examples such that ${\cal N}$ is not irreducible.
In \cite{sek}, Sekiguchi classified (over $k={\mathbb C}$) the involutions for which ${\cal N}$ is not irreducible.
Our analysis of associated cocharacters, combined with the
classification of involutions (see for example \cite{springer}), simplifies the task of extending Sekiguchi's results to positive characteristic.
We begin with the following observation.

\begin{theorem}\label{gthetaorbsintro}
Let $e,\lambda$ be as above and let $C=Z_G(\lambda)\cap Z_G(e)$ (the reductive part of $Z_G(e)$, see \cite[Thm. 2.3]{premnil}).
Let $Z=Z(G)$ and $P=\{ g^{-1}\theta(g)\,|\, g\in G\}$.
Denote by $\tau:G\longrightarrow P$ the morphism $g\mapsto g^{-1}\theta(g)$, and by $\Gamma$ the set of $G^\theta$-orbits in ${\cal N}_{reg}$.

(a) The map from $K^*$ to $\Gamma$ given by $g\mapsto gG^\theta\cdot e$ is surjective and induces a one-to-one correspondence $K^*/G^\theta C\longrightarrow\Gamma$.

(b) The morphism $\tau$ induces an isomorphism $K^*/G^\theta C\longrightarrow (Z\cap A)/{\tau(C)}$.
Since $Z\subseteq C$, there is a surjective map $(Z\cap A)/{\tau(Z)}\longrightarrow (Z\cap A)/{\tau(C)}$.

(c) The embedding $F^*\hookrightarrow K^*$ induces a surjective map $F^*/{F(Z\cap A)}\rightarrow\Gamma$.

(d) The map $F^*\rightarrow Z\cap A$, $a\mapsto a^2$ induces an isomorphism $F^*/{F(Z\cap A)}$ $\longrightarrow Z\cap A/{(Z\cap A)^2}$.
\end{theorem}

Thm. \ref{gthetaorbsintro} holds for an arbitrary reductive group $G$.
If $G$ is semisimple and simply-connected, then $G^\theta=K$ by a result of Steinberg, hence the $G^\theta$-orbits in ${\cal N}_{reg}$ are in one-to-one correspondence with the irreducible components of ${\cal N}$.
We can use this observation together with Thm. \ref{gthetaorbsintro} to describe the number of irreducible components of ${\cal N}$ for any involution of an almost simple group.
An involution $\theta$ of $G$ is of {\it maximal rank} if the maximal $\theta$-split torus $A$ is a maximal torus of $G$.
If $G$ is almost simple and $\theta$ is of maximal rank, then $(Z\cap A)/{\tau(C)}=Z/Z^2$.
For example, Thm. \ref{gthetaorbsintro} implies immediately that the variety of $n\times n$ symmetric nilpotent matrices has two irreducible components if $n$ is even, and is irreducible if $n$ is odd.
(See Sect. \ref{sec:5.5} and Sect. \ref{sec:6.3} for further details.)

In Sect. \ref{sec:6} we generalise Kostant-Rallis' construction of a reductive subalgebra ${\mathfrak g}^*\subset{\mathfrak g}$ such that ${\mathfrak a}$ is a Cartan subalgebra of ${\mathfrak g}^*$.

\begin{theorem}\label{subalgintro}
Let $\omega$ be as in Cor. \ref{omegaintro} and let $E\in{\mathfrak g}(2;\omega)$ be such that $[{\mathfrak g}^\omega,E]={\mathfrak g}(2;\omega)$.
Let ${\mathfrak g}^*$ be the Lie subalgebra of ${\mathfrak g}$ generated by $E,d\theta(E)$ and ${\mathfrak a}$.

(a) ${\mathfrak a}$ is a Cartan subalgebra of ${\mathfrak g}^*$.
There exists a reductive group $G^*$ satisfying the standard hypotheses (A)-(C), such that $\Lie(G^*)={\mathfrak g}^*$.

(b) There is an involutive automorphism $\theta^*$ of $G^*$ such that $d\theta^*=d\theta|_{{\mathfrak g}^*}$.
\end{theorem}

In \cite{kostrall}, it was proved that each fibre of the quotient morphism $\pi:{\mathfrak p}\longrightarrow{\mathfrak p}\quot K$ has a dense open (regular) $K^*$-orbit.
Let $K^*$ act on $P$ by conjugation (this is valid by \cite[8.2]{rich2}).
In \cite{rich2}, Richardson conjectured that there is a dense open $K^*$-orbit on each fibre of the quotient morphism $\pi_P:P\longrightarrow P\quot K=P\quot {K^*}\cong A/ W_A$ (see \cite[8.3-4]{rich2}).

\begin{theorem}
(a) There is a dense open $K^*$-orbit in each fibre of $\pi_{\mathfrak p}$.

(b) The corresponding statement for $\pi_P$ is false.
\end{theorem}

We draw some further conclusions from Thm. \ref{subalgintro}.
Let ${\mathfrak k}^*={\mathfrak g}^*\cap{\mathfrak k}$, ${\mathfrak p}^*={\mathfrak g}^*\cap{\mathfrak p}$.
By Thm. \ref{Chevintro} and the Chevalley Restriction Theorem, $k[{\mathfrak p}]^K\cong k[{\mathfrak a}]^{W_A}\cong k[{\mathfrak g}^*]^{G^*}$.

\begin{lemma}
If two elements of ${\mathfrak g}^*$ are $G^*$-conjugate, then they are $G$-conjugate.
\end{lemma}

This allows us to establish the following equivalence:

\begin{lemma}
Let $x\in{\mathfrak p}^*$.
The following are equivalent:
(i) $x$ is a regular element of ${\mathfrak p}$,
(ii) $x$ is a regular element of ${\mathfrak g}^*$,
(iii) ${\mathfrak z}_{{\mathfrak k}^*}(x)=0$,
(iv) $\dim{\mathfrak z}_{{\mathfrak p}^*}(x)=r=\dim{\mathfrak a}$.
\end{lemma}

Let $e\in{\mathfrak p}^*$ be a regular nilpotent element.
By Cor. \ref{cocharintro} there exists a cocharacter $\lambda:k^\times\longrightarrow (G^*)^{\theta^*}$ which is associated to $e$.
Hence we can choose a $\lambda$-graded subspace ${\mathfrak v}$ of ${\mathfrak p}^*$ such that $[e,{\mathfrak g}^*]\oplus{\mathfrak v}={\mathfrak g}^*$.
Then we also have $[e,{\mathfrak k}]\oplus{\mathfrak v}={\mathfrak p}$.
It is known (\cite{veldkamp,premtang}) that every element of $e+{\mathfrak v}$ is regular in ${\mathfrak g}^*$, that the embedding $e+{\mathfrak v}\hookrightarrow{\mathfrak g}^*$ induces an isomorphism $e+{\mathfrak v}\longrightarrow{\mathfrak g}^*\quot G^*$, and that each regular orbit in ${\mathfrak g}^*$ intersects $e+{\mathfrak v}$ in precisely one point.

\begin{lemma}
Let $j$ be the composite of the isomorphisms $k[{\mathfrak p}]^K\rightarrow k[{\mathfrak a}]^{W}\rightarrow k[{\mathfrak g}^*]^{G^*}$ and let $f\in k[{\mathfrak p}]^K,g\in k[{\mathfrak g}^*]^{G^*}$.
Then $j(f)=g$ if and only if $f|_{e+{\mathfrak v}}=g|_{e+{\mathfrak v}}$.
Hence the embedding $e+{\mathfrak v}\hookrightarrow {\mathfrak p}$ induces an isomorphism $e+{\mathfrak v}\longrightarrow{\mathfrak p}\quot K$, and each regular $K^*$-orbit in ${\mathfrak p}$ intersects $e+{\mathfrak v}$ in exactly one point.
\end{lemma}

In particular we have:

\begin{corollary}
Let $k[{\mathfrak p}]^K=k[u_1,u_2,\ldots,u_r]$.
The differentials $(du_i)_x$ are linearly independent for any regular $x\in{\mathfrak p}$.
\end{corollary}

The above observation allows us to apply Skryabin's theorem on infinitesimal invariants to show that:

\begin{theorem}
Let $k[{\mathfrak p}]^{(p^i)}$ denote the ring of all $p^i$-th powers of elements of $k[{\mathfrak p}]$ and let $K_i$ denote the $i$-th Frobenius kernel of $K$.

(a) $k[{\mathfrak p}]^{K_i}=k[{\mathfrak p}]^{(p^i)}[u_1,u_2,\ldots,u_r]$, and is free of rank $p^{ir}$ over $k[{\mathfrak p}]^{(p^i)}$.

(b) $k[{\mathfrak p}]^{K_i}$ is a locally complete intersection.
\end{theorem}

{\it Notation.}
The connected component of an algebraic group $G$ (containing the identity element) will be denoted $G^\circ$.
If $\theta$ is an automorphism of $G$, then we denote by $G^\theta$ the isotropy subgroup $\{ g\in G\,|\,\theta(g)=g\}$.
We use similar notation for the fixed points of an algebra or Lie algebra with respect to an automorphism or group of automorphisms.
If $x\in G$, then $Z_G(x)$ (resp. ${\mathfrak g}^x$) will denote the centralizer of $x$ in $G$ (resp. in ${\mathfrak g}$).
Similar notation will be used, where appropriate, for the centralizers in $K,{\mathfrak k},{\mathfrak p}$, etc.
We write $x=x_s x_u$ (resp. $x=x_s+x_n$) for the Jordan-Chevalley decomposition of $x\in G$ (resp. $x\in{\mathfrak g}$), where $x_s$ is the semisimple part and $x_u$ is the unipotent part (resp. $x_n$ is the nilpotent part) of $x$.
Throughout the paper we write $g\cdot x$ (resp. $g\cdot\lambda$) for $\Ad g(x)$ (resp. $\Ad g\circ\lambda$), where $g\in G$ and $x\in{\mathfrak g}$ (resp. $\lambda$ is a cocharacter in $G$).\linebreak[2]

{\it Acknowledgement.}

I would like to thank Alexander Premet for his advice and
encouragement.
I also thank Dmitri Panyushev for directing me towards
the paper of Sekiguchi and alerting me to some mistakes in the
calculation of the number of irreducible components of ${\cal N}$.
This paper was written in part with the support of an Engineering and
Physical Sciences Research Council PhD studentship.

\section{Preliminaries}
\label{sec:1}

Let $G$ be a reductive algebraic group over the algebraically closed field $k$ of characteristic not equal to $2$.
We assume throughout that $\Char k=p$ is good for $G$.
(Let $\Delta$ be a basis for the root system $\Phi$ of $G$, let $\hat\alpha$ be the longest element of $\Phi$ relative to $\Delta$, and let $\hat\alpha=\sum_{\beta\in\Delta}m_\beta \beta$.
Then $p$ is good for $G$ if and only if $p>m_\beta$ for all $\beta\in\Delta$.)
Let $\theta:G\longrightarrow G$ be an involutive automorphism and let $K$ denote the connected component of the isotropy subgroup $G^\theta$.
Let ${\mathfrak g}=\Lie(G)$.
Then ${\mathfrak g}={\mathfrak k}\oplus{\mathfrak p}$, where ${\mathfrak k}=\{ x\in{\mathfrak g}|\;d\theta(x)=x\}$, ${\mathfrak p}=\{ x\in{\mathfrak g}|\;d\theta(x)=-x\}$.
Clearly $[{\mathfrak k},{\mathfrak k}]\subseteq{\mathfrak k}$, $[{\mathfrak k},{\mathfrak p}]\subseteq{\mathfrak p}$, and $[{\mathfrak p},{\mathfrak p}]\subseteq{\mathfrak k}$.
Hence we have a ${\mathbb Z}/2{\mathbb Z}$-grading of ${\mathfrak g}$.
By \cite[8.1]{steinberg}, $K$ is reductive.
Moreover, $\Lie(K)={\mathfrak k}$ by \cite[\S 9.1]{bor}.
The following result is due to Steinberg \cite[7.5]{steinberg}:

{\it - There exists a Borel subgroup $B$ of $G$ and a maximal torus $T$ contained in $B$ such that $\theta(B)=B,\,\theta(T)=T$.}

Following Springer \cite{springer} we call such a pair $(B,T)$ a {\it fundamental pair}.
Let $(B,T)$ be a fundamental pair and let $\Delta$ be the basis of the root system $\Phi=\Phi(G,T)$ corresponding to $B$.
Let $\{ h_\alpha,e_\beta\,:\,\alpha\in\Delta,\beta\in\Phi\}$ be a Chevalley basis for ${\mathfrak g}'=\Lie(G^{(1)})$.
There exist constants $\{c(\alpha)\in k^\times:\alpha\in\Phi\}$ and an automorphism $\gamma$ of $\Phi$ with $\gamma(\Delta)=\Delta$ such that $d\theta(e_\alpha)=c(\alpha)e_{\gamma(\alpha)}$ for each $\alpha\in\Phi$.
It is easy to see that:

{\it - $c(\alpha)c(\gamma(\alpha))=1$,

 - If $\gamma(\alpha)\neq\alpha$, then either $\gamma(\alpha)$ and $\alpha$ are orthogonal, or they generate a root system of type $A_2$,

 - $c(\alpha)c(-\alpha)=1$,

 - $\theta(h_\alpha)=h_{\gamma(\alpha)}$ for all $\alpha\in\Delta$.}

If $G$ is semisimple, then the data $\gamma$ and $\{ c(\alpha),\alpha\in\Delta\}$ fully determine $d\theta$.
In the general reductive case, we need a little more preparation.

Recall that ${\mathfrak g}$ is a {\it restricted} Lie algebra.
Thus there is a canonical $p$-operation on ${\mathfrak g}$, denoted $x\mapsto x^{[p]}$.
If $G$ is a closed subgroup of some $\GL(V)$, then ${\mathfrak g}$ is a subalgebra of $\mathfrak{gl}(V)$ and the $p$-operation is just the restriction to ${\mathfrak g}$ of the $p$-th power map of matrices.
An element $t\in{\mathfrak g}$ is a {\it toral element} if $t^{[p]}=t$.
A subalgebra of ${\mathfrak g}$ is a {\it toral algebra} if it is commutative and has a basis of toral elements.
If $T$ is a torus in $G$ then $\Lie(T)$ is a toral algebra in ${\mathfrak g}$.
For a toral algebra ${\mathfrak s}\subseteq{\mathfrak g}$, we denote by ${\mathfrak s}^{tor}$ the set of all toral elements in ${\mathfrak s}$: ${\mathfrak s}^{tor}$ is a vector space over the prime subfield ${\mathbb F}_p$ of $k$, and ${\mathfrak s}\cong{\mathfrak s}^{tor}\otimes_{{\mathbb F}_p} k$.

\begin{lemma} \label{redcase}
Let $\theta$ be an automorphism of $G$ of order $m$, $p\nmid m$, let $T$ be a $\theta$-stable torus in $G$ and let ${\mathfrak t}=\Lie(T)$, ${\mathfrak t}'=\Lie(T\cap G^{(1)})$.
There exists a $\theta$-stable toral algebra ${\mathfrak s}$ such that ${\mathfrak t}={\mathfrak t}'\oplus{\mathfrak s}$, and hence ${\mathfrak g}={\mathfrak g}'\oplus{\mathfrak s}$ (vector space direct sum).

If $m|(p-1)$, then we can choose a toral basis for ${\mathfrak s}$ consisting of eigenvectors for $d\theta$.
\end{lemma}

\begin{proof}
As $d\theta$ is a restricted Lie algebra automorphism, the sets ${\mathfrak t}^{tor}$ and $({\mathfrak t}')^{tor}$ are $d\theta$-stable.
Therefore by Maschke's Theorem there is a $d\theta$-stable ${\mathbb F}_p$-vector space ${\mathfrak s}^{tor}$ such that ${\mathfrak t}^{tor}=({\mathfrak t}')^{tor}\oplus{\mathfrak s}^{tor}$.
Let ${\mathfrak s}$ be the toral algebra generated by ${\mathfrak s}^{tor}$.
Then ${\mathfrak t}={\mathfrak t}'\oplus{\mathfrak s}$.

To prove the second assertion, we consider the action of $d\theta$ on ${\mathfrak s}^{tor}$.
As $\theta$ has order $m$, the minimal polynomial $m(t)$ of $d\theta|_{{\mathfrak s}^{tor}}$ divides $(t^m-1)$.
But if $m$ divides $(p-1)$ then there is a primitive $m$-th root of unity in ${\mathbb F}_p$, hence $m(t)$ splits over ${\mathbb F}_p$ as a product of distinct linear factors.
In other words $d\theta|_{{\mathfrak s}^{tor}}$ is diagonalizable.
Choose a basis for ${\mathfrak s}^{tor}$ consisting of eigenvectors for $d\theta$.
This completes the proof.
%\qed
\end{proof}

Let us return now to the case where $\theta$ is an involution.
It may be illustrative at this point to give explicit bases for ${\mathfrak k}$ and ${\mathfrak p}$.

For ${\mathfrak k}$: $\left\{
\begin{array}{ll} h_{\alpha_i} & \alpha_i\in\Delta,\gamma(\alpha_i)=\alpha_i, \\
h_{\alpha_i}+h_{\gamma(\alpha_i)} & \alpha_i\in\Delta,\gamma(\alpha_i)\neq\alpha_i, \\
e_\alpha & \mbox{$\alpha\in\Phi,\gamma(\alpha)=\alpha$ and $c(\alpha)=1$}, \\
e_\alpha+c(\alpha)e_{\gamma(\alpha)} & \alpha\in\Phi,\gamma(\alpha)\neq\alpha, \\
t_i & 1\leq i\leq l.
\end{array}
\right.
$

For ${\mathfrak p}$: $\left\{
\begin{array}{ll} h_{\alpha_i}-h_{\gamma(\alpha_i)} & \alpha_i\in\Delta,\gamma(\alpha_i)\neq\alpha_i, \\
e_\alpha & \mbox{$\alpha\in\Phi,\gamma(\alpha)=\alpha$ and $c(\alpha)=-1$}, \\
e_\alpha-c(\alpha)e_{\gamma(\alpha)} & \alpha\in\Phi,\gamma(\alpha)\neq\alpha, \\
t_j' & 1\leq j\leq h.
\end{array}
\right.
$

The elements $t_i,t_j'$ are toral elements spanning the toral algebra ${\mathfrak s}$ of Lemma \ref{redcase}.
With this description we can prove the following useful lemma:

\begin{lemma}\label{nontriv}
The following are equivalent:

(i) ${\mathfrak p}$ is a toral algebra contained in ${\mathfrak z}({\mathfrak g})$,

(ii) There are no non-central semisimple elements in ${\mathfrak p}$,

(iii) There are no non-zero nilpotent elements in ${\mathfrak p}$,

(iv) $\theta|_{G^{(1)}}$ is trivial.
\end{lemma}

\begin{proof}
Clearly (i) $\Rightarrow$ (ii) and (i) $\Rightarrow$ (iii).
Suppose (iv) holds.
Then, by the above remarks ${\mathfrak p}$ is a toral algebra contained in ${\mathfrak t}$.
Let $t\in{\mathfrak p}$ and let $\alpha\in\Phi$, hence $e_\alpha\in{\mathfrak g}'\subseteq{\mathfrak k}$.
Then $[t,e_\alpha]=d\alpha(t)e_\alpha\in{\mathfrak p}\Rightarrow d\alpha(t)=0$.
Thus $t\in{\mathfrak z}({\mathfrak g})$, and (i) holds.

To complete the proof we will show that (ii)$\Rightarrow$(iv) and (iii)$\Rightarrow$(iv).
Keep the notation from above, and suppose that $\theta|_{G^{(1)}}$ is non-trivial.
We will show that (ii) cannot hold.
Assume first of all that $\theta|_{G^{(1)}}$ is inner.
There is some $\alpha\in\Delta$ such that $e_\alpha\in{\mathfrak p}$.
Moreover $e_{-\alpha}\in{\mathfrak p}$ also, since $c(\alpha)c(-\alpha)=1$.
Hence $s=e_{\alpha}+e_{-\alpha}$ is a semisimple element of ${\mathfrak p}$.
But $s$ is not in ${\mathfrak h}$ and therefore $s\notin{\mathfrak z}$ (see \cite[2.3]{me}).
Assume therefore that $\gamma$ is non-trivial.
Then $\alpha\neq\gamma(\alpha)$ for some $\alpha\in\Delta$.
Hence $h=h_\alpha - h_{\gamma(\alpha)}\in{\mathfrak p}$.
If (ii) holds then $h\in{\mathfrak z}$, hence $\charac k=3$ and $\alpha,\gamma(\alpha)$ generate a subsystem of $\Phi$ of type $A_2$.
Thus $[e_\alpha,e_{\gamma(\alpha)}]=Ne_{\alpha+\gamma(\alpha)}\in{\mathfrak p}$, $N\neq 0$.
Therefore $e_{\alpha+\gamma(\alpha)}\in{\mathfrak p}$, and by the same argument $e_{-(\alpha+\gamma(\alpha))}\in{\mathfrak p}$.
Let $s= e_{\alpha+\gamma(\alpha)} + e_{-(\alpha+\gamma(\alpha))}$.
Then $s$ is a semisimple element of ${\mathfrak p}$ not in ${\mathfrak z}({\mathfrak g})$.

We have shown that (ii) $\Rightarrow$ (iv).
It remains to prove that if $\theta|_{G^{(1)}}$ is non-trivial then there is a non-zero nilpotent element of ${\mathfrak p}$.
If $\gamma$ is non-trivial, then we choose $\alpha$ with $\gamma(\alpha)\neq\alpha$ and set $n= e_\alpha - d\theta(e_\alpha)=e_\alpha-c(\alpha)e_{\gamma(\alpha)}$.
If $\theta|_{G^{(1)}}$ is inner, then we can choose $\alpha\in\Phi$ with $e_\alpha\in{\mathfrak p}$.
This completes the proof.
%\qed
\end{proof}

We will require the following observation of Steinberg:

\begin{lemma}\label{sccover}
Let $G$ be a semisimple group and let $\theta$ be an automorphism of $G$.
Let $\pi:G_{sc}\rightarrow G$ be the universal covering of $G$.
Then there exists a unique automorphism $\theta_{sc}$ of $G_{sc}$ such that the following diagram is commutative:
\begin{diagram}
G_{sc} & \rTo^{\theta_{sc}} & G_{sc} \\
\dTo^\pi & & \dTo^\pi \\
G & \rTo^{\theta} & G \\
\end{diagram}

If $\theta$ is an involution, then so is $\theta_{sc}$.
\end{lemma}

\begin{proof}
The first statement follows from \cite[9.16]{steinberg}.
But now by uniqueness, if $\theta$ is of order 2 then so is $\theta_{sc}$.
%\qed
\end{proof}

Finally, we make the following observation for later reference.

\begin{lemma}\label{GLautos}
Let $G=\GL(n,k),{\mathfrak g}=\Lie(G),{\mathfrak g}'=\Lie(G^{(1)})$.
We denote by $\Aut G$ (resp. $\Aut{\mathfrak g}$) the (abstract) group of algebraic automorphisms of $G$ (resp. restricted Lie algebra automorphisms of ${\mathfrak g}$).

(i) $\Aut G$ contains $\Int G$, the inner automorphisms, as a normal subgroup of index 2.
For $n\geq 3$ (resp. $n=2$) let $\phi:G\longrightarrow G$ be the involution given by $g\mapsto {^t}g^{-1}$ (resp. $g\mapsto g/{(\det g)}$) and let $C$ be the subgroup of $\Aut G$ generated by $\phi$.
Then $\Aut G$ is the semidirect product of $\Int G$ by $C$ (resp. the direct product of $\Int G$ and $C$).

(ii) The natural map $\Aut G\rightarrow \Aut(G^{(1)})$ is bijective if $n\geq 3$, and surjective with kernel $C$ for $n=2$.

(iii) For any $\theta\in\Aut G$, the differential $d\theta$ is a restricted Lie algebra automorphism of $G$.
The map $d:\Aut G\longrightarrow\Aut{\mathfrak g}$ is injective and $d:\Aut G^{(1)}\longrightarrow\Aut{\mathfrak g}'$ is bijective for all $n$ and $p$.

(iv) If $p\nmid n$ then $\Aut{\mathfrak g}\cong\Aut{\mathfrak g}'\times{\mathbb F}_p^\times$.
If $p\, |\, n$ then $\Aut{\mathfrak g}\cong\Aut{\mathfrak g}'\times B$, where $B$ is the cyclic group of order $p$ generated by the automorphism $x\mapsto x+(\tr x)I$ and $I$ is the identity matrix.

(v) If $2\neq p\,|\,n$ then for any involution $\eta$ of the restricted Lie algebra ${\mathfrak g}'$ there is a unique involutive automorphism $\theta$ of $G$ (resp. $\psi$ of ${\mathfrak g}$) such that $d\theta|_{{\mathfrak g}'}=\eta$ (resp. $\psi|_{{\mathfrak g}'}=\eta$).
\end{lemma}

\begin{proof}
If $n=2$, then all automorphisms of $G^{(1)}$ are inner.
Otherwise, $\Aut G^{(1)}$ is generated by $\Int G^{(1)}$ together with the outer automorphism $g\mapsto {^t}g^{-1}$ (\cite[\S 14.9]{bor}).
Hence the restriction map $\Aut G\rightarrow\Aut G^{(1)}$ is surjective for any $n$.
Suppose $\theta\in\Aut G$ is such that $\theta(g)=g\;\forall g\in G^{(1)}$.
Then $\theta$ is trivial unless $\theta(z)=z^{-1}$ for all $z\in Z(G)$.
This possibility clearly only occurs if $n=2$ and $\theta:g\mapsto g/{(\det g)}$.
Hence we have proved (i) and (ii).

The automorphism group of the abstract Lie algebra ${\mathfrak g}'$ is given in \cite{hog}.
We can see easily from the tables in \cite{hog} that $d:\Aut G^{(1)}\longrightarrow\Aut{\mathfrak g}'$ is bijective (and that any automorphism of the abstract Lie algebra ${\mathfrak g}'$ is a restricted Lie algebra automorphism) unless $n=p=2$.
We deal with this case as follows:
Let $\{ h,e,f\}$ be the standard basis for ${\mathfrak g}'$.
Then $h$ is the identity matrix, and in fact is the only non-zero toral element of ${\mathfrak g}'$.
Hence any $\theta\in\Aut{\mathfrak g}'$ satisfies $\theta(h)=h$.
Suppose $\theta(e)=x$.
Then, since any two non-zero nilpotent elements of ${\mathfrak g}'$ are conjugate, there exists $g\in G^{(1)}$ such that $\Ad g(e)=x$.
But there is a unique nilpotent element $y\in {\mathfrak g}'$ such that $[x,y]=h$.
Hence $\Ad g(f)=y=\theta(f)$.
It follows that $\theta=\Ad g$.
Thus differentiation $d:\Aut G^{(1)}\longrightarrow\Aut{\mathfrak g}'$ is surjective.
Injectivity follows from the fact that $\ker\Ad=Z(G)$.

We have shown that $d:\Aut G^{(1)}\longrightarrow\Aut{\mathfrak g}'$ is bijective for all $n$ and $p$.
Therefore $d:\Aut G\longrightarrow\Aut{\mathfrak g}$ is injective for all $n\geq 3$.
Injectivity for $n=2$ will follow from (iv), since $d\phi:x\mapsto x-(\tr x)I$.
Suppose first of all that $p\nmid n$.
Then ${\mathfrak g}={\mathfrak z}({\mathfrak g})\oplus{\mathfrak g}'$, hence $\Aut{\mathfrak g}\cong\Aut{\mathfrak g}'\times\Aut{\mathfrak z}$.
The toral algebra ${\mathfrak z}$ is generated by the identity matrix.
Hence $\Aut{\mathfrak z}$ consists of the maps $\lambda I\longrightarrow m\lambda I$ with $m\in{\mathbb F}_p^\times$.
Thus $\Aut{\mathfrak g}\cong\Aut{\mathfrak g}'\times{\mathbb F}_p^\times$.
Assume therefore that $p\, |\, n$.
As $\Aut G\longrightarrow\Aut G^{(1)}$ is surjective and $\Aut G^{(1)}\cong\Aut{\mathfrak g}'$, any automorphism of ${\mathfrak g}'$ can be extended to an automorphism of ${\mathfrak g}$.
Therefore $\Aut{\mathfrak g}\longrightarrow\Aut{\mathfrak g}'$ is surjective.
Let $\phi\in\Aut{\mathfrak g}$ be such that $\phi(x)=x\;\forall x\in{\mathfrak g}'$.
Let $e_{ij}$ be the matrix with 1 in the $(i,j)$-th position and 0 elsewhere.
By considering the values $d\alpha(d\theta(e_{11}))$ for $\alpha\in\Phi$, we see that $d\theta(e_{11})=e_{11}+\lambda I$ for some $\lambda\in k$.
Moreover $e_{11}^{[p]}=e_{11}$, hence $\lambda\in{\mathbb F}_p$.
It follows that $\theta$ must be of the form $\theta_\lambda:x\mapsto x+\lambda(\tr x)I$ for some $\lambda\in{\mathbb F}_p$.
Moreover $\theta_\lambda$ is a valid automorphism of ${\mathfrak g}$ for each $\lambda\in{\mathbb F}_p$.
The description of $\Aut{\mathfrak g}$ follows.

To prove (v), suppose $2\neq p\, |\, n$.
Then $\Aut G\longrightarrow\Aut{\mathfrak g}'$ is bijective, hence for each involution $\eta$ of ${\mathfrak g}'$ there is a unique automorphism $\theta$ of $G$, necessarily involutive, such that $d\theta|_{{\mathfrak g}'}=\eta$.
Moreover, $\Aut{\mathfrak g}\cong \Aut{\mathfrak g}'\times B$, where $B$ is a cyclic group of order $p$.
Hence there is a unique element $\psi\in\Aut{\mathfrak g}$ of order 2 such that $\psi|_{{\mathfrak g}'}=\eta$.
%\qed
\end{proof}

\section{Cartan Subspaces}
\label{sec:2}

\subsection{Maximal Toral Algebras}
\label{sec:2.1}

In \cite{kostrall}, Kostant and Rallis defined Cartan subspaces of ${\mathfrak p}$ and showed that any two Cartan subspaces are $K$-conjugate.
In this section we will show that this extends to positive characteristic.
We follow \cite{kostrall}, although Lemma \ref{sep} and Cor. \ref{cart} are new.

We begin with two easy lemmas.

\begin{lemma}
Let $x\in{\mathfrak g}$, and denote the Jordan-Chevalley decomposition of $x$ by $x_s+x_n$.
Then $x\in{\mathfrak k}$ (resp. ${\mathfrak p}$) if and only if $x_s,x_n\in{\mathfrak k}$ (resp. ${\mathfrak p}$).
\end{lemma}

\begin{proof}
Any automorphism of ${\mathfrak g}$ maps semisimple (resp. nilpotent) elements to semisimple (resp. nilpotent) elements.
Thus $\theta(x)=\theta(x_s)+\theta(x_n)$ is the Jordan-Chevalley decomposition of $\theta(x)$ for any $x\in{\mathfrak g}$.
Hence $\theta(x)=\lambda x$ if and only if $\theta(x_s)=\lambda x_s,\theta(x_n)=\lambda x_n$.
%\qed
\end{proof}

The following lemma is in \cite{rich2}.
For completeness, we reproduce a proof here.

\begin{lemma}\label{stabletori}
Let $T$ be a $\theta$-stable torus of $G$.
Let $T_+=(T\cap K)^\circ$ and $T_-=\{t\in T|\theta(t)=t^{-1}\}^\circ$.
Then $T=T_+\cdot T_-$ and the intersection is finite.
Let ${\mathfrak t}=\Lie(T)$.
Then ${\mathfrak t}\cap{\mathfrak k}=\Lie(T_+)$ and ${\mathfrak t}\cap{\mathfrak p}=\Lie(T_-)$.
\end{lemma}

\begin{proof}
Clearly $T_+$ and $T_-$ are subtori of $T$.
We consider the surjective morphism $p_+:T\longrightarrow T_+$, $t\mapsto t\theta(t)$.
Evidently $T_-$ is the connected component of $\ker p_+$ containing the identity element.
Hence $\dim T_-+\dim T_+=\dim T$.
Moreover $T_+\cap T_-$ is clearly finite.
Thus $T_+\cdot T_-=T$.
Clearly $\Lie(T_+)\subseteq{\mathfrak k}$ and $\Lie(T_-)\subseteq{\mathfrak p}$.
Therefore ${\mathfrak t}\supseteq\Lie(T_+)\oplus\Lie(T_-)$.
By equality of dimensions ${\mathfrak t}=\Lie(T_+)\oplus\Lie(T_-)$, from which the second part of the lemma follows immediately.
%\qed
\end{proof}

We call a toral algebra ${\mathfrak a}$ a {\it maximal torus} of ${\mathfrak p}$ if it is maximal in the collection of toral algebras contained in ${\mathfrak p}$.

\begin{lemma}\label{torcent}
Let ${\mathfrak a}$ be a maximal torus of ${\mathfrak p}$.
Then ${\mathfrak z}_{\mathfrak p}({\mathfrak a})={\mathfrak a}$.
\end{lemma}

\begin{proof}
Let $L=Z_G({\mathfrak a})$.
Then $L$ is a $\theta$-stable Levi subgroup of $G$, hence $p$ is good for $G$.
Moreover ${\mathfrak l}=\Lie(L)={\mathfrak z}_{\mathfrak g}({\mathfrak a})={\mathfrak z}_{\mathfrak k}({\mathfrak a})\oplus{\mathfrak z}_{\mathfrak p}({\mathfrak a})$ by \cite[\S 9.1]{bor}.
Since ${\mathfrak a}$ is maximal all semisimple elements of ${\mathfrak l}\cap{\mathfrak p}$ are in ${\mathfrak a}$.
Applying Lemma \ref{nontriv}, we see that ${\mathfrak z}_{\mathfrak p}({\mathfrak a})$ is a toral algebra.
Thus ${\mathfrak z}_{\mathfrak p}({\mathfrak a})={\mathfrak a}$.
%\qed
\end{proof}

A torus $A$ in $G$ is {\it $\theta$-split} or {\it $\theta$-anisotropic} if $\theta(a)=a^{-1}$ for all $a\in A$.

\begin{lemma}\label{maxsplittori}
Let ${\mathfrak a}$ be a maximal torus of ${\mathfrak p}$.
Then there is a unique maximal $\theta$-split torus $A$ of $G$ such that ${\mathfrak a}=\Lie(A)$.
\end{lemma}

\begin{proof}
Let $L=Z_G({\mathfrak a})$ and let ${\mathfrak l}=\Lie(L)={\mathfrak z}_{\mathfrak g}({\mathfrak a})$.
Since ${\mathfrak l}\cap{\mathfrak p}={\mathfrak a}\subseteq {\mathfrak z}({\mathfrak l})$, $\theta|_{L^{(1)}}$ is trivial by Lemma \ref{nontriv}.
Let $S$ be any maximal torus of $L$: then $S=(S\cap L^{(1)})\cdot Z(L)^\circ$.
Hence $A=S_-\subset Z(L)$.
Moreover, ${\mathfrak a}\subseteq{\mathfrak z}({\mathfrak l})\subseteq\Lie(S)$ by \cite[2.3]{me}.
It follows that $\Lie(A)={\mathfrak a}$.
It remains to prove uniqueness.
But $A\subset Z(L)$, hence $A$ is the unique maximal $\theta$-split torus of $L$.
%\qed
\end{proof}

\subsection{Summary of Results On Maximal $\theta$-split Tori}
\label{sec:2.2}

The main idea of \cite{kostrall} is that the pair $(G^\theta,{\mathfrak p})$ (with $G^\theta$ acting on ${\mathfrak p}$ via the adjoint representation) can be thought of as a generalised version of the pair $(G,{\mathfrak g})$.
In the new setting the role of Cartan subalgebra of ${\mathfrak g}$ is taken by the maximal toral algebra ${\mathfrak a}$ of ${\mathfrak p}$.
By Lemma \ref{maxsplittori} there exists a maximal $\theta$-split torus $A$ of $G$ such that $\Lie(A)={\mathfrak a}$.
Hence it is useful to recall some results of Vust, Richardson, and Springer concerning maximal $\theta$-split tori.

By Vust we have (\cite[\S 1]{vust}):

{\it - Any two maximal $\theta$-split tori of $G$ are conjugate by an element of $K$.}

It follows immediately from Lemma \ref{maxsplittori} that any two maximal tori in ${\mathfrak p}$ are conjugate by an element of $K$ (this also follows from Thm. \ref{carts} below).
Let $F$ be the finite group of all $a\in A$ satisfying $a^2=e$, the identity element of $G$.
It is easy to see that $F\subset G^\theta$, hence that $F$ normalizes $K$.
Moreover:

{\it - $G^\theta=F\cdot K$ (\cite[\S 1]{vust})}.

If $G$ is not adjoint, we are in fact more interested in the group $K^*=\{ g\in G\,|\,g^{-1}\theta(g)\in Z(G)\}$ introduced by Richardson in \cite{rich2}.
Let $\pi:G\longrightarrow G/Z(G)=\overline{G}$ be the projection onto the adjoint quotient $\overline{G}$, and let $\overline\theta$ be the unique involutive automorphism of $\overline{G}$ making the following diagram commutative:

\begin{diagram}
G & \rTo^{\theta} & G \\
\dTo^\pi & & \dTo^\pi \\
\overline{G} & \rTo^{\overline\theta} & \overline{G} \\
\end{diagram}

Then $K^*=\pi^{-1}(\overline{G}^{\overline\theta})$.
We have (see \cite[8.1]{rich2}):

{\it - $F^*$ normalizes $K$ and $K^*=F^*\cdot K$.}

Let $\Phi_A=\Phi(G,A)$, the roots of $G$ relative to $A$, let $S$ be a maximal torus of $G$ containing $A$, let $\Phi_S=\Phi(G,S)$ and let $W_S=W(G,S)$.
By \cite[2.6(iv)]{rich2} $S$ is $\theta$-stable.
Denote by $\theta^*$ the automorphism of $\Phi_S$ induced by $\theta$.
A parabolic subgroup $P$ of $G$ is {\it $\theta$-split} if $P\cap\theta(P)$ is a Levi subgroup of $P$ (and therefore also of $\theta(P)$).
By Vust \cite[\S 1]{vust}:

{\it - Let $P\supset A$ be a $\theta$-split parabolic subgroup of $G$.
Then $P$ is a minimal $\theta$-split parabolic if and only if $P\cap\theta(P)=Z_G(A)$.
Any two minimal $\theta$-split parabolic subgroups of $G$ are conjugate by an element of $K$.}

Fix a minimal $\theta$-split parabolic subgroup $P$ of $G$ containing $S$ and let $B$ be a Borel subgroup of $G$ such that $S\subset B\subset P$.
Let $\Delta_S$ be the corresponding basis of simple roots in $\Phi_S$.
For a subset $I$ of $\Delta_S$, denote by $\Phi_I$ the subsystem of $\Phi_S$ generated by $\{\alpha:\alpha\in I\}$, by $W_I$ the subgroup of $W_S$ generated by $\{s_\alpha:\alpha\in I\}$, and by $w_I$ the longest element of $W_I$ relative to this Coxeter basis.
By \cite[1.3-4]{springer} (established in \cite{springer2}) we have:

\begin{lemma}\label{basis}
There is a subset $I$ of $\Delta_S$ and a graph automorphism $\psi$ of $\Phi_S$ such that:

(i) $\psi(\Delta_S)=\Delta_S$ and $\psi(I)=I$,

(ii) $\theta^*(\alpha)=-w_I(\psi(\alpha))=-\psi(w_I(\alpha))$ for all $\alpha\in\Phi_S$,

(iii) $\theta^*(\alpha)=\alpha$ for any $\alpha\in\Phi_I$.

The maximal $\theta$-split torus $A'=A\cap G^{(1)}$ of $G^{(1)}$ can be characterised as follows: $A'=\{ s\in S\cap G^{(1)}\,|\,\alpha(s)=1,\beta(s)=\psi(\beta)(s):\alpha\in I,\beta\in\Delta_S\setminus I\}^\circ$.
\end{lemma}

It follows that $\Pi=\{\alpha|_A\, :\, \alpha\in\Delta_S\setminus I\}$ is a basis for $\Phi_A$.
Note that for $\alpha,\beta\in\Delta_S\setminus I$, $\alpha|_A=\beta|_A$ if and only if $\beta\in\{\alpha,\psi(\alpha)\}$.
(We will use $\Delta$ or $\Delta_T$ to denote a basis of roots relative to a maximal torus $T$ of $G$, and $\Pi$ to denote a basis of simple roots in $\Phi_A$, where $A$ is a maximal $\theta$-split torus of $G$.)

The `baby Weyl group' $W_A=N_G(A)/Z_G(A)$ was described by Richardson \cite[\S 4]{rich2}:

{\it - Let $W_1=\{ w\in W_S\,|\, w(A)=A\}$, $W_2=\{ w\in W_1\,|\, w|_A=1|_A\,\}$.
Then the restriction $w\mapsto w|_A$ induces an isomorphism $W_1/W_2\rightarrow W_A$.}

Let $\Gamma$ be the group of automorphisms of $S$ generated by $W$ and $\theta$, let $X(S)$ be the group of characters of $S$ and let $E=X(S)\otimes_{\mathbb Z}{\mathbb R}$.
There exists a $\Gamma$-equivariant inner product $(.\, ,.):E\times E\rightarrow{\mathbb R}$.
Let $E_-$ be the $(-1)$ eigenspace for $\theta$: $E_-$ identifies naturally with $X(A)\otimes_{\mathbb Z}{\mathbb R}$.
Hence $(.\, ,.)$ restricts to a $W_A$-equivariant inner product on $E_-$.
Let $Y(S)$ be the group of cocharacters in $S$.
The dual space $E^*$ to $E$ identifies naturally with $Y(S)\otimes_{\mathbb Z}{\mathbb R}$, and the $(-1)$ eigenspace $E_-^*$ identifies with $Y(A)\otimes_{\mathbb Z}{\mathbb R}$.
Hence the inner product $(.\, ,.)$ induces a $\Gamma$-equivariant isomorphism $E\rightarrow E^*$, which restricts to a $W_A$-equiviarant isomorphism $E_-\rightarrow E_-^*$.
Let $\langle .\, ,.\rangle:X(A)\times Y(A)\longrightarrow{\mathbb Z}$ be the natural pairing.
For $\beta\in\Phi_A$, denote by $s_\beta$ the reflection in the hyperplane orthogonal to $\beta$.
If $\alpha,\beta\in\Phi_A$, then by abuse of notation we write $\langle\alpha,\beta\rangle$ for $2(\alpha,\beta)/(\beta,\beta)$: hence $s_\beta(\alpha)=\alpha-\langle\alpha,\beta\rangle\beta$.

{\it - The set $\Phi_A$ is a (non-reduced) root system in $X(A)$ with Cartan integers $\langle\alpha,\beta\rangle\in{\mathbb Z}$.
The Weyl group $W_A$ is generated by the reflections $\{ s_\alpha\,:\, \alpha\in\Phi_A\}$, hence by the set $\{ s_\alpha\,:\,\alpha\in\Pi\}$.
Each element of $W_A$ has a representative in $K$.
Thus $W_A\cong N_K(A)/Z_K(A)$ (\cite[\S 4]{rich2}).}

Note that it follows from Lemma \ref{maxsplittori} that $N_G(A)=N_G({\mathfrak a})$ and $Z_G(A)=Z_G({\mathfrak a})$.
Let $\Phi_A^*$ be the set of $\alpha\in\Phi_A$ such that $\alpha/m\in\Phi_A\Rightarrow m=\pm 1$.
It follows from the above that $\Phi_A^*$ is a reduced root system.
Finally, we observe using the classification of involutions (see
Springer, \cite{springer}):

\begin{lemma}\label{pisgood}
If $p$ is good for $G$, then it is also good for $\Phi_A$.
If $\alpha\in\Phi_A$, then $3\alpha\notin\Phi_A$.
\end{lemma}

\subsection{Cartan subspaces}
\label{sec:2.3}

Let ${\mathfrak h}$ be a nilpotent subalgebra of ${\mathfrak g}$.
We recall (Fitting's Lemma, see \cite[II.4]{jac}) that there is a decomposition ${\mathfrak g}={\mathfrak g}^0({\mathfrak h})\oplus{\mathfrak g}^1({\mathfrak h})$ and a Zariski open subset $U$ of ${\mathfrak h}$ such that $(\ad u)$ is nilpotent on ${\mathfrak g}^0({\mathfrak h})$ and is non-singular on ${\mathfrak g}^1({\mathfrak h})$ for all $u\in U$.

The following lemma appears in \cite{kostrall}.
We include the proof (which is identical to Kostant-Rallis') for the readers' convenience.

\begin{lemma}\label{fitting}
Let ${\mathfrak h}$ be a nilpotent subalgebra of ${\mathfrak g}$ contained in ${\mathfrak p}$.
Then
$${\mathfrak g}^i({\mathfrak h})=({\mathfrak g}^i({\mathfrak h})\cap{\mathfrak k})\oplus({\mathfrak g}^i({\mathfrak h})\cap{\mathfrak p})\;\; \mbox{for}\;\;\; i=0,1.$$
\end{lemma}

\begin{proof}
Let $y\in U\subseteq {\mathfrak h}$, where $U$ is the subset of ${\mathfrak h}$ defined above.
Since $(\ad y)$ is nilpotent (resp. non-singular) on ${\mathfrak g}^0({\mathfrak h})$ (resp. ${\mathfrak g}^1({\mathfrak h})$), then the same is true of $(\ad y)^2$.
But $(\ad y)^2$ also stabilises ${\mathfrak k}$ and ${\mathfrak p}$.
Hence ${\mathfrak g}^i(k(\ad y)^2)={\mathfrak g}^i(k(\ad y)^2)\cap{\mathfrak k}\oplus{\mathfrak g}^i(k(\ad y)^2)\cap{\mathfrak p}$ for $i=0,1$.
%\qed
\end{proof}

Following \cite{kostrall}, we define a {\it Cartan subspace} of ${\mathfrak p}$ to be a nilpotent algebra ${\mathfrak h}\subseteq{\mathfrak p}$ such that ${\mathfrak g}^0({\mathfrak h})\cap{\mathfrak p}={\mathfrak h}$.

\begin{lemma}\label{torcart}
Let ${\mathfrak a}$ be a maximal torus of ${\mathfrak p}$.
Then ${\mathfrak a}$ is a Cartan subspace.
\end{lemma}

\begin{proof}
As ${\mathfrak a}$ is a toral algebra, ${\mathfrak g}$ is a completely reducible $(\ad{\mathfrak a})$-module.
Thus ${\mathfrak g}^0({\mathfrak a})={\mathfrak z}_{\mathfrak g}({\mathfrak a})$.
By Lemma \ref{torcent}, ${\mathfrak z}_{\mathfrak p}({\mathfrak a})={\mathfrak a}$.
Hence by Lemma \ref{fitting}, ${\mathfrak g}^0({\mathfrak a})\cap{\mathfrak p}={\mathfrak a}$.
%\qed
\end{proof}

Let $x\in{\mathfrak p}$.
Then $kx$ is a nilpotent subalgebra of ${\mathfrak g}$.
We write ${\mathfrak g}^i(x)$ for ${\mathfrak g}^i(kx)$.
Let $q=\min\{ \dim({\mathfrak g}^0(x)\cap{\mathfrak p})\}$, and let $Q=\{ x\in{\mathfrak p}\,|\,\dim({\mathfrak g}^0(x)\cap{\mathfrak p})=q\}$.
It is easy to see that $\dim({\mathfrak g}^0(x)\cap{\mathfrak p})$ is the degree of the first non-zero term in the characteristic polynomial of $(\ad x)^2|_{\mathfrak p}$.
Hence $Q$ is a non-empty open subset of ${\mathfrak p}$.
The following result follows immediately from the proof of \cite[Lemma 3]{kostrall}, although it is not explicitly stated there.
The proof is similar to Richardson's proof of \cite[3.3]{rich2}.

\begin{lemma}\label{sep}
Let $x\in{\mathfrak p}$.
Then the map $\pi:K\times ({\mathfrak g}^0(x)\cap{\mathfrak p})\longrightarrow{\mathfrak p}$ given by  $(k,y)\mapsto \Ad k(y)$ is a separable morphism.
\end{lemma}

\begin{proof}
We consider the differential of $\pi$ at $(e,x)$, where $e$ is the identity element of $G$.
Identify the tangent spaces $T_x({\mathfrak g}^0(x)\cap{\mathfrak p})$ and $T_x({\mathfrak p})$ with $({\mathfrak g}^0(x)\cap{\mathfrak p})$ and ${\mathfrak p}$ respectively.
Hence $d\pi_{(e,x)}:{\mathfrak k}\oplus({\mathfrak g}^0(x)\cap{\mathfrak p})\longrightarrow{\mathfrak p}$, $(U,V)\mapsto [U,x]+V$.
Therefore $d\pi_{(e,x)}({\mathfrak k}\oplus({\mathfrak g}^0(x)\cap{\mathfrak p}))=[x,{\mathfrak k}]+({\mathfrak g}^0(x)\cap{\mathfrak p})$.
By the properties of the Fitting decomposition, $(\ad x)$ is non-singular on ${\mathfrak g}^1(x)$, hence $(\ad x)^2$ is non-singular on $({\mathfrak g}^1(x)\cap{\mathfrak p})$.
Thus $[x,{\mathfrak k}]\supseteq[x,[x,{\mathfrak p}]]\supseteq[x,[x,({\mathfrak g}^1(x)\cap{\mathfrak p})]]=({\mathfrak g}^1(x)\cap{\mathfrak p})$.
It follows that $d\pi_{(e,x)}$ is surjective.
By \cite[AG. 17.3]{bor} $\pi$ is separable.
%\qed
\end{proof}

\begin{corollary}\label{cart}
Let ${\mathfrak h}$ be a Cartan subspace of ${\mathfrak p}$.
The map $\pi:K\times{\mathfrak h}\longrightarrow{\mathfrak p}$ given by $(g,h)\mapsto \Ad g(h)$ is separable, and $K\cdot{\mathfrak h}$ contains a dense open subset of ${\mathfrak p}$.
\end{corollary}

We can now prove the main theorem of this section.
Our proof is somewhat shorter than the proof given in \cite{kostrall}.

\begin{theorem}\label{carts}
Any two Cartan subspaces of ${\mathfrak p}$ are $K$-conjugate.
The Cartan subspaces are just the maximal tori of ${\mathfrak p}$.
An element $x\in{\mathfrak p}$ is semisimple if and only if it is contained in a Cartan subspace of ${\mathfrak p}$.
\end{theorem}

\begin{proof}
Let ${\mathfrak h}$ be a Cartan subspace.
Let $U$ be the open subset of elements $h\in{\mathfrak h}$ such that ${\mathfrak g}^i(h)={\mathfrak g}^i({\mathfrak h})$ for $i=0,1$.
By Cor. \ref{cart}, $K\cdot U$ contains a dense open subset of ${\mathfrak p}$.
Hence $(K\cdot U)\cap Q$ is non-empty.
But $Q$ is $K$-stable, hence $U\cap Q$ is non-empty.
Let $u\in U\cap Q$.
Then ${\mathfrak g}^0(u)\cap{\mathfrak p}={\mathfrak h}$.
Therefore $\dim{\mathfrak h}=q$.
On the other hand, if $u\in{\mathfrak h}\cap Q$, then ${\mathfrak g}^0(u)\cap{\mathfrak p}\supseteq{\mathfrak h}$, hence ${\mathfrak g}^0(u)\cap{\mathfrak p}={\mathfrak h}$.
It follows that $U=Q\cap{\mathfrak h}$.

Let ${\mathfrak h}'$ be any other Cartan subspace.
Then $K\cdot(Q\cap{\mathfrak h})$ and $K\cdot(Q\cap{\mathfrak h}')$ contain non-empty open subsets of ${\mathfrak p}$, hence their intersection is non-empty.
Therefore $(K\cdot(Q\cap{\mathfrak h}))\cap{\mathfrak h}'$ is non-empty.
It follows that $g\cdot{\mathfrak h}={\mathfrak h}'$ for some $g\in K$.
The remaining statements of the theorem follow at once.
%\qed
\end{proof}

\section{A $\theta$-stable reduction}
\label{sec:3}

We assume from this point on that $G$ has the following three properties:

(A) $p$ is good for $G$.

(B) The derived subgroup $G^{(1)}$ is simply-connected.

(C) There exists a symmetric $G$-invariant non-degenerate bilinear form $B:{\mathfrak g}\times{\mathfrak g}\longrightarrow k$.

In this section we will prove a $\theta$-stable analogue of a result of Gordon and Premet (\cite[6.2]{gandp}).
An important corollary is that the trace form in (C) may be chosen so that it is invariant with respect to $\theta$.

Let $G_i\, (1\leq i\leq l)$ be the minimal normal subgroups of $G^{(1)}$ and let ${\mathfrak g}_i=\Lie(G_i)$.
As $G^{(1)}$ is simply-connected, $G^{(1)}=G_1\times G_2\times\ldots\times G_l$ and ${\mathfrak g}'={\mathfrak g}_1\oplus{\mathfrak g}_2\oplus\ldots\oplus{\mathfrak g}_l$.
We introduce new groups $\tilde{G}_i$, defined as follows:

$$\tilde{G}_i=\left\{
\begin{array}{ll}
GL(V_i) & \mbox{if $G_i$ is isomorphic to $\SL(V_i)$ and $p\,|\dim V_i$,} \\
G_i & \mbox{otherwise.}
\end{array}
\right.
$$

Let $\tilde{G}=\tilde{G}_1\times\tilde{G}_2\times\ldots\times\tilde{G}_l,\tilde{\mathfrak g}_i=\Lie(\tilde{G}_i),\tilde{\mathfrak g}=\Lie(\tilde{G})$.
Identify $G_i$ with the derived subgroup of $\tilde{G}_i$, hence consider $G^{(1)}$ as a subgroup of both $G$ and $\tilde{G}$.

Let $(T',B')$ be a fundamental pair for $\theta|_{G^{(1)}}$ (see Sect. \ref{sec:1}) and let $T$ (resp. $\tilde{T}$) be the unique maximal torus of $G$ (resp. $\tilde{G}$) containing $T'$.
Let ${\mathfrak h}'=\Lie(T'),{\mathfrak h}=\Lie(T),\tilde{\mathfrak h}=\Lie(\tilde{T})$, ${\mathfrak h}_i={\mathfrak h}\cap{\mathfrak g}_i$, $\tilde{\mathfrak h}_i=\tilde{\mathfrak h}\cap\tilde{\mathfrak g}_i$.

\begin{theorem}\label{redthm}
There exists a torus $T_0$, an involution $\hat{\theta}$ of $\hat{G}=\tilde{G}\times T_0$, and an injective restricted Lie algebra homomorphism $\psi:{\mathfrak g}\longrightarrow\hat{\mathfrak g}=\Lie(\hat{G})$ such that:

(i) $\psi({\mathfrak g}_i)\subseteq\tilde{\mathfrak g}_i$ for all $i\in\{1,2,\ldots, l\}$ and $\psi({\mathfrak h}')\subseteq\tilde{\mathfrak h}$.

(ii) $\hat{\theta}|_{G^{(1)}}=\theta|_{G^{(1)}}$, and the following diagram is commutative:
\begin{diagram}
{\mathfrak g} & \rTo^{d\theta} & {\mathfrak g} \\
\dTo^\psi & & \dTo^\psi \\
\hat{\mathfrak g} & \rTo^{d\hat{\theta}} & \hat{\mathfrak g} \\
\end{diagram}

(iii) There exists a toral algebra ${\mathfrak t}_1$ such that $\hat{\mathfrak g}=\psi({\mathfrak g})\oplus{\mathfrak t}_1$ (Lie algebra direct sum) and $d\hat{\theta}(t)=t\;\forall\, t\in{\mathfrak t}_1$.

(iv) $\theta(G_i)=G_j$ implies $\hat{\theta}(\tilde{G}_i)=\tilde{G}_j$.
\end{theorem}

\begin{proof}
The existence of a torus $T_0$, an injective restricted Lie algebra homomorphism $\eta:{\mathfrak g}\longrightarrow\hat{\mathfrak g}=\Lie(\tilde{G}\times T_0)=\tilde{\mathfrak g}\oplus{\mathfrak t}_0$, and a toral algebra ${\mathfrak s}_1$ such that $\hat{\mathfrak g}=\eta({\mathfrak g})\oplus{\mathfrak s}_1$ was proved by Premet \cite[Lemma 4.1]{comp} and Gordon-Premet \cite[6.2]{gandp}.
Identify each ${\mathfrak g}_i$ with its image $\eta({\mathfrak g}_i)\subseteq\tilde{\mathfrak g}_i$.
Define an automorphism $\phi$ of the restricted Lie algebra $\hat{\mathfrak g}$ by $\phi(\eta(x))=\eta(d\theta(x))$ for $x\in{\mathfrak g}$, $\phi(s)=s$ for $s\in{\mathfrak s}_1$ and linear extension to all of $\hat{\mathfrak g}$.

The main idea of our proof is to find $\phi$-stable restricted subalgebras $\overline{\mathfrak g}_i,{\mathfrak s}_0$, and $\overline{\mathfrak g}_i\oplus\overline{\mathfrak g}_j$ of $\hat{\mathfrak g}$ with ${\mathfrak g}_i\subseteq\overline{\mathfrak g}_i\cong\tilde{\mathfrak g}_i,{\mathfrak s}_0\cong{\mathfrak t}_0$ and $\hat{\mathfrak g}=\sum\overline{\mathfrak g}_i\oplus{\mathfrak s}_0$.

{\bf Step 1. The toral algebra ${\mathfrak s}_0$.}

Let $\hat{\mathfrak z}={\mathfrak z}(\hat{\mathfrak g}),\tilde{\mathfrak z}={\mathfrak z}(\tilde{\mathfrak g})$ and ${\mathfrak z}_i={\mathfrak z}({\mathfrak g}_i)$.
Clearly $\hat{\mathfrak z}=\tilde{\mathfrak z}\oplus{\mathfrak t}_0=\eta({\mathfrak z})\oplus{\mathfrak s}_1$ and $\tilde{\mathfrak z}=\sum{\mathfrak z}_i={\mathfrak z}({\mathfrak g}')$.
Hence $\tilde{\mathfrak z}\subseteq\hat{\mathfrak z}$ are $\phi$-stable toral algebras.
The restriction of $\phi$ to ${\hat{\mathfrak z}}^{tor}$ has order 1 or 2.
Therefore by Maschke's theorem there is a $\phi$-stable ${\mathbb F}_p$-vector space ${\mathfrak s}_0^{tor}$ such that ${\hat{\mathfrak z}}^{tor}={\tilde{\mathfrak z}}^{tor}\oplus{\mathfrak s}_0^{tor}$.

Let ${\mathfrak s}_0$ be the toral algebra in $\hat{\mathfrak g}$ generated by ${\mathfrak s}_0^{tor}$.
Using the same argument as in the proof of Lemma \ref{redcase} we can choose a toral basis for ${\mathfrak s}_0$ consisting of eigenvectors for $\phi$.
This basis can be used to construct an isomorphism of toral algebras $f_0:{\mathfrak s}_0\longrightarrow{\mathfrak t}_0$ and an involutive automorphism $\theta_0:T_0\longrightarrow T_0$ such that the following diagram commutes:

\begin{diagram}
{\mathfrak s}_0 & \rTo^{f_0} & {\mathfrak t}_0 \\
\dTo^\phi & & \dTo^{d\theta_0} \\
{\mathfrak s}_0 & \rTo^{f_0} & {\mathfrak t}_0 \\
\end{diagram}

{\bf Step 2. The subalgebra $\overline{\mathfrak g}_i$, for $\theta$-stable $G_i$}

If $\tilde{G}_i=G_i$, there is nothing to prove.
So assume $\tilde{G}_i=\GL(V_i)$ and $p\,|\dim V_i$.
Let $\Delta_i$ be the subset of $\Delta$ corresponding to $G_i$.
We define ${\mathfrak m}_i=\sum_{j\neq i}{\mathfrak g}_j$ and ${\mathfrak n}_i={\mathfrak z}_{\hat{\mathfrak g}}({\mathfrak m}_i)\cap\hat{\mathfrak h}$.
Clearly $\sum_{j\neq i}{\mathfrak z}_j\oplus{\mathfrak s}_0\subseteq{\mathfrak n}_i=\sum_{j\neq i}{\mathfrak z}_j\oplus\tilde{\mathfrak h}_i\oplus{\mathfrak s}_0\subseteq\hat{\mathfrak h}$ are $\phi$-stable toral algebras.
Hence there is a $\phi$-stable toral algebra $\overline{\mathfrak h}_i$ containing ${\mathfrak h}_i$ such that ${\mathfrak n}_i=\overline{\mathfrak h}_i\oplus\sum_{j\neq i}{\mathfrak z}_j\oplus{\mathfrak s}_0$.

By \cite[4.2]{me}, the maps $d\alpha|_{\overline{\mathfrak h}_i}$ with $\alpha\in\Delta_i$ are linearly independent.
It follows that $\overline{\mathfrak h}_i$ and ${\mathfrak g}_i$ together generate a restricted Lie algebra isomorphic to $\tilde{\mathfrak g}_i$.
Let $f_i:\overline{\mathfrak g}_i\longrightarrow\tilde{\mathfrak g}_i$ be an isomorphism such that $f_i(x)=x$ for all $x\in{\mathfrak g}_i$.
Then by Lemma \ref{GLautos} there exists a unique involutive automorphism $\theta_i:\tilde{G}_i\longrightarrow\tilde{G}_i$ such that the following diagram commutes:

\begin{diagram}
\overline{\mathfrak g}_i & \rTo^{f_i} & \tilde{\mathfrak g}_i \\
\dTo^\phi & & \dTo^{d\theta_i} \\
\overline{\mathfrak g}_i & \rTo^{f_i} & \tilde{\mathfrak g}_i \\
\end{diagram}

{\bf Step 3. The subalgebras $\overline{\mathfrak g}_i,\overline{\mathfrak g}_j$ when $\theta(G_i)=G_j$.}

Once again we may assume that $\tilde{G}_i=\GL(V_i)$ and $p|\dim V_i$.
We set $\overline{\mathfrak g}_i=\tilde{\mathfrak g}_i,\overline{\mathfrak g}_j=\phi(\tilde{\mathfrak g}_i)$.
We have only to show that $\hat{\mathfrak g}=\overline{\mathfrak g}_i\oplus\overline{\mathfrak g}_j\oplus\sum_{k\neq i,j}\tilde{\mathfrak g}_k\oplus{\mathfrak s}_0$.
Let $\Delta_i,\Delta_j$ be the subsets of $\Delta$ corresponding respectively to $G_i,G_j$ and let ${\mathfrak n}_{(i,j)}=\{ h\in\hat{\mathfrak h}|\,d\alpha(h)=0\,\forall\,\alpha\in\Delta\setminus(\Delta_i\cup\Delta_j)\}$.
Clearly ${\mathfrak n}_{(i,j)}=\tilde{\mathfrak h}_i\oplus\tilde{\mathfrak h}_j\oplus\sum_{k\neq i,j}{\mathfrak z}_k\oplus{\mathfrak s}_0$.
The automorphism of $\Phi$ induced by $\theta$ sends $\Delta_i$ onto $\Delta_j$.
Hence the differentials $d\alpha|_{\tilde{\mathfrak h}_i\oplus d\theta(\tilde{\mathfrak h}_i)}$ for $\alpha\in\Delta_i\cup\Delta_j$ are linearly independent.
It follows by dimensional considerations that $\tilde{\mathfrak h}_i\oplus d\theta(\tilde{\mathfrak h}_i)\oplus\sum_{k\neq i,j}{\mathfrak z}_k\oplus{\mathfrak s}_0 = {\mathfrak n}_{(i,j)}$.
Therefore $\tilde{\mathfrak g}_i\oplus d\theta(\tilde{\mathfrak g}_i)\oplus\sum_{k\neq i,j}\tilde{\mathfrak g}_k\oplus{\mathfrak s}_0=\hat{\mathfrak g}$.

It is now easy to see that there are isomorphisms $f_j:\overline{\mathfrak g}_j\longrightarrow\tilde{\mathfrak g}_j,\tau_j:\tilde{G}_i\longrightarrow\tilde{G}_j$ and $\theta_{(i,j)}:\tilde{G}_i\times\tilde{G}_j$ such that $f_j(x)=x\;\forall x\in{\mathfrak g}_i$ and the following diagram is commutative:

\begin{diagram}
\overline{\mathfrak g}_i\oplus\overline{\mathfrak g}_j & \rTo^{(\Id,f_j)} & \tilde{\mathfrak g}_i\oplus\tilde{\mathfrak g}_j \\
\dTo^\phi & & \dTo^{d\theta_{(i,j)}} \\
\overline{\mathfrak g}_i\oplus\overline{\mathfrak g}_j & \rTo^{(\Id,f_j)} & \tilde{\mathfrak g}_i\oplus\tilde{\mathfrak g}_j \\
\end{diagram}

where $\theta_{(i,j)}:\tilde{G}_i\times\tilde{G}_j\longrightarrow\tilde{G}_i\times\tilde{G}_j$ is given by $(g_i,g_j)\mapsto (\tau^{-1}(g_j),\tau(g_i))$.

We now let $f:\sum\overline{\mathfrak g}_i\oplus{\mathfrak s}_0=\hat{\mathfrak g}\longrightarrow\sum\tilde{\mathfrak g}_i\oplus{\mathfrak t}_0=\hat{\mathfrak g}$ and $\hat\theta:\tilde{G}\times T_0\longrightarrow\tilde{G}\times T_0$ be the maps obtained in the obvious way from the $f_i$ and the $\theta_i,\theta_{(i,j)}$ respectively.
Then the following diagram is commutative:

\begin{diagram}
\hat{\mathfrak g} & \rTo^{f} & \hat{\mathfrak g} \\
\dTo^\phi & & \dTo^{d\hat{\theta}} \\
\hat{\mathfrak g} & \rTo^{f} & \hat{\mathfrak g} \\
\end{diagram}

Let $\psi=f\circ\eta:{\mathfrak g}\longrightarrow\hat{\mathfrak g}$ and let ${\mathfrak t}_1=f({\mathfrak s}_1)$.
Then $\psi,\tilde{\mathfrak g}_i,T_0,{\mathfrak t}_1$ satisfy the requirements of the theorem.
%\qed
\end{proof}

\begin{corollary}\label{trace}
Let $G$ satisfy the standard hypotheses (A),(B),(C).
Suppose that $\Char k\neq 2$ and that $\theta$ is an involutive automorphism of $G$.
Then the trace form in (C) may be chosen to be $\theta$-equivariant.
\end{corollary}

\begin{proof}
To prove the corollary we construct a $\hat{\theta}$-equivariant trace form on $\hat{\mathfrak g}$ which restricts to a non-degenerate form on ${\mathfrak g}$.
Recall that $\hat{\mathfrak g}=\tilde{\mathfrak g}\oplus{\mathfrak t}_0=\psi({\mathfrak g})\oplus{\mathfrak t}_1$.
Identify ${\mathfrak g}$ with its image $\psi({\mathfrak g})$.
Let $G_i$ be a minimal normal subgroup of $G$.
As is well-known (see for example \cite[I.5]{sands}) there exists a non-degenerate trace form $\kappa_i:\tilde{\mathfrak g}_i\times\tilde{\mathfrak g}_i\longrightarrow k$ associated to a rational representation of $\tilde{G}_i$.
Moreover, as $\tilde{\mathfrak g}_i$ is an indecomposable $\tilde{G}_i$-module, $\kappa_i$ is unique up to multiplication by a non-zero scalar.
We will prove that $\kappa_i$ is invariant under any automorphism of $\tilde{G}_i$.

By Lemma \ref{GLautos} it suffices to prove this for a set of graph automorphisms $\gamma$ generating $\Aut\tilde{G}_i / \Int\tilde{G}_i$.
Let $\gamma$ be such an automorphism and define a new trace form $\kappa_i^\gamma:(x,y)\mapsto \kappa_i(d\gamma(x),d\gamma(y))$.
Then $\kappa_i^\gamma$ is a scalar multiple of $\kappa_i$.
Hence it will suffice to find $(x,y)\in\tilde{\mathfrak g}_i\times\tilde{\mathfrak g}_i$ such that $\kappa_i(x,y)=\kappa_i^\gamma(x,y)\neq 0$.
Assume first of all that $G_i$ is not of type $A$ (therefore $\tilde{G}_i=G_i)$.
Let $(B_i,T_i)$ be a fundamental pair for $\gamma$ and let $\Delta_i$ be the basis of the roots $\Phi_i=\Phi(G_i,T_i)$ corresponding to $B_i$.
Let $\{h_{\alpha_i},e_\alpha|\alpha_i\in\Delta_i,\alpha\in\Phi_i\}$ be a Chevalley basis for ${\mathfrak g}_i$.

We observe first of all that there exists $\alpha\in\Delta_i$ such that $\gamma(\alpha)=\alpha$.
For type $D_n$ we choose $\alpha=\alpha_{n-2}$, and for type $E_6$ we choose $\alpha=\alpha_2$ (we use Bourbaki's numbering conventions \cite{bourbaki}).
We have $d\gamma(e_\alpha)=ce_\alpha$ and $d\gamma(e_{-\alpha})=c'e_{-\alpha}$.
But $[e_\alpha,e_{-\alpha}]=h_\alpha$, hence $cc'=1$.
Therefore $\kappa_i^\gamma(e_\alpha,e_{-\alpha})=\kappa_i(e_\alpha,e_{-\alpha})$.
$\kappa_i$ is non-degenerate and $T_i$-invariant.
Thus $\kappa_i(e_\alpha,e_{-\alpha})\neq 0$.

Assume now that $G_i$ is of type $A$.
In this case $G_i$ is isomorphic to $\SL(V_i)$ and it will be sufficient to prove $\kappa_i^\gamma=\kappa_i$ for $\gamma:g\mapsto {^t}g^{-1}$.
Recall that the ordinary trace form $\kappa_i(x,y)=\tr(xy)$ is non-degenerate on $\tilde{\mathfrak g}_i$.
Hence $\kappa_i^\gamma(x,y)=\kappa_i(-{^t}x,-{^t}y)=\tr({^t}x {^t}y)=\tr({^t}(yx))=\tr(yx)=\tr(xy)=\kappa_i(x,y)$.

To construct the form $\hat{\kappa}$ we proceed as follows.
For $d\theta$-stable $\tilde{\mathfrak g}_i$ we choose a trace form $\kappa_i$ as above.
For each pair $\tilde{\mathfrak g}_i,\tilde{\mathfrak g}_j$ with $d\theta(\tilde{\mathfrak g}_i)=\tilde{\mathfrak g}_j$ we let $\kappa_i$ be a non-degenerate trace form on $\tilde{\mathfrak g}_i$, and define $\kappa_j$ on $\tilde{\mathfrak g}_j$ by $\kappa_j(x,y)=\kappa_i(d\theta(x),d\theta(y))$.

Let $\hat{\mathfrak z}={\mathfrak z}(\hat{\mathfrak g}),\tilde{\mathfrak z}={\mathfrak z}(\tilde{\mathfrak g}),{\mathfrak z}={\mathfrak z}({\mathfrak g})$.
It is easy to see that $\hat{\mathfrak z}=\tilde{\mathfrak z}\oplus{\mathfrak t}_0={\mathfrak z}\oplus{\mathfrak t}_1$.
Moreover $\tilde{\mathfrak z}={\mathfrak z}({\mathfrak g}')\subseteq{\mathfrak z}$.
Hence ${\mathfrak z}=\tilde{\mathfrak z}\oplus({\mathfrak z}\cap{\mathfrak t}_0)$.
By the same argument as used in the proof of Lemma \ref{redcase} there exists a $\hat{\theta}$-stable toral algebra ${\mathfrak t}_2$ such that ${\mathfrak t}_0={\mathfrak z}\cap{\mathfrak t}_0\oplus{\mathfrak t}_2$.
Let $\kappa_z$ be a non-degenerate $\hat{\theta}$-invariant form on ${\mathfrak z}\cap{\mathfrak t}_0$, and let $\kappa_t$ be such a form on ${\mathfrak t}_2$.
Any $x\in\hat{\mathfrak g}$ can be expressed uniquely as $(\sum x_i)+ x_z + x_t$, with $x_i\in\tilde{\mathfrak g}_i,x_z\in{\mathfrak z}\cap{\mathfrak t}_0$, and $x_t\in{\mathfrak t}_2$.
We define $\kappa(x,y)=\sum\kappa_i(x_i,y_i)+\kappa_z(x_z,y_z)+\kappa_t(x_t,y_t)$

It remains to show that the restriction of $\kappa$ to ${\mathfrak g}$ is non-degenerate.
Let $x\in\hat{\mathfrak g}$ be such that $\kappa(x,y)=0\,\forall y\in{\mathfrak g}$.
Then $\kappa_i(x_i,{\mathfrak g}_i)=0\,\forall i$, hence $x_i\in{\mathfrak z}_i$.
Moreover $\kappa_z(x_z,{\mathfrak z}\cap{\mathfrak t}_0)=0$, hence $x_z=0$.
Suppose $x_i\neq 0$.
Let $\Delta_i=\{ \alpha_1,\alpha_2,\ldots\}$ be the subset of $\Delta$ corresponding to $G_i$, ordered in the standard way.
We have $x_i=\lambda([e_{\alpha_1},e_{-\alpha_1}]+2[e_{\alpha_2},e_{-\alpha_2}]+\ldots)$ and $\lambda\neq 0$.
By \cite[3.3]{me} there exists $h\in\tilde{\mathfrak h}_i$ such that $d\alpha_1(h)=1$, and $d\alpha(h)=0\;\forall\alpha\in\Delta\setminus\{\alpha_1\}$.
Then $\kappa_i(x_i,h)=\lambda\kappa_i(e_{\alpha_1},e_{-\alpha_1})\neq 0$.
This is a contradiction, hence $x_i=0\,\forall i$.

It follows that $x\in{\mathfrak t}_2$.
Therefore the restriction of $\kappa$ to ${\mathfrak g}$ is non-degenerate.
%\qed
\end{proof}

\section{Centralizers and Invariants}
\label{sec:4}

\subsection{Centralizers}
\label{sec:4.1}

The following lemma is an important step in \cite{kostrall}.
Cor. \ref{trace} allows us to prove it by the same argument.

\begin{lemma}\label{centdim}
Let $x\in{\mathfrak p}$.
Then $\dim{\mathfrak z}_{\mathfrak k}(x)-\dim{\mathfrak z}_{\mathfrak p}(x)=\dim{\mathfrak k}-\dim{\mathfrak p}$.
\end{lemma}

\begin{proof}
Let $\kappa:{\mathfrak g}\times{\mathfrak g}\longrightarrow k$ be a non-degenerate $(\theta,G)$-equivariant symmetric bilinear form.
By the $\theta$-equivariance $\kappa({\mathfrak k},{\mathfrak p})=0$.
Let $x\in{\mathfrak p}$ and let $\kappa_x:{\mathfrak g}\times{\mathfrak g}\longrightarrow k$ be the alternating bilinear form defined by $\kappa_x(y,z)=\kappa([x,y],z)=\kappa(y,[z,x])$.
Clearly $\kappa_x(y,z)=0$ for all $z\in{\mathfrak g}$ if and only if $y\in{\mathfrak z}_{\mathfrak g}(x)$.
Hence $\kappa_x$ induces a non-degenerate alternating bilinear form $\overline{\kappa_x}:{\mathfrak g}/{{\mathfrak z}_{\mathfrak g}(x)}\times{\mathfrak g}/{{\mathfrak z}_{\mathfrak g}(x)}\longrightarrow k$.
But now ${\mathfrak g}/{{\mathfrak z}_{\mathfrak g}(x)}={\mathfrak k}/{{\mathfrak z}_{\mathfrak k}(x)}\oplus{\mathfrak p}/{{\mathfrak z}_{\mathfrak p}(x)}$.
Furthermore ${\mathfrak k}/{{\mathfrak z}_{\mathfrak k}(x)}$ and ${\mathfrak p}/{{\mathfrak z}_{\mathfrak p}(x)}$ are $\overline{\kappa_x}$-isotropic subspaces, hence are maximal such, and their dimensions are equal.
%\qed
\end{proof}

The following result will also be useful.

\begin{lemma}\label{globalinf}
Let $x\in{\mathfrak k}$ or ${\mathfrak p}$.
Then $\Lie(Z_G(x)^\circ)={\mathfrak z}_{\mathfrak g}(x)$ and $\Lie(Z_K(x)^\circ)={\mathfrak z}_{\mathfrak k}(x)$.
\end{lemma}

\begin{proof}
Clearly $\Lie(Z_G(x)^\circ)\subseteq{\mathfrak z}_{\mathfrak g}(x)$.
To show that $\Lie(Z_G(x)^\circ)={\mathfrak z}_{\mathfrak g}(x)$, it will therefore suffice to show equality of dimensions.
To do this we use the homomorphism $\psi:{\mathfrak g}\longrightarrow\hat{\mathfrak g}$ of Thm. \ref{redthm}.
It is easy to see that $\dim Z_G(x)^\circ=\dim{\mathfrak z}_{\mathfrak g}(x)$ if and only if $\dim Z_{\hat{G}} (d\psi(x))=\dim {\mathfrak z}_{\hat{\mathfrak g}}(d\psi(x))$.
But equality is known for each of the components $\tilde{G}_i$ (see \cite[I.5.3]{sands}) hence for $\hat{G}$.
Therefore $\Lie(Z_G(x)^\circ)={\mathfrak z}_{\mathfrak g}(x)$.

Now let $L=Z_G(x)^\circ,{\mathfrak l}=\Lie(L)$.
The restriction of $\theta$ to $L$ is a semisimple automorphism, hence $\Lie((L\cap K)^\circ)={\mathfrak l}\cap{\mathfrak k}$ by \cite[\S 9.1]{bor}.
%\qed
\end{proof}

\subsection{Regular Elements}
\label{sec:4.2}

We say that $x\in{\mathfrak p}$ is {\it regular} if $\dim{\mathfrak z}_{\mathfrak k}(x)\leq\dim{\mathfrak z}_{\mathfrak k}(y)$ for all $y\in{\mathfrak p}$.
We denote by ${\cal R}$ the open subset of regular elements in ${\mathfrak p}$.
Let ${\mathfrak a}$ be a Cartan subspace of ${\mathfrak p}$ and let $A$ be a maximal $\theta$-split torus of $G$ such that $\Lie(A)={\mathfrak a}$ (Lemma \ref{maxsplittori}).
We recall (\cite[3.1,3.2]{rich2}) that $Z_G(A)=M\cdot A$ (almost direct product) and ${\mathfrak g}^A={\mathfrak M}\oplus{\mathfrak a}$, where $M=Z_K(A)^\circ,{\mathfrak M}={\mathfrak k}^A=\Lie(M)$, and that $\dim{\mathfrak M}-\dim{\mathfrak a}=\dim{\mathfrak k}-\dim{\mathfrak p}$.

\begin{lemma}\label{regs}
Let $x\in{\mathfrak p}$.
The following are equivalent:

(i) $x$ is regular,

(ii) $\dim{\mathfrak z}_{\mathfrak g}(x)=\dim{\mathfrak a}+\dim{\mathfrak M}$,

(iii) $\dim{\mathfrak z}_{\mathfrak k}(x)=\dim{\mathfrak M}$,

(iv) $\dim{\mathfrak z}_{\mathfrak p}(x)=\dim{\mathfrak a}$.
\end{lemma}

\begin{proof}
Let ${\cal S}$ be the set of semisimple elements in ${\mathfrak p}$, which is a non-empty open subset by Cor. \ref{sep} and Thm. \ref{carts}.
Hence ${\cal S}\cap{\cal R}$ is non-empty.
The equivalence of the four conditions now follows immediately from Lemma \ref{centdim}.
%\qed
\end{proof}

\begin{lemma}
Let $x\in{\mathfrak p}$.
The following are equivalent:

(i) $x$ is regular,

(ii) $K\cdot x$ is a $K$-orbit of maximal dimension in ${\mathfrak p}$,

(iii) $\codim_{\mathfrak p} K\cdot x=\dim{\mathfrak a}$,

(iv) $\codim_{\mathfrak g} G\cdot x = \dim{\mathfrak a}+\dim{\mathfrak M}$.
\end{lemma}

\begin{proof}
This follows immediately from Lemma \ref{globalinf} and Lemma \ref{regs}.
%\qed
\end{proof}

\subsection{Geometric Invariant Theory}
\label{sec:4.3}

Here we briefly recall the definitions and some important facts concerning Mumford's Geometric Invariant Theory.
In positive characteristic this requires the fact that reductive groups are geometrically reductive, proved by Haboush in \cite{haboush}.
For details we refer the reader to \cite{mum,luna,haboush}.

Let $R$ be an affine algebraic group such that the connected component $R^\circ$ is reductive.
Let $X$ be an affine variety on which $R$ acts.
Denote the action by $r\cdot x$ for $r\in R,x\in X$.
We always assume that the map $R\times X\longrightarrow X$, $(r,x)\mapsto r\cdot x$ is a morphism of varieties.
There is an induced action of $R$ on the coordinate ring $k[X]$.
The algebra of invariants $k[X]^R$ is finitely generated.
Hence we can construct the affine variety $X\quot R=\Spec(k[X]^R)$.
The embedding $k[X]^R\hookrightarrow k[X]$ induces a morphism $\pi:X\longrightarrow X\quot R$.

The affine variety $X\quot R$ is the {\it quotient} (of $X$ by $R$) and the map $\pi$ is called the {\it quotient morphism}.
If there is possible ambiguity, we will use the notation $\pi_{X,R}$ or $\pi_X$ for the quotient morphism from $X$ to $X\quot R$.
We have the following facts (see \cite{mum,luna,haboush}):

{\it - $\pi$ is surjective.

 - If $X_1$ and $X_2$ are disjoint closed $R$-stable subsets of $X$, then there exists $f\in k[X]^R$ such that $f(x)=0$ for $x\in X_1$, and $f(x)=1$ for $x\in X_2$.

 - Let $\xi\in X\quot R$. The fibre $\pi^{-1}(\xi)$ is $R$-stable and contains a unique closed $R$-orbit, $T(\xi)$, which is also the unique minimal $R$-orbit in $\pi^{-1}(\xi)$.
Hence $\pi$ determines a bijection between the set of closed $R$-orbits in $X$ and the ($k$-rational) points of $X\quot R$.

 - Let $x\in X$ and let $\xi\in X\quot R$.
Then $\pi(x)=\xi$ if and only if $T(\xi)$ is contained in the closure of $R\cdot x$ in $X$.

 - Suppose $X$ is irreducible, and that there exists $x\in X$ such that $R\cdot x$ is closed and $\dim R\cdot x\geq\dim R\cdot y$ for all $y\in X$.
Then $\pi$ is separable (\cite[9.3]{rich2}).

 - If $X$ is normal, then $X\quot R$ is normal.

 - Let $X,Y$ be two affine varieties admitting (algebraic) $R$-actions and let $f:X\longrightarrow Y$ be an $R$-equivariant morphism of varieties.
There exists a unique morphism $\pi(f):X\quot R\longrightarrow Y\quot R$ such that the following diagram commutes:}

\begin{diagram}
 X & \rTo^f & Y \\
 \dTo^{\pi_{X,R}} & & \dTo^{\pi_{Y,R}} \\
 X\quot R & \rTo^{\pi(f)} & Y\quot R
\end{diagram}

\begin{rk}\label{geoquot}
Let $H$ be a reductive group and let $L_1,L_2$ be commuting reductive subgroups of $H$ such that $H=L_1\cdot L_2$.
Let $X$ be an affine variety on which $H$ acts.
Since $L_1$ commutes with $L_2$, it stabilizes the subring $k[X]^{L_2}$.
Hence $L_1$ acts on the quotient $X\quot L_2$.
Clearly $(k[X]^{L_2})^{L_1}=k[X]^H$.
The quotient $(X\quot L_2)\quot L_1$ therefore identifies naturally with $X\quot H$. We will use the notation $\pi_{X,H/L_2}$ for the morphism $X\quot L_2\rightarrow X\quot H$ induced by the inclusion $k[X]^H\hookrightarrow k[X]^{L_2}$.
(Using the notation above, $\pi_{X,H/L_2}=\pi_{X\quot L_2,L_1}$.)
The following diagram is commutative:
\begin{diagram}
X & \rTo^{\pi_{X,L_2}} & X\quot L_2 \\
\dTo^{\pi_{X,L_1}} & & \dTo^{\pi_{X,H/L_2}} \\
X\quot L_1 & \rTo^{\pi_{X,H/L_1}} & X\quot H
\end{diagram}
\end{rk}

\subsection{Unstable and closed $K$-orbits}
\label{sec:4.4}

Let $\rho:G\longrightarrow\GL(V)$ be a rational representation.
For $U\subset V$, we denote by $\overline{U}$ the closure of $U$ in $V$ (in the Zariski topology).
Recall that an element $v\in V$ is {\it $G$-unstable} if $0\in\overline{\rho(G)(v)}$.
It is well-known that if $\rho$ is the adjoint representation then an element of ${\mathfrak g}$ is $G$-unstable if and only if it is nilpotent.
(This is true even if the characteristic is bad, see \cite[9.2.1]{barrich}.)

\begin{lemma}\label{unstable}
Let $x\in{\mathfrak p}$.
Then $x$ is $K$-unstable if and only if it is nilpotent.
\end{lemma}

\begin{proof}
Let $x\in{\mathfrak p}$ be $K$-unstable.
Then $0\in\overline{K\cdot x}\subseteq\overline{G\cdot x}$, hence $x$ is $G$-unstable, therefore nilpotent.
Suppose on the other hand that $x$ is nilpotent.
Let $(B,T)$ be a fundamental pair for $\theta$, let $\Phi=\Phi(G,T)$, let $\Delta$ be the basis of $\Phi$ corresponding to $B$ and let $H=H(\Phi,\Delta)$ be the group of ${\mathbb Z}$-linear maps from the root lattice of $\Phi$ to ${\mathbb Z}$.
By Kawanaka \cite{kawanaka} there exists a $\theta$-stable element $h\in H$ such that $x\in{\mathfrak g}(2;h)$ (see Sect. \ref{sec:5.2} for a more detailed account of Kawanaka's theorem).
But for any $\theta$-stable $h\in H$ there is some $m\in{\mathbb N}$ and a cocharacter $\lambda:k^\times\longrightarrow (T\cap K)$ such that $(\Ad\lambda(t))(e_\alpha)=t^{mh(\alpha)}e_\alpha$ for all $\alpha\in\Phi$.
Hence $0\in\overline{(\Ad\lambda(t))(x)}$.
%\qed
\end{proof}

This allows us to describe the closed $K$-orbits in ${\mathfrak p}$.

\begin{lemma}\label{closed}
Let $x\in{\mathfrak p}$ and let $x=x_s+x_n$ be the Jordan-Chevalley decomposition of $x$.
Then $K\cdot x_s$ is the unique closed (resp. minimal) orbit in $\overline{K\cdot x}$.
\end{lemma}

\begin{proof}
By standard results of geometric invariant theory there is a unique closed orbit in $\overline{K\cdot x}$, which is also the unique minimal orbit.
Let $y\in\overline{K\cdot x}$.
Clearly $y$ is in the minimal orbit if and only if $\dim Z_K(y)\geq \dim Z_K(y')$ for all $y'\in\overline{K\cdot x}$.
But by Lemmas \ref{centdim} and \ref{globalinf} this is true if and only if $\dim Z_G(y)\geq Z_G(y')$ for all $y'\in\overline{K\cdot x}$.
It is well-known that $G\cdot x_s$ is the unique closed orbit in $\overline{G\cdot x}$.
Thus $\dim Z_G(x_s)\geq\dim Z_G(y)$ for all $y\in\overline{G\cdot x}$.
It remains to show that $x_s\in\overline{K\cdot x}$.

Let $L=Z_G(x_s)^\circ$.
Then $L$ is a $\theta$-stable reductive group satisfying the standard hypotheses (A)-(C), and $x_s,x_n\in{\mathfrak l}=\Lie(L)={\mathfrak z}_{\mathfrak g}(x_s)$.
By Lemma \ref{unstable}, $x_n$ is $(K\cap L)^\circ$-unstable.
Hence $x_s\in\overline{(K\cap L)^\circ\cdot x}$.
Therefore $x_s\in\overline{K\cdot x}$.
This completes the proof.
%\qed
\end{proof}

\subsection{Chevalley Restriction Theorem}
\label{sec:4.5}

We now present a variant of the Chevalley Restriction Theorem.
The proof follows Richardson's proof of the corresponding result for the group $G$.
We begin with the following lemma, which is a direct analogue of \cite[11.1]{rich2}.
Fix a maximal $\theta$-split torus $A$ of $G$ with `baby Weyl group' $W=N_G(A)/Z_G(A)\cong N_K(A)/Z_K(A)$ (\cite[\S 4]{rich2}).
Let ${\mathfrak a}=\Lie(A)$.

\begin{lemma}\label{toriconj}
Suppose $\Ad g(Y)$ $\subseteq{\mathfrak a}$ for some $g\in K$.
Then there exists $w\in W$ such that $w\cdot y=\Ad g(y)$ for all $y\in Y$.
\end{lemma}

\begin{proof}
Let $L=Z_G(\Ad g(Y))^\circ$ and ${\mathfrak l}=\Lie(L)={\mathfrak z}_{\mathfrak g}(\Ad g(Y))$: $L$ is $\theta$-stable, reductive and satisfies the standard hypotheses (A)-(C) of \S 3.
Since ${\mathfrak a}\subseteq{\mathfrak l}$ and $\Ad g({\mathfrak a})\subseteq{\mathfrak l}$ there exists $l\in (K\cap L)^\circ$ such that $\Ad l(\Ad g({\mathfrak a}))={\mathfrak a}$.
Thus $n=(lg)\in N_K({\mathfrak a})=N_K(A)$ by Lemma \ref{maxsplittori}.
But $\Ad n(y)=\Ad g(y)$ for all $y\in Y$.
%\qed
\end{proof}

Since any finite set of points is closed, the set ${\mathfrak a}/W$ of $W$-orbits in ${\mathfrak a}$ has the structure of an affine variety with coordinate ring $k[{\mathfrak a}]^{W}$.

\begin{theorem}\label{Chev}
Let $A$ be a maximal $\theta$-split torus of $G$, and let $W=N_G(A)/Z_G(A)$.
Let ${\mathfrak a}=\Lie(A)$.
Then the natural embedding $j:{\mathfrak a}\longrightarrow{\mathfrak p}$ induces an isomorphism of affine varieties $j':{\mathfrak a}/W\longrightarrow{\mathfrak p}\quot K$.
Hence $k[{\mathfrak p}]^K$ is isomorphic to $k[{\mathfrak a}]^W$.
\end{theorem}

\begin{proof}
Let $\pi_{\mathfrak p}=\pi_{{\mathfrak p},K}:{\mathfrak p}\longrightarrow{\mathfrak p}\quot K$ and let $\pi_{\mathfrak a}=\pi_{{\mathfrak a},W}:{\mathfrak a}\longrightarrow{\mathfrak a}/W$.
Any $K$-invariant function on ${\mathfrak p}$ restricts to a $W$-invariant function on ${\mathfrak a}$.
Hence there is a well-defined $k$-algebra homomorphism from $k[{\mathfrak p}]^K$ to $k[{\mathfrak a}]^W$.
Taking the induced map on prime ideal spectra we have a morphism $j'$ making the following diagram commutative:

\begin{diagram}
 {\mathfrak a} & \rTo^j & {\mathfrak p} \\
 \dTo^{\pi_{\mathfrak a}} & & \dTo^{\pi_{\mathfrak p}} \\
 {\mathfrak a}/W & \rTo^{j'} & {\mathfrak p}\quot K
\end{diagram}

By a standard result of geometric invariant theory the varieties ${\mathfrak a}/W$ and ${\mathfrak p}\quot K$ are normal.
Thus by \cite[\S AG. 18.2]{bor} it will suffice to show that $j'$ is bijective and separable.
Recall that the points of ${\mathfrak p}\quot K$ correspond bijectively with the set of closed $K$-orbits in ${\mathfrak p}$.
Moreover by Lemma \ref{closed} the closed $K$-orbits in ${\mathfrak p}$ are precisely the semisimple orbits.
But by Thm. \ref{carts} any semisimple orbit meets ${\mathfrak a}$.
Hence $j'$ is surjective.
Let $a,a'\in{\mathfrak a}$ be such that $\pi_{\mathfrak p}(a)=\pi_{\mathfrak p}(a')$.
As $a,a'$ are semisimple they must be in the same $K$-orbit.
But by Lemma \ref{toriconj} this implies that $w\cdot a=a'$ for some $w\in W$.
Hence $\pi_{\mathfrak a}(a)=\pi_{\mathfrak a}(a')$.
Therefore $j'$ is injective.

It remains to show that $j'$ is separable.
As ${\mathfrak p}$ is irreducible and the set of regular semisimple elements is non-empty, the quotient morphism $\pi=\pi_{{\mathfrak p},K}$ is separable (\cite[9.3]{rich2}).
Moreover $\phi:K\times{\mathfrak a}\longrightarrow{\mathfrak p}$, $\phi(g,a)=\Ad g(a)$ is a separable morphism by Cor. \ref{sep}.
Thus $\pi\circ\phi:K\times{\mathfrak a}\longrightarrow{\mathfrak p}\quot K$ is separable.
We consider the action of $K$ on $K\times{\mathfrak a}$ in which $g'\cdot(g,a)=(g'g,a)$.
Since $\pi(\Ad g(A))=\pi(a)$, the composition $\pi\circ\phi$ factors through the action of $K$ on $K\times{\mathfrak a}$.
Note that ${\mathfrak p}\quot K$ can be thought of as a $K$-variety with the trivial action.
Hence there is a morphism $\sigma$ making the following diagram commutative:
\begin{diagram}
 K\times {\mathfrak a} & \rTo^{\pi\circ\phi} & {\mathfrak p}\quot K \\
 \dTo^{\pi_{K\times{\mathfrak a},K}} & & \dTo^{\Id_{{\mathfrak p}\quot K}} \\
 (K\times {\mathfrak a})\quot K & \rTo^{\sigma} & {\mathfrak p}\quot K
\end{diagram}

Since $\pi\circ\phi$ is separable, so is $\sigma$.
Let $i:{\mathfrak a}\longrightarrow K\times{\mathfrak a}$, $i(a)=(e,a)$.
Then it is easy to see that $\mu=\pi_{{K\times{\mathfrak a}},K}\circ i:{\mathfrak a}\rightarrow (K\times{\mathfrak a})\quot K$ is an isomorphism of varieties, hence that $\sigma\circ\mu:{\mathfrak a}\rightarrow{\mathfrak p}\quot K$ is separable.
But $\sigma\circ\mu=j'\circ\pi_{\mathfrak a}$.
Hence $j'$ is separable.
This completes the proof of the theorem.
%\qed
\end{proof}

Recall that $K^*=\{ g\in G\,|\;g^{-1}\theta(g)\in Z(G)\}$ normalizes $K$, and that $K^*= K\cdot F^*$, where $F^*=\{ a\in A| a^2\in Z(G)\}$ \cite[8.1]{rich2}.

\begin{corollary}\label{invext}
$k[{\mathfrak p}]^{K^*}=k[{\mathfrak p}]^K$.
\end{corollary}

\begin{proof}
Clearly $k[{\mathfrak p}]^{K^*}\subseteq k[{\mathfrak p}]^K$.
Hence we have to prove that any element of $k[{\mathfrak p}]^K$ is $K^*$-invariant.
As $K$ is normal in $K^*$, $K^*$ acts on $k[{\mathfrak p}]^K$.
Let $f\in k[{\mathfrak p}]^K$.
To show that $f\in k[{\mathfrak p}]^{K^*}$ it will suffice to show that $a\cdot f=f$ for any $a\in F^*$.
But $(a\cdot f)(x)=f(a^{-1}\cdot x)=f(x)$ for all $x\in{\mathfrak a}$, hence $(j')^*(a\cdot f)=(j')^*(f)$.
Taking inverses under $(j')^*$, we see that $a\cdot f=f$.
Thus $f\in k[{\mathfrak p}]^{K^*}$.
%\qed
\end{proof}

\subsection{$k[{\mathfrak p}]^K$ is a polynomial ring}
\label{sec:4.6}

Let ${\mathfrak g}$ be a complex semisimple Lie algebra, let ${\mathfrak t}$ be a Cartan subspace of ${\mathfrak g}$, and let $W$ be the Weyl group of ${\mathfrak g}$ acting on ${\mathfrak t}$.
It is well-known from the classical theory that the algebra of invariants ${\mathbb C}[{\mathfrak t}]^W$ is generated by $r=\dim{\mathfrak t}$ algebraically independent homogeneous generators of degrees $(m_1+1),(m_2+1),\ldots,(m_r+1)$, where the $m_i$ are the {\it exponents} of ${\mathfrak g}$.

We will now show that an analogous statement is true for ${\mathfrak a}$.
It is a straightforward application of Demazure's theorem \cite{dem} on Weyl group invariants.

\begin{lemma}\label{demaz}
Let $A$ be a maximal $\theta$-split torus of $G$ and let ${\mathfrak a}=\Lie(A)$.
Let $W=N_G({\mathfrak a})/Z_G({\mathfrak a})=N_G(A)/Z_G(A)$.
There are $r$ algebraically independent homogeneous polynomials $f_1,f_2,\ldots,f_r$ (where $r$ is the rank of ${\mathfrak a}$) such that $k[{\mathfrak a}]^W=k[f_1,f_2,\ldots,f_r]$.
Moreover
$$\sum_{w\in W} t^{l(w)}=\prod_{i=1}^r{\frac{1-t^{\deg f_i}}{1-t}}$$
where $l$ is the length function on $W$ corresponding to a basis of simple roots in $\Phi_A$.
\end{lemma}

\begin{proof}
Let $T$ be a torus of rank $n$ and let ${\mathfrak t}=\Lie(T)$.
The character group $X(T)$ is a free abelian group of rank $n$.
There is a natural isomorphism $X(T)\otimes_{\mathbb Z} k\rightarrow{\mathfrak t}^*$ induced by the map $\alpha\otimes 1\mapsto d\alpha$, which is equivariant with respect to any group $H$ of automorphisms of $T$.
Hence $k[{\mathfrak t}]^H\cong S(X(T)\otimes_{\mathbb Z} k)^H$.
In particular, $k[{\mathfrak a}]^{W}\cong S(X(A)\otimes_{\mathbb Z} k)^{W}$.

We recall that, according to Demazure's definition, a reduced root system is a triple ${\cal R}=(M,R,\rho)$, where $M$ is a free ${\mathbb Z}$-module of finite type, $R$ is a subset of $M$, and $\rho:\alpha\mapsto \alpha^\vee$ is a map from $R$ into the dual $M^*$ of $M$ such that:

(a) $R$ is finite and $R\cap(2R)=\emptyset$,

(b) For every $\alpha\in R, \alpha^\vee(\alpha)=2$,

(c) If $\alpha,\beta\in R$, then $\beta-\alpha^\vee(\beta)\alpha\in R$, and $\beta^\vee-\beta^\vee(\alpha)\alpha^\vee\in R^\vee$, where $R^\vee=\rho(R)$.

Let $\Phi_A^*$ be the subset of $\Phi_A$ consisting of all roots $\alpha$ such that $\alpha/m\in\Phi_A\Rightarrow m=\pm 1$.
By \cite[\S 4]{rich2} there exists a map $\rho$ such that $(X^*(A),\Phi_A^*,\rho)$ is a root system in this sense.
Moreover by \cite[4.3]{rich2}, $W$ is generated by the reflections $s_\alpha$ with $\alpha\in\Phi_A^*$.
Finally, by Lemma \ref{pisgood}, $p$ is good for $\Phi_A^*$.
Hence by \cite[Cor. to Thm. 2, Thm. 3]{dem} $S(X_*(A)\otimes k)^W$ is generated by $r$ algebraically independent homogeneous polynomials, of degrees $d_1,d_2,\ldots,d_r$ such that
$$\sum_{w\in W}t^{l(w)}=\prod_{i=1}^r \frac{1-t^{d_i}}{1-t}.$$
%\qed
\end{proof}

We remark that the product $\prod_{i=1}^r \frac{1-t^{d_i}}{1-t}$ here may include a number of factors of the form $(1-t)/(1-t)=1$.

\section{The nilpotent cone}
\label{sec:5}

\subsection{Equidimensionality}
\label{sec:5.1}

Let ${\cal N}={\cal N}({\mathfrak p})$ be the set of nilpotent elements of ${\mathfrak p}$.
In general ${\cal N}$ is not irreducible (see for example Cor. \ref{splitcmpts}).
However, we have the following straightforward result (Thm. 3 in \cite{kostrall}).
We include the proof, which is similar to Kostant-Rallis', for the convenience of the reader.

\begin{theorem}\label{nil1}
Let ${\cal N}$ be the affine variety of all nilpotent elements in ${\mathfrak p}$, and let ${\cal N}_1,{\cal N}_2,\ldots,{\cal N}_m$ be the irreducible components of ${\cal N}$.
The number of $K$-orbits in ${\cal N}$ is finite.
For each $i$, $\codim_{\mathfrak p}{\cal N}_i = r = \rank A$, where $A$ is a maximal $\theta$-split torus of $G$.
Moreover, $K$ normalizes ${\cal N}_i$, and there is an open $K$-orbit in ${\cal N}_i$.
An element of ${\cal N}_i$ is in the open $K$-orbit if and only if it is regular.
\end{theorem}

\begin{proof}
Let $e\in{\cal N}$.
Then $g\cdot e\in{\cal N}$ for any $g\in K$ (in fact for any $g\in K^*$).
Hence $K$ normalizes ${\cal N}$.
But $K$ is connected, therefore $K\cdot {\cal N}_i={\cal N}_i$ for each $i$.
By \cite[Thm. D]{rich3} there are finitely many $K$-orbits in ${\cal N}$.
Hence each irreducible component of ${\cal N}$ contains a unique open orbit.
If $x\in{\mathfrak p}$, then $\codim_{\mathfrak p} (K\cdot x)\geq r$ by Lemma \ref{regs}.
Therefore $\codim_{\mathfrak p}{\cal N}_i\geq r$.
But by Lemmas \ref{unstable} and \ref{demaz}, ${\cal N}$ is the set of zeros of $r$ homogeneous polynomials $u_1,u_2,\ldots,u_r$, where $k[{\mathfrak p}]^K=k[u_1,u_2,\ldots,u_r]$.
Therefore $\codim_{\mathfrak p}{\cal N}_i\leq r$.
The remaining statements follow at once.
%\qed
\end{proof}

\subsection{Kawanaka's Theorem}
\label{sec:5.2}

In \cite{kawanaka}, Kawanaka generalised the Bala-Carter theory to classify nilpotent orbits in eigenspaces for automorphisms of semisimple Lie algebras.
We now recall Kawanaka's theorem as it applies to the case of an involution.
Let $(B,T)$ be a fundamental pair for $\theta$, let $\Delta$ be the basis of the roots $\Phi=\Phi(G,T)$ corresponding to $B$, and let $W_T=W(G,T)$ be the Weyl group.
Let $\Lambda_r$ be the root lattice of $\Phi$ and let $H=H(\Phi,\Delta)$ be the abelian group of all homomorphisms from $\Lambda_r$ to ${\mathbb Z}$.
An element $h\in H$ is uniquely determined by the values $h(\alpha_i)$ for $\alpha_i\in\Delta$.
Hence we may describe an element of $H$ by means of a copy of the Dynkin diagram on $\Delta$ with weights attached to each node.

Let $X(T)=\Hom(T,k^\times)$ and let $Y(T)=\Hom(k^\times,T)$.
Denote by $\langle .\, ,.\rangle:X(T)\times Y(T)\longrightarrow{\mathbb Z}$ the natural $W$-equivariant, ${\mathbb Z}$-bilinear map.
Hence $\alpha(\lambda(t))=t^{\langle\alpha,\lambda\rangle}$ for all $t\in k^\times$.
The pairing induces a homomorphism $Y(T)\rightarrow H$.
We denote by $\overline\lambda$ the element of $H$ corresponding to $\lambda\in Y(T)$.
Hence $\overline\lambda(\alpha)=\langle\alpha,\lambda\rangle$ for all $\alpha\in\Phi$.
The image of $Y(T)$ is of finite index in $H$.
Thus, for any $h\in H$ there exists a positive integer $m$ and a cocharacter $\lambda$ such that $\overline\lambda=mh$.

Let $H^+$ be the positive Weyl chamber associated to $\Delta$: $h\in H^+\Leftrightarrow h(\alpha_i)\geq 0\;\forall \,\alpha_i\in\Delta$.
The Weyl group $W_T$ acts naturally on $H$, and $\overline{w(\lambda)}=w(\overline\lambda)$ for any $\lambda\in Y(T)$.
For any $h\in H$ there exists $w\in W_T$ and $h_+\in H^+$ such that $w(h)=h_+$.
Moreover, $h_+$ is unique.
For $h\in H$, let ${\mathfrak g}(i;h)=\sum_{h(\alpha)=i}{\mathfrak g}_\alpha$, $i\neq 0$, and ${\mathfrak g}(0;h)={\mathfrak t}\oplus\sum_{h(\alpha)=0}{\mathfrak g}_\alpha$.
The decomposition ${\mathfrak g}=\oplus{\mathfrak g}(i;h)$ is a ${\mathbb Z}$-grading of ${\mathfrak g}$, and the $\overline\lambda$-grading coincides with the $(\Ad\lambda)$-grading for $\lambda\in Y(T)$.

If $k={\mathbb C}$, there is a straightforward classification of nilpotent orbits via conjugacy classes of $\mathfrak{sl}(2)$-triples:
any nilpotent element $e\in{\mathfrak g}$ can be embedded as the nilpositive element of an $\mathfrak{sl}(2)$-triple $\{ h,e,f\}$; moreover, there is a unique $G$-conjugate $h'$ of $h$ such that $h'\in{\mathfrak t}$ and $\alpha(h')\geq 0$ for all $\alpha\in\Delta$.
(It was proved by Dynkin that $\alpha(h')\in\{ 0,1,2\}$ for all $\alpha\in\Delta$.)
In this way one can associate to $e$ a unique element of $H(\Phi,\Delta)^+$, called the {\it weighted Dynkin diagram} associated to $e$.
We denote the set of all weighted Dynkin diagrams by $H(\Phi,\Delta)_n$.
Hence there is a one-to-one correspondence between the elements of $H(\Phi,\Delta)_n$ and the nilpotent conjugacy classes in ${\mathfrak g}$.

This argument using $\mathfrak{sl}(2)$-triples is only valid if the characteristic is zero or large.
However, Pommerening proved in \cite{pom1,pom2} that the nilpotent orbit structure is essentially the same in all good characteristics.
Let $h\in H(\Phi,\Delta)_n$ and let $G(0)_h$ be the unique closed connected subgroup of $G$ such that $\Lie(G(0)_h)={\mathfrak g}(0;h)$.
There is an open $G(0)_h$-orbit in ${\mathfrak g}(2;h)$: let $N_h$ be a representative for the open orbit and set ${\mathfrak o}_h=G\cdot N_h$.
The correspondence $h\mapsto{\mathfrak o}_h$ is one-to-one between the elements of $H(\Phi,\Delta)_n$ and nilpotent conjugacy classes in ${\mathfrak g}$.

In good characteristic Pommerening replaced weighted Dynkin diagrams with {\it associated characters}.
A cocharacter $\lambda$ is associated to $e$ if $e\in{\mathfrak g}(2;\lambda)$ and there is a Levi subgroup $L$ of $G$ such that $\lambda(k^\times)\subset L^{(1)}$ and $e$ is distinguished in $\Lie(L)$.
(A nilpotent element $x\in{\mathfrak g}$ is distinguished if $Z_{G^{(1)}}(x)^\circ$ is a unipotent group.)
If $\lambda$ is an associated cocharacter for $e$ and $g\in Z_G(e)$, then $g\cdot\lambda$ is also associated to $e$; moreover, any two associated cocharacters for $e$ are conjugate by an element of $Z_G(e)^\circ$ (\cite[Prop. 11]{mcninch3}).

Premet has recently given a short conceptual proof of Pommerening's theorem, valid in all good characteristics.
The proof uses the theory of optimal cocharacters for $G$-unstable elements, also called the Kempf-Rousseau theory.
If $\rho:G\longrightarrow\GL(V)$ is a rational representation, then the Kempf-Rousseau theory attaches to a $G$-unstable vector $v\in V$ a collection of {\it optimal} cocharacters.
In general the optimal cocharacters depend on the choice of a length function on the set of cocharacters in $G$.
(See Sect. \ref{sec:6.2} for the details concerning the Kempf-Rousseau theory.)
Let $h\in H(\Phi,\Delta)_n$.
As observed in \cite[\S 2.4]{premnil}, there exists a (unique) cocharacter $\lambda:k^\times\longrightarrow T\cap G^{(1)}$ such that $\overline\lambda=h$.
(Since this holds for simply-connected $G^{(1)}$, it holds for any isogenous image of $G^{(1)}$, hence for arbitrary reductive groups.)
Let $U$ be the unique closed connected $T$-stable subgroup of $G$ such that $\Lie(U)=\sum_{i>0}{\mathfrak g}(i;h)$ and let $P=P(\lambda)=Z_G(\lambda)\cdot U$ (a parabolic subgroup of $G$).
Then, after choosing a suitable length function on the set of cocharacters in $G$, we have (see \cite[Thm. 2.3]{premnil} and \cite[3.5]{mcninch3}):

\begin{theorem}[Premet]\label{premorbits}
(a) $\lambda$ is optimal for $N_h$.

(b) Let $C=Z_G(e)\cap Z_G(\lambda)$.
Then $Z_G(e)\subset P$ and $Z_G(e)=C\cdot Z_U(e)$ (semidirect product): $C$ is the reductive part and $Z_U(e)$ the unipotent radical of $Z_G(e)$.

(c) Let $S$ be a maximal torus of $C$ and let $L=Z_G(S)$.
Then $e$ is a distinguished nilpotent element of $\Lie(L)$ and $\lambda(k^\times)\subset L^{(1)}$.
\end{theorem}

Note that by (c) $\lambda$ is associated to $N_h$.
It follows that the decomposition in (b) holds for arbitrary $e,\lambda$, where $e$ is nilpotent and $\lambda$ is associated to $e$.
We wish to restate Kawanaka's theorem (for the case of an involution) in the language of associated cocharacters.
Let $h\in H$ be $\theta$-stable.
Define a subalgebra $\overline{\mathfrak g}_h$ of ${\mathfrak g}$ with graded components $\overline{\mathfrak g}_h(i)$ as follows:
$\overline{\mathfrak g}_h(i)=
\left\{
\begin{array}{ll}
{\mathfrak k}(i;h) & \mbox{if $i = 0$ (mod 4),} \\
{\mathfrak p}(i;h) & \mbox{if $i = 2$ (mod 4),} \\
\{ 0\} & \mbox{otherwise.}
\end{array}
\right.$

Suppose further that $h_+\in H(\Phi,\Delta)_n$.
Since $h$ is $W$-conjugate to $h_+$, there exists a unique cocharacter $\lambda:k^\times\longrightarrow T\cap G^{(1)}$ such that $\overline\lambda=h$.
But $\theta(h)=h$, hence $\lambda(k^\times)\subset T\cap K\cap G^{(1)}$.
The Lie algebra $\overline{\mathfrak g}_h$ is equal to ${\mathfrak g}^{d\psi}=\{ x\in{\mathfrak g}\,|\,d\psi(x)=x\}$, where $t_0=\lambda(\sqrt{-1})$ and $\psi=\Int t_0\circ\theta$.
Moreover, $\psi$ is of order 1,2 or 4, hence is semisimple.
It follows that $\overline{\mathfrak g}_h=\Lie((G^\psi)^\circ)$ and $\overline{G}_h=(G^\psi)^\circ$ is reductive.
(This is true for any $\theta$-stable $h\in H$, see \cite{kawanaka}.)
Let $\overline{G}_h(0)=Z_K(\lambda)$.
Then $T(0)=(T\cap K)^\circ$ is a maximal torus of $\overline{G}_h$, and $\Lie(\overline{G}_h(0))={\mathfrak k}^\lambda=\overline{\mathfrak g}_h(0)$.
Following Kawanaka, $h$ is {\it slim} (with respect to $\theta$) if $\lambda(k^\times)\subset \overline{G}_h^{(1)}$.

Let $\alpha\in\Phi=\Phi(G,T)$.
Recall that $\theta$ induces an automorphism $\gamma$ of $\Phi$ stabilizing $\Delta$.
Denote by ${\mathfrak g}_{(\alpha)}$ the span of the root spaces ${\mathfrak g}_\alpha$ and ${\mathfrak g}_{\gamma(\alpha)}$.
If $\gamma(\alpha)\neq\alpha$, then ${\mathfrak g}_{(\alpha)}=({\mathfrak g}_{(\alpha)}\cap{\mathfrak k})\oplus({\mathfrak g}_{(\alpha)}\cap{\mathfrak p})$ and the dimension of each summand is 1.
Let $\overline{\alpha}$ denote the restriction of $\alpha$ to $T(0)$.
Note that $\overline\alpha=\overline\beta$ if and only if $\beta\in\{\alpha,\gamma(\alpha)\}$.
We have $\Phi_h=\Phi(\overline{G}_h,T(0))=\{\overline\alpha\,|\,{\mathfrak g}_{(\alpha)}\cap\overline{\mathfrak g}_h\neq\{ 0\}\}$.

Let $\alpha\in\Phi$.
There are three possibilities: (i) $\gamma(\alpha)=\alpha$, (ii) $\gamma(\alpha)$ and $\alpha$ are orthogonal, (iii) $\gamma(\alpha)$ and $\alpha$ generate a root system of type $A_2$.
Introduce corresponding elements $s_{(\alpha)}$ of $W$: (i) $s_{(\alpha)}=s_\alpha$, (ii) $s_{(\alpha)}=s_\alpha s_{\gamma(\alpha)}$, (iii) $s_{(\alpha)}=s_\alpha s_{\gamma(\alpha)} s_\alpha=s_{\gamma(\alpha)}s_\alpha s_{\gamma(\alpha)}=s_{\alpha+\gamma(\alpha)}$.
We can embed the Weyl group $W_h=W(\Phi_h)$ in $W$: $W_h$ is generated by all $s_{(\alpha)}$ with $\overline\alpha\in\Phi_h$.
Let $\Phi^+$ be the positive system in $\Phi$ determined by $\Delta$ and let $\Phi_h^+=\{\overline\alpha\in\Phi_h\,|\,\alpha\in\Phi^+\}$.
Then $\Phi_h^+$ is a positive system in $\Phi_h$.
We let $\Delta_h$ be the corresponding basis.
Any $\theta$-stable element $h'$ of $H(\Phi,\Delta)$ gives rise to a well-defined element $\overline{h'}$ of $H(\Phi_h,\Delta_h)$.

Kawanaka introduced a subset $H(\Phi,\Delta,\theta)'_n$ of $H$ in order to parametrise the nilpotent $K$-orbits in ${\mathfrak p}$: $h\in H(\Phi,\Delta,\theta)'_n$ if and only if:

(i) $h_+\in H(\Phi,\Delta)_n$,

(ii) $h$ is $\theta$-invariant,

(iii) $h$ is slim with respect to $\theta$,

(iv) $\overline{h}_+\in H(\Phi_h,\Delta_h)_n$.

Let $W(0)=N_K(T)/Z_K(T)$ and let $W^\theta=\{ w\in W|\theta(w)=w\}$.
Let $H(\Phi,\Delta,\theta)_n$ be a set of representatives for the $W(0)$-orbits in $H(\Phi,\Delta,\theta)'_n$.

Kawanaka's theorem states that \cite[(3.1.5)]{kawanaka}:

\begin{theorem}[Kawanaka]
For each $h\in H(\Phi,\Delta,\theta)_n$ choose a representative $N_h$ of the open $\overline{G}_h(0)$-orbit in $\overline{\mathfrak g}_h(2)={\mathfrak p}(2;h)$.
Then the correspondence $h\mapsto K\cdot N_h$ is one-to-one between elements of $H(\Phi,\Delta,\theta)_n$ and nilpotent $K$-orbits in ${\mathfrak p}$.
We have $K\cdot N_h\subset{\mathfrak o}_{h_+}$, the $G$-orbit determined by $h_+$.
Two orbits $K\cdot N_h$ and $K\cdot N_{h'}$ are contained in the same $G$-orbit if and only if $h_+=h'_+$.
\end{theorem}

(Kawanaka's theorem is stated in a much more general setting, which includes the case of an automorphism of $G$ of finite order prime to $p$.)
In view of the remarks above, we have the following:

\begin{corollary}\label{assoc}
Let $e\in{\cal N}$.
Then there exists a cocharacter $\lambda:k^\times\longrightarrow K$ which is associated to $e$.
Any two such cocharacters are conjugate by an element of $Z_K(e)^\circ$.
\end{corollary}

\begin{proof}
By Kawanaka's theorem there exists $g\in K$ and $h\in H(\Phi,\Delta,\theta)_n$ such that $g\cdot e=N_h$.
But as we have already seen, there exists a unique cocharacter $\lambda:k^\times\longrightarrow T\cap K\cap G^{(1)}$ such that $\overline\lambda=h$.
Moreover, $\lambda$ is associated to $N_h$.
It follows that $g^{-1}\cdot\lambda$ is associated to $e$.

Suppose $\lambda,\mu$ are associated cocharacters for $e$ such that $\lambda(k^\times),\mu(k^\times)\subset K$.
There exists $g\in Z_G(e)^\circ$ such that $g\cdot\lambda=\mu$ (\cite[Prop. 11]{mcninch3}).
Let $C=Z_G(e)\cap Z_G(\lambda)$: then $Z_G(e)^\circ=C^\circ\cdot Z_U(e)$ (semidirect product), where $Z_U(e)$ is the unipotent radical of $Z_G(e)$.
Hence there exists $u\in Z_U(e)$ such that $u\cdot\lambda=\mu$.
Since $e\in{\mathfrak p}$, $Z_U(e)$ is $\theta$-stable.
But now $u^{-1}\theta(u)\in Z_G(\lambda)\cap Z_U(e)\;\Rightarrow\; u^{-1}\theta(u)=1$.
By \cite[III.3.12]{sands} $u\in Z_K(e)^\circ$.
%\qed
\end{proof}

This observation allows us to replace the notion of weighted Dynkin diagrams with that of associated cocharacters.
If $e\in{\cal N}$ and $\lambda:k^\times\longrightarrow K$ is an associated cocharacter for $e$, we use the notation $\overline{G}_\lambda= (G^{\psi})^\circ,\overline{\mathfrak g}_\lambda=\Lie(\overline{G}_\lambda)$, where $\psi=\Int\lambda(\sqrt{-1})\circ\theta$.

\begin{rk}\label{nonsc}
The theorems of Kawanaka, Pommerening and Premet are true for arbitrary reductive $G$ such that $p$ is good.
Hence Cor. \ref{assoc} is true without the assumptions (B),(C) of \S 3.
If we assume only that $p$ is good for $G$, then we can define $x\in{\mathfrak p}$ to be regular if $\dim Z_G(x)$ is minimal: then $\dim G-\dim Z_G(x)=\dim{\mathfrak g}^A$ by Lemma \ref{maxsplittori} and \cite[3.2]{rich2}.
(We don't in general have $\dim {\mathfrak k}-\dim {\mathfrak p}
=\dim {\mathfrak z}_{\mathfrak k}(x)-\dim {\mathfrak z}_{\mathfrak p}(x))$ for all $x\in{\mathfrak p}$.)
Let $G$ be simply-connected and semisimple and let $\tilde{G}$ be the group defined in \S 3.
Then we can lift an involution of $G$ to $\tilde{G}$ by Lemma \ref{GLautos}.
Hence Thm. \ref{nil1} is true for any semisimple simply-connected group.
Let $G$ be an arbitrary semisimple group and let $\pi:G_{sc}\rightarrow G$ be the universal cover of $G$.
Then by the argument in \cite[2.3]{premnil} $\pi$ induces a $G/Z(G)$-equivariant bijection ${\cal N}({\mathfrak g}_{sc})\longrightarrow{\cal N}({\mathfrak g})$.
Moreover, any involutive automorphism of $G$ can be lifted to an involutive automorphism of $G_{sc}$.
It follows that Thm. \ref{nil1} holds for any semisimple group with involution (assuming $p$ is good).
Note that if $p$ is good for $G$ then it is good for $\overline{G}_\lambda$.
(This is immediate since $p\neq 2$, therefore $p$ can only be bad for $\overline{G}_\lambda$ if it is of exceptional type: but if $\overline{G}_\lambda$ is of exceptional type then so is $G$, and the semisimple rank of $G$ is greater than that of $\overline{G}_\lambda$.)
\end{rk}

\subsection{Semiregular Elements in Type $D_n$}
\label{sec:5.3}

Let $G$ be almost simple, simply-connected of type $D_n$, let $T$ be a maximal torus of $G$ and let $\Delta=\{\alpha_1,\alpha_2,\ldots,\alpha_n\}$ be a basis for $\Phi=\Phi(G,T)$, numbered in the standard way.
Let ${\mathfrak g}=\Lie(G)$ and let $\{ h_{\alpha_i},e_\alpha\,|\,\alpha_i\in\Delta,\alpha\in\Phi\}$ be a Chevalley basis for ${\mathfrak g}$.
Let $\gamma$ be the graph automorphism which sends $\alpha_{n-1}\mapsto\alpha_n$, $\alpha_n\mapsto\alpha_{n-1}$, and fixes all other elements of $\Delta$.
The following lemma is due to Premet.

\begin{lemma}\label{sigmaexists}
There exists an automorphism $\sigma$ of $G$ satisfying $d\sigma(e_\alpha)=e_{\gamma(\alpha)}$ for all $\alpha\in\Phi$.
\end{lemma}

\begin{proof}
Since any automorphism of ${\mathfrak g}$ gives rise to an automorphism of the adjoint group, and hence by Lemma \ref{sccover} to an automorphism of $G$, it will suffice to show that there is an automorphism of ${\mathfrak g}$ satisfying $e_\alpha\mapsto e_{\gamma(\alpha)}$ for all $\alpha\in\Phi$.
Let $\phi$ be the (unique) automorphism of ${\mathfrak g}$ which sends $e_\alpha$ to $e_{\gamma(\alpha)}$ for $\alpha\in\pm\Delta$.
Let $I=\{\alpha_1,\alpha_2,\ldots,\alpha_{n-2}\}$ and let $\Phi_I$ be the subsystem of $\Phi$ generated by the elements of $I$.
It is easy to see that $\phi(e_\alpha)=e_\alpha$ for any $\alpha\in\Phi_I$.

Let $\alpha\in\Phi^+\setminus\Phi_I$.
There are four possibilities:

(i) $\alpha=\beta+\alpha_{n-1}$ for some $\beta\in\Phi_I$,

(ii) $\alpha=\beta+\alpha_n$ for some $\beta\in\Phi_I$,

(iii) $\alpha=\beta+\alpha_{n-1}+\alpha_n$ for some $\beta\in\Phi_I^+$,

(iv) $\alpha=(\beta+\alpha_{n-1})+(\gamma+\alpha_n)$ for some $\beta,\gamma\in\Phi_I^+$ with $\beta+\alpha_{n-1},\gamma+\alpha_n\in\Phi$.

For case (i), $e_\alpha=[e_\beta,e_{\alpha_{n-1}}]\mapsto [e_\beta,e_{\alpha_n}]=e_{\gamma(\alpha)}$.
Similarly for case (ii).
For (iii), $e_\alpha=[[e_\beta,e_{\alpha_{n-1}}],e_{\alpha_n}]= [[e_\beta,e_{\alpha_n}],e_{\alpha_{n-1}}]$.
Hence $\phi(e_\alpha)=e_\alpha=e_{\gamma(\alpha)}$.
Finally, if (iv) holds then $e_\alpha=\pm [e_{\beta+\alpha_{n-1}+\alpha_n},e_\gamma]$.
But $\phi(e_{\beta+\alpha_{n-1}+\alpha_n})=e_{\beta+\alpha_{n-1}+\alpha_n}$ and $\phi(e_\gamma)=e_\gamma$, by the above.
Hence $\phi(e_\alpha)=e_\alpha$.

We have proved that $\phi(e_\alpha)=e_{\gamma(\alpha)}$ for any $\alpha\in\Phi^+$.
But then by properties of the Chevalley basis $\phi(e_\alpha)=e_{\gamma(\alpha)}$ for any $\alpha\in\Phi$.
%\qed
\end{proof}

\begin{rk}
The existence of $\sigma$ clearly also holds if $G$ is of adjoint type.
However, if $n$ is even and $G$ is intermediate (that is, neither simply-connected nor adjoint) then $\sigma$ does not in general exist.
\end{rk}

Recall that a nilpotent element $e\in{\mathfrak g}$ is {\it distinguished} if $Z_{G^{(1)}}(e)^\circ$ is a unipotent group, and $e$ is {\it semiregular} if $Z_G(e)$ is the product of $Z(G)$ and a (connected) unipotent group.
(Hence a semiregular element is distinguished.)
Let $h\in H(\Phi,\Delta)_n$ be the weighted Dynkin diagram corresponding to a semiregular orbit, and let $\lambda:k^\times\longrightarrow T$ be the unique cocharacter satisfying $\langle\alpha,\lambda\rangle=h(\alpha)$ for all $\alpha\in\Phi$ (this exists by \cite[2.4]{premnil}).
Let $Y_\lambda$ be the open $Z_G(\lambda)$-orbit in ${\mathfrak g}(2;\lambda)$ and let $E\in Y_\lambda$.
It follows from \cite[III.4.28(ii)]{sands} that $\sigma(\lambda(t))=\lambda(t)$, and that $E$ is $Z_G(\lambda)$-conjugate to an element of the form $\sum_{\beta\in\Gamma}e_\gamma$, where $\Gamma$ is a $\gamma$-stable subset of $\{\alpha\in\Phi\,|\,h(\alpha)=2\}$.

Hence:

\begin{lemma}\label{typed}
Let $e$ be a semiregular nilpotent element of ${\mathfrak g}$ and let $\mu$ be an associated cocharacter for $G$.
After conjugating $e$ and $\mu$ by an element of $G$, if necessary, we may assume that $\mu(k^\times)\subset G^\sigma$ and $e\in{\mathfrak g}^\sigma$.
\end{lemma}

We also record the following result to be used in the next subsection.

\begin{lemma}\label{semiregred}
Let $G$ be any reductive group such that $p$ is good for $G$, and let $e$ be a distinguished nilpotent element of ${\mathfrak g}$.
Then there exists a reductive subgroup $L$ of $G$ such that (i) $e$ is a semiregular element of $\Lie(L)$, (ii) $p$ is good for $L$.
\end{lemma}

\begin{proof}
For any $x\in{\mathfrak g}$, $Z_G(x)=Z(G)\cdot Z_{G^{(1)}}(x)$.
Moreover, $e\in\Lie(G^{(1)})$.
Hence, after replacing $G$ by $G^{(1)}$, we may assume that $G$ is semisimple.
We now prove the lemma by induction on the order of the group $\overline{A}(e)=Z_G(e)/Z(G)Z_G(e)^\circ$.
If $\overline{A}(e)$ is trivial, then we are done.
Otherwise, let $x$ be any element of $Z_G(e)\setminus Z(G)Z_G(e)^\circ$, and let $x=x_s x_u$ be the Jordan-Chevalley decomposition of $x$.
Then $x_u\in Z_G(e)^\circ$, hence after replacing $x$ by $x_s$ we may assume that $x$ is semisimple.
Let $L'=Z_G(x)^\circ$.
Then $L'$ is a pseudo-Levi subgroup of $G$, hence is reductive and $p$ is good for $L'$.
Since $e$ is distinguished in $G$ (hence in $L'$), $Z(L')^\circ$ is trivial and $Z_{L'}(e)^\circ=(Z_G(e)^\circ)^x$.
Let $H=Z_{L'}(e)/Z(G)Z_{L'}(e)^\circ$ and let $\overline{A}_{L'}(e)=Z(L')(e)/Z(L')Z_{L'}(e)^\circ$.
Then $H\hookrightarrow Z_G(e)^x/Z(G)(Z_G(e)^\circ)^x$, hence $H$ can be considered as a subgroup of $\overline{A}(e)$.
Moreover, $H$ maps surjectively onto $\overline{A}_{L'}(e)$, and the kernel is non-trivial; thus the order of $\overline{A}_{L'}(e)$ is strictly less than that of $\overline{A}(e)$.
By the induction hypothesis, there exists a subgroup $L$ of $L'$ satisfying the conditions of the Lemma.
%\qed
\end{proof}

\subsection{Regular Nilpotent Elements}
\label{sec:5.4}

Our goal is to prove that the regular nilpotent elements form a single $K^*$-orbit, where $K^*=\{ g\in G|\,g^{-1}\theta(g)\in Z(G)\}$.
The following lemma is the key step.
In view of Remark \ref{nonsc}, we assume until further notice only that $p$ is good for $G$.
We use Bourbaki's numbering conventions on roots \cite{bourbaki}.

\begin{lemma}\label{splitconj}
Let $e$ be a nilpotent element of ${\mathfrak p}$ and let $\lambda:k^\times\longrightarrow K$ be associated to $e$.
Then there exists $g\in G$ such that $(\Int g)\circ\lambda$ is $\theta$-split.
Equivalently $\Int n(\lambda)=-\lambda$ where $n=g^{-1}\theta(g)$.
\end{lemma}

\begin{proof}
Recall that if $p$ is good for $G$ then it is good for $\overline{G}_\lambda$ (resp. a pseudo-Levi subgroup of $G$).
Hence, after replacing $G$ by $\overline{G}_\lambda$, we have only to prove the lemma under the assumption that $\theta=\Int\lambda(\sqrt{-1})$ and that all weights of $\lambda$ on ${\mathfrak g}$ are even.
Let $S$ be a maximal torus of $Z_G(\lambda)\cap Z_G(e)$.
Then $Z_G(S)$ is a $\theta$-stable Levi subgroup of $G$ and $e$ is a distinguished element of $Z_G(S)$ (\cite[Prop. 2.5]{premnil}).
Hence, after replacing $G$ by $Z_G(S)$, we may assume that $e$ is distinguished.
Let $L$ be a reductive subgroup of $G$ such that $p$ is good for $L$ and $e$ is a semiregular element of $\Lie(L)$ (Lemma \ref{semiregred}).
Let $\mu$ be an associated cocharacter for $e$ in $L$: then $\mu$ is also an associated cocharacter for $e$ in $G$.
Hence $\mu$ is $Z_G(e)$-conjugate to $\lambda$.
Conjugating $L$ by some element of $Z_G(e)$, if necessary, we may assume that $\lambda(k^\times)\subset L$.
It is well-known that $e\in\Lie(L^{(1)})$ (see \cite[\S 2.3]{premnil}, for example).
Replacing $G$ by $L^{(1)}$, we may assume that $G$ is semisimple and that $e$ is semiregular in ${\mathfrak g}$.

Now if $\eta:G_{sc}\rightarrow G$ is the universal covering, then by Lemma \ref{sccover} there exists a unique involutive automorphism $\theta_{sc}$ of $G_{sc}$ which lifts $\theta$.
By \cite[Rk. 1]{premnil} there is a (unique) cocharacter $\lambda_{sc}$ such that $\eta\circ\lambda_{sc}=\lambda$.
Hence $\theta_{sc}=\Int\lambda_{sc}(\sqrt{-1})$.
To prove that $\lambda$ is $G$-conjugate to a $\theta$-split cocharacter, it will clearly suffice to prove that $\lambda_{sc}$ is $G_{sc}$-conjugate to a $\theta$-split cocharacter.
Note that the statement of the Lemma does not depend on the choice of $e$: let $e_{sc}$ be any representative for the open $Z_{G_{sc}}(\lambda_{sc})$-orbit in ${\mathfrak g}_{sc}(2;\lambda_{sc})$.
Replacing $G,\lambda$, and $e$ respectively by $G_{sc},\lambda_{sc}$ and $e_{sc}$, we may assume that $G$ is semisimple and simply-connected, and that $e$ is semiregular in ${\mathfrak g}$.
Finally, let $G_1,G_2,\ldots,G_l$ be the minimal normal subgroups of $G$ and let ${\mathfrak g}_i=\Lie(G_i)$, $1\leq i\leq l$.
There is a unique expression $e=\sum e_i$, where each $e_i\in{\mathfrak g}_i$; thus $e_i$ is semiregular in ${\mathfrak g}_i$.
Moreover $\theta$ is inner, hence each component $G_i$ is $\theta$-stable.
We may assume therefore that $G$ is almost simple.

Any regular nilpotent element is semiregular.
In fact, there are no non-regular semiregular nilpotent elements except when $G$ is of type $D$ or $E$.
If $G$ is of type $D_n$, then by Lemma \ref{typed} above there exists a non-trivial involutive automorphism $\sigma:G\longrightarrow G$ such that $\lambda(k^\times)\subset G^\sigma$ and $e\in{\mathfrak g}^\sigma$.
Since $\theta=\Ad\lambda(t_0)$, $G^\sigma$ is also $\theta$-stable.
The group $G^\sigma$ is semisimple, of type $B_{n-1}$.
By Lemma \ref{sccover} we can replace $G$ by the universal covering of $G^\sigma$.
(In fact this is unnecessary, as our argument below doesn't require the assumption of simply-connectedness.)
Hence it will suffice to prove the lemma in the case where $e$ is semiregular and $G$ is not of type $D$.
For type $E$ the semiregular orbits are as follows: $E_6(reg),E_6(a_1)$; $E_7(reg),E_7(a_1),E_7(a_2)$; $E_8(reg),E_8(a_1),E_8(a_2)$ (\cite{som,premnil,mcninchsom}).

For each $\alpha\in\Phi$ denote by $U_\alpha$ be the unique connected, unipotent $T$-stable subgroup of $G$ satisfying $\Lie(U_\alpha)={\mathfrak g}_\alpha$.
Let $\epsilon_\alpha:k\longrightarrow U_\alpha$, $\alpha\in\Phi$ be isomorphisms such that $t\epsilon_\alpha(y)t^{-1}=\epsilon_\alpha(\alpha(t)y)$ for all $t\in T$, $y\in k$, and such that $n_\alpha=\epsilon_\alpha(1)\epsilon_{-\alpha}(-1)\epsilon_\alpha(1)\in N_G(T)$, $n_\alpha$ represents the reflection $s_\alpha\in W$.

Note that $\theta(\epsilon_\alpha(t))=
\left\{
\begin{array}{ll}
\epsilon_\alpha(t) & e_\alpha\in{\mathfrak k}, \\
\epsilon_\alpha(-t) & e_\alpha\in{\mathfrak p}.
\end{array}
\right.$

Let $w_0$ be the longest element of $W$ with respect to the Coxeter basis $s_{\alpha}$, $\alpha\in\Delta$.
Let $\hat\alpha$ be the longest root in $\Phi^+$ and let $\Phi_0$ be the set of roots in $\Phi$ which are orthogonal to $\hat\alpha$.
Then $\Phi_0$ is a root subsystem of $\Phi$ with basis $\Delta_0=\{\alpha\in\Delta\,|\,\alpha\bot\hat\alpha\}$.
Moreover $w_0=s_{\hat\alpha}w_0(\Phi_0)$, where $w_0(\Phi_0)$ is the longest element of $W(\Phi_0)$ with respect to the Coxeter basis $\{ s_{\alpha}:\alpha\in\Delta_0\}$.
Inductive application of this statement gives a description of $w_0$ as a product of orthogonal reflections $s_\alpha$ with $\alpha\in\Phi$.

We can now prove the lemma by means of the following observation.
Suppose $\beta_1,\beta_2,\ldots,\beta_t$ are orthogonal roots with $e_{\beta_i}\in{\mathfrak p}$ for all $i$.
Let $$g=\epsilon_{-\beta_1}(1/2)\epsilon_{-\beta_2}(1/2)\ldots\epsilon_{-\beta_t}(1/2)\epsilon_{\beta_1}(-1)\epsilon_{\beta_2}(-1)\ldots\epsilon_{\beta_t}(-1).$$
Then $g^{-1}\theta(g)=\prod_{i=1}^t \epsilon_{\beta_i}(1)\epsilon_{-\beta_i}(-1)\epsilon_{\beta_i}(1)=\prod_{i=1}^t n_i$, where $n_i = n_{\beta_i}$ for each $i$.
Moreover $\theta=\Int t_0$ and $t_0\in T$, hence the induced action of $\theta$ on $W$ is trivial.
To show that $\lambda$ is conjugate to a $\theta$-split torus, therefore, it will suffice to show that there is an element $w\in W$ which is conjugate to a product $s_{\beta_1}s_{\beta_2}\ldots s_{\beta_t}$, where the $\beta_i$ are orthogonal, $e_{\beta_i}\in{\mathfrak p}$, and such that $w\cdot\lambda=-\lambda$.
Recall that $e$ is regular unless $G$ is of type $E$.

{\it Type $A_n$.}
In this case $w_0$ is conjugate to
$\left\{
\begin{array}{ll}
s_{\alpha_1}s_{\alpha_3}\ldots s_{\alpha_n} & \mbox{if $n$ is odd,} \\
s_{\alpha_1}s_{\alpha_3}\ldots s_{\alpha_{n-1}} & \mbox{if $n$ is even.}
\end{array}
\right.$

But $\langle\lambda,\alpha_i\rangle =2$, hence $e_{\alpha_i}\in{\mathfrak p}$ for all $i$.
This proves the lemma in this case.

{\it Type $B_n$.}
Let $\beta_i=
\left\{
\begin{array}{ll}
\alpha_i+2\alpha_{i+1}+2\alpha_{i+2}+\ldots+2\alpha_{n} & \mbox{if $i$ is odd, $1\leq i\leq n$,} \\
\alpha_{i-1} & \mbox{if $i$ is even, $2\leq i\leq n$.}
\end{array}
\right.$

Then the $\beta_i$ are orthogonal, $e_{\beta_i}\in{\mathfrak p}$ for each $i$ and $w_0=s_{\beta_1}s_{\beta_2}\ldots s_{\beta_n}$.

{\it Type $C_n$.}
Let $\beta_i=2\alpha_i+2\alpha_{i+1}+\ldots+2\alpha_{n-1}+\alpha_n$ for $1\leq i\leq {n-1}$ and let $\beta_n=\alpha_n$.
Then the $\beta_i$ are orthogonal, $e_{\beta_i}\in{\mathfrak p}$, and $w_0=s_{\beta_1}s_{\beta_2}\ldots s_{\beta_n}$.

{\it Type $F_4$.}
Let $\beta_1=\hat{\alpha}=2\alpha_1+3\alpha_2+4\alpha_3+2\alpha_4,\beta_2=\alpha_2+2\alpha_3+2\alpha_4,\beta_3=\alpha_2+2\alpha_3$ and $\beta_4=\alpha_2$.
Clearly $e_{\beta_i}\in{\mathfrak p}$, the $\beta_i$ are orthogonal and $w_0=s_{\beta_1}s_{\beta_2}s_{\beta_3}s_{\beta_4}$.

{\it Type $G_2$.}
Let $\beta_1=3\alpha_1+2\alpha_2$ and $\beta_2=\alpha_1$.
Then $w_0=s_{\beta_1}s_{\beta_2}$ is the required expression for $w_0$.

{\it Type $E_6$.}
Let $\beta_1=\hat{\alpha}=\alpha_1+2\alpha_2+2\alpha_3+3\alpha_4+2\alpha_5+\alpha_6,\beta_2=\alpha_1+\alpha_3+\alpha_4+\alpha_5+\alpha_6,\beta_3=\alpha_3+\alpha_4+\alpha_5,\beta_4=\alpha_4$.
Then $w_0=s_{\beta_1}s_{\beta_2}s_{\beta_3}s_{\beta_4}$.
If $e$ is regular, then $\langle\lambda,\alpha_i\rangle =2\;\forall i$, hence $e_{\beta_i}\in{\mathfrak p}$ for all $i$.
This proves the lemma for $E_6(reg)$.

Suppose therefore that $e$ is in the semiregular orbit $E_6(a_1)$.
Then $\langle\lambda,\alpha\rangle =2$ for $\alpha_4\neq\alpha\in\Delta$, and $\langle\lambda,\alpha_4\rangle =0$.
Thus $w_0 s_{\alpha_4}\cdot\lambda=-\lambda$.
Hence it will suffice in this case to show that $s_{\beta_1}s_{\beta_2}s_{\beta_3}$ is conjugate to some element $s_{\gamma_1}s_{\gamma_2}s_{\gamma_3}\in W$ with $e_{\gamma_1},e_{\gamma_2},e_{\gamma_3}\in{\mathfrak p}$.
Let $\alpha=\hat\alpha-\alpha_2$.

Then $\alpha\in\Phi$ and $s_\alpha(\beta_i)=
\left\{
\begin{array}{ll}
\alpha_2 & \mbox{if $i=1$,} \\
-(\alpha_2+\alpha_3+2\alpha_4+\alpha_5) & \mbox{if $i=2$,} \\
-(\alpha_1+\alpha_2+\alpha_3+2\alpha_4+\alpha_5+\alpha_6) & \mbox{if $i=3.$}
\end{array}
\right.$

Therefore $s_\alpha(w_0s_{\alpha_4})s_\alpha^{-1}$ has the required form.
This completes the $E_6$ case.

{\it Type $E_7$.}
Let $\beta_1=\hat{\alpha},\beta_2=\alpha_2+\alpha_3+2\alpha_4+2\alpha_5+2\alpha_6+\alpha_7,\beta_3=\alpha_7$, $\beta_4=\alpha_2+\alpha_3+2\alpha_4+\alpha_5,\beta_5=\alpha_2,\beta_6=\alpha_3,\beta_7=\alpha_5.$
We have $w_0=s_{\beta_1}s_{\beta_2}\ldots s_{\beta_7}$.
If $e$ is regular, then $\langle \lambda,\alpha\rangle =2\;\forall\alpha\in\Delta$.
If $e$ is of type $E_7(a_1)$, then $\langle\lambda,\alpha\rangle =2$ for $\alpha_4\neq\alpha\in\Delta$ and $\langle\lambda,\alpha_4\rangle =0$.
If $e$ is of type $E_7(a_2)$ then $\langle\lambda,\alpha\rangle =
\left\{
\begin{array}{ll}
2 & \mbox{if $\alpha\in\Delta\setminus\{\alpha_4,\alpha_6\}$,} \\
0 & \mbox{if $\alpha=\alpha_4,\alpha_6$.}
\end{array}
\right.$

In each case we can see that $e_{\beta_i}\in{\mathfrak p}$ for all $i$.
Hence by our earlier observation there exists $g$ such that $n_0=g^{-1}\theta(g)\in N_G(T)$ and $n_0 T=w_0$.

{\it Type $E_8$.}
For regular $e$ we have $\langle\lambda,\alpha\rangle =2\;\forall\,\alpha\in\Delta$, for subregular $e$ (type $E_8(a_1)$) $\langle\lambda,\alpha\rangle =2$ for all $\alpha_4\neq\alpha\in\Delta$, and $\langle\lambda,\alpha_4\rangle =0$, while for the final case $E_8(a_2)$, we have
$$\langle\lambda,\alpha\rangle =
\left\{
\begin{array}{ll}
2 & \mbox{if $\alpha\in\Delta\setminus\{\alpha_4,\alpha_6\}$,} \\
0 & \mbox{if $\alpha=\alpha_4,\alpha_6$.}
\end{array}
\right.$$

Let $\hat\alpha$ be the longest element of $\Phi^+$ and let $\Phi_0$ be the subsystem of all roots orthogonal to $\hat\alpha$.
Then $\Phi_0$ is a subsystem of $\Phi$ isomorphic to $E_7$, and $\{\alpha_1,\alpha_2,\ldots,\alpha_7\}$ is a basis for $\Phi_0$.
Identify $\Phi_0$ with $E_7$ and let $\beta_1,\beta_2,\ldots,\beta_7$ be the orthogonal roots given for the $E_7$ case above.
Then $w_0=s_{\hat\alpha}s_{\beta_1}s_{\beta_2}\ldots s_{\beta_7}$.
Moreover, it is easy to see that $e_{\hat\alpha},e_{\beta_1},e_{\beta_2},\ldots,e_{\beta_7}\in{\mathfrak p}$.
Hence there exists $g\in G$ such that $g^{-1}\theta(g)\in N_G(T)$ represents $w_0$.
This completes the proof.
%\qed
\end{proof}

Let $A$ be a maximal $\theta$-split torus of $G$.
The roots $\Phi_A=\Phi(G,A)$ form a non-reduced root system \cite[4.7]{rich2}.
Let $\Pi$ be a basis for $\Phi_A$.
We can now use Lemma \ref{splitconj} to give a criterion for $e\in{\cal N}$ to be regular.

\begin{lemma}\label{regconj}
There exists a cocharacter $\omega:k^\times\longrightarrow A\cap{G^{(1)}}$ such that $\langle\omega,\alpha\rangle =2\;\forall\alpha\in\Pi$.
Let $e\in{\cal N}$ and let $\lambda:k^\times\longrightarrow K$ be associated to $e$.
Then $e$ is regular if and only if $\lambda$ is $G$-conjugate to $\omega$.
Hence the set ${\cal N}_{reg}$ of regular nilpotent elements is contained in a single $G$-orbit.
\end{lemma}

\begin{proof}
By Lemma \ref{splitconj}, $\lambda$ is $G$-conjugate to a $\theta$-split cocharacter $\mu$.
But any two maximal $\theta$-split tori are conjugate by an element of $K$, hence we may assume that $\mu(k^\times)\subset A$.
Moreover, we may assume after conjugating further by an element of $N_K(A)$, if necessary, that $\langle\mu,\alpha\rangle\geq 0$ for all $\alpha\in\Pi$.

It follows from the properties of associated cocharacters (see for example \cite[Thm. 2.3(iv)]{premnil}) that $\dim{\mathfrak z}_{\mathfrak g}(e)=\dim{\mathfrak g}(0;\lambda)+\dim{\mathfrak g}(1;\lambda)=\dim{\mathfrak g}(0;\mu)+\dim{\mathfrak g}(1;\mu)$.
But $\mu(k^\times)\subset A$, hence $\dim{\mathfrak g}(0;\mu)\geq \dim{\mathfrak z}_{\mathfrak g}({\mathfrak a})$.
Thus by Lemma \ref{regs}, $e$ is regular if and only if $\mu$ is regular in $A$ and all weights of $\Ad\mu$ on ${\mathfrak g}$ are even.
Let $S$ be a maximal torus of $G$ containing $A$.
By Lemma \ref{basis} there exists a basis $\Delta_S$ for $S$ such that every element of $\Pi$ can be written in the form $\beta|_A$ for some $\beta\in\Delta_S$.
Hence by properties of weighted Dynkin diagrams, $\langle\mu,\alpha\rangle\in\{0,1,2\}$ for each $\alpha\in\Pi$.
It follows that $e$ is regular if and only if $\langle\mu,\alpha\rangle=2$ for all $\alpha\in\Pi$.
But there exists some regular nilpotent element; hence $\omega$ exists.
%\qed
\end{proof}

\begin{rk}\label{eiseven}
Let $S$ be a maximal torus of $G$ containing $A$ and let $\Delta_S$ be a basis for $\Phi_S=\Phi(G,S)$, such that $\{\alpha|_A\, :\,\alpha\in\Delta_S,\alpha|_A\neq 1\}$ is a basis for $\Phi_A$.
Let $I=\{\alpha\in\Delta_S\, :\,\alpha|_A=1\}$.
Then $\omega$ satisfies $\langle\alpha,\omega\rangle=
\left\{
\begin{array}{ll}
0 & \mbox{if $\alpha\in I$,} \\
2 & \mbox{if $\alpha\in\Delta_S\setminus I$.}
\end{array}
\right.$
\end{rk}

\begin{corollary}
Let $e$ be a regular nilpotent element of ${\mathfrak p}$.
Then $e$ is even.
\end{corollary}

\begin{proof}
Let $\lambda$ be an associated cocharacter for $e$.
Then $\lambda$ is conjugate to $\omega$.
But now by the remark above $\omega$ is even.
%\qed
\end{proof}

Fix a cocharacter $\omega$ as in Lemma \ref{regconj} and denote by $Y_\omega$ the open $Z_G(\omega)$-orbit in ${\mathfrak g}(2;\omega)$.

\begin{lemma}\label{ainZ}
Let $E\in Y_\omega$.
Suppose $a\in A$ and $a\cdot E=E$.
Then $a\in Z(G)$.
\end{lemma}

\begin{proof}
Since $Z_G(\omega)\cdot E=Y_\omega$ and $Z_G(\omega)=Z_G(A)$, it follows that $a\cdot E'=E'$ for all $E'\in Y_\omega$.
Therefore $a\cdot E'=E'$ for all $E'\in{\mathfrak g}(2;\omega)$, which implies that $\alpha(a)=1\;\forall \alpha\in\Pi$.
It follows that $a\in Z(G)$.
%\qed
\end{proof}

\begin{lemma}\label{omega}
Let $e\in{\cal N}$ be regular and let $\lambda:k^\times\longrightarrow K$ be associated to $e$.
Let $g\in G$ be such that $g\cdot e\in{\mathfrak p}$ and $(g\cdot \lambda)(k^\times)\subset K$.
Then $g\in K^*$.
In particular $C=Z_G(\lambda)\cap Z_G(E)\subseteq K^*$.
\end{lemma}

\begin{proof}
Let $g$ be such that $g\cdot e\in{\mathfrak p}$ and $(g\cdot\lambda)(k^\times)\subset K$, and let $x=g^{-1}\theta(g)$.
Assume first of all that $x$ is semisimple.
By \cite[6.3]{rich2} there exists a maximal $\theta$-split torus of $G$ containing $x$.
Hence, after conjugating $e,\lambda$, and $g$ by a suitable element of $K$, we may assume that $x\in A$.
Let $H=Z_G(x)^\circ$ and let ${\mathfrak h}=\Lie(H)$.
We claim that $\lambda$ is an associated cocharacter for $e$ in $H$.
Let $d=\min_{y\in{\mathfrak h}\cap{\mathfrak p}}\dim Z_H(y)$; since $Z_G(A)\subset H$, $d=\dim Z_G(A)$.
Thus $Z_G(e)^\circ\subset H$.
In particular, $C^\circ\subset H$.
Recall that $C^\circ$ is a ($\theta$-stable) reductive subgroup of $G$.
Hence we can choose a $\theta$-stable maximal torus $S$ of $C^\circ$ (\cite[7.5]{steinberg}).
Let $L=Z_G(S)$, a $\theta$-stable Levi subgroup of $G$.
By \cite[Prop. 2.5]{premnil}, $e$ is distinguished in ${\mathfrak l}=\Lie(L)$ and $\lambda(k^\times)\subset L^{(1)}$.
Clearly $x\in L$.
Hence $e$ is distinguished in $Z_L(x)^\circ=Z_H(S)=Z_G(x,S)^\circ$.
Let $T$ be a maximal torus of $Z_H(S)$ containing $\lambda(k^\times)$.
Then $T=(T\cap L^{(1)})\cdot Z(L)^\circ = (T\cap Z_H(S)^{(1)})\cdot Z(L)^\circ$.
Therefore $\lambda(k^\times)\subset Z_H(S)^{(1)}$, that is, $\lambda$ is an associated cocharacter for $e$ in $H$.

Since $A\subset H$, we can consider $\Phi(H,A)$ as a subset of $\Phi_A$.
Let $\Phi(H,A)^+=\Phi(H,A)\cap \Phi_A^+$ and let $\Pi_H$ be the corresponding basis for $\Phi(H,A)$.
By Lemma \ref{regconj} there exists $\omega_H:k^\times\rightarrow A\cap H^{(1)}$ such that $\langle\alpha,\omega\rangle=2$ for all $\alpha\in\Pi_H$, and $h\in H$ such that $h\cdot\lambda=\omega_H$.
But $\lambda$ is $G$-conjugate to $\omega$: hence, since $\omega$ and $\omega_H$ are in the same Weyl chamber in $Y(A)$, we must have $\omega=\omega_H$.
Thus $h\cdot\lambda=\omega$ and $E=h\cdot e\in Y_\omega$.
Moreover, $x\cdot E=E$.
Now by Lemma \ref{ainZ}, $x\in Z(G)$.

Suppose therefore that $x$ is not semisimple.
Let $x=su$ be the Jordan-Chevalley decomposition of $x$.
Since $x\in C$, $s,u\in C$ also.
By \cite[III.3.15]{sands}, all unipotent elements of $Z_G(e)$ are in $Z_G(e)^\circ$.
Hence by \cite[Pf. of Thm. 2.3, p.347]{premnil}, $u\in C^\circ$.
But now $\theta$ acts non-trivially on the derived subgroup of (the reductive group) $C^\circ$, hence there exists a non-central $\theta$-split torus in $C^\circ$ (\cite[\S 1]{vust}).
This contradicts the assumption that $e$ is regular, by the above.
%\qed
\end{proof}

Thus we have our desired reward.

\begin{theorem}\label{regorbit}
The set ${\cal N}_{reg}$ of regular nilpotent elements of ${\mathfrak p}$ is a single $K^*$-orbit.
Hence $K^*$ permutes the irreducible components of ${\cal N}$ transitively and ${\cal N}$ is the closure of the regular nilpotent $K^*$-orbit.
\end{theorem}

\begin{proof}
Let $e\in{\cal N}_{reg}$ and let $\lambda:k^\times\longrightarrow K$ be an associated cocharacter for $e$.
By Cor. \ref{regconj}, ${\cal N}_{reg}=G\cdot e\cap{\mathfrak p}$.
Suppose $g\in G$ and $e'=g\cdot e\in{\mathfrak p}$.
By Lemma \ref{assoc} there exists an associated cocharacter $\mu:k^\times\longrightarrow K$ for $e'$.
Moreover $\mu$ is $Z_G(e')^\circ$-conjugate to $g\cdot\lambda$.
Hence there exists $h\in G$ such that $h\cdot e=e'=g\cdot e$ and $h\cdot\lambda=\mu$.
But now by Lemma \ref{omega}, $h\in K^*$.

We have proved that any element of ${\cal N}_{reg}$ is $K^*$-conjugate to $e$.
The regular elements are dense in each irreducible component by Thm. \ref{nil1}.
But therefore $\overline{{\cal N}_{reg}}={\cal N}$.
This completes the proof.
%\qed
\end{proof}

Thm. \ref{regorbit} generalises \cite[Thm. 6]{kostrall} to good positive characteristic.
In \cite{sek}, Sekiguchi determined (for $k={\mathbb C}$) the
involutions for which the set of nilpotent elements is non-irreducible.
The proof comes down to checking which elements of the group $F=\{
a\in A\,|\,a^2\in Z(G)\}$ stabilize a particular irreducible component
of ${\cal N}$.
The calculations in the classical case were omitted.
Fortunately, our analysis of associated cocharacters, together with
the classification of involutions (\cite{springer}), considerably simplify the task of generalizing
Sekiguchi's results.
We begin with the following:

\begin{theorem}\label{gthetaorbs}
Let $e,\lambda,C$ be as above.
Let $Z=Z(G)$, $P=\{ g^{-1}\theta(g)\,|\, g\in G\}$, $\tau:G\longrightarrow P$, $g\mapsto g^{-1}\theta(g)$ and denote by $\Gamma$ the set of $G^\theta$-orbits in ${\cal N}_{reg}$.

(a) The map from $K^*$ to $\Gamma$ given by $g\mapsto gG^\theta\cdot e$ is surjective and induces a one-to-one correspondence $K^*/G^\theta C\longrightarrow\Gamma$.

(b) The morphism $\tau$ induces an isomorphism $K^*/G^\theta C\longrightarrow (Z\cap A)/{\tau(C)}$.
Since $Z\subseteq C$, there is a surjective map $(Z\cap A)/{\tau(Z)}\longrightarrow (Z\cap A)/{\tau(C)}$.

(c) The embedding $F^*\hookrightarrow K^*$ induces a surjective map $F^*/{F(Z\cap A)}\rightarrow\Gamma$.

(d) The map $F^*\rightarrow Z\cap A$, $a\mapsto a^2$ induces an isomorphism of finte groups $F^*/{F(Z\cap A)}\longrightarrow Z\cap A/{(Z\cap A)^2}$.
\end{theorem}

\begin{proof}
Since $K^*$ permutes the elements of ${\cal N}_{reg}$ transitively, the map in (a) from $K^*$ to $\Gamma$ is surjective and factors through $G^\theta C$.
Suppose $g,g'\in K^*$ and $gG^\theta\cdot e=g' G^\theta\cdot e$.
Then there exists $x\in G^\theta$ such that $g^{-1}g'\cdot e=x\cdot e$.
Moreover, since $g^{-1}g'\cdot\lambda$ is an associated cocharacter for $x\cdot e$ and $g^{-1}g'\cdot\lambda(k^\times)\subset K$, there exists $y\in Z_K(e)^\circ$ such that $yx\cdot\lambda = g^{-1}g'\cdot\lambda$ by Cor. \ref{assoc}.
Thus $g\in g'CG^\theta=g' G^\theta C$.
Hence the map $K^*/G^\theta C\rightarrow\Gamma$ is one-to-one.
This proves (a).

Since $K^*=\tau^{-1}(Z\cap A)$, the induced map $\overline{\tau}$ from $K^*$ to $Z\cap A/{\tau(C)}$ is surjective.
Suppose therefore that $g\in K^*$ and that there exists $c\in C$ such that $g^{-1}\theta(g)=c^{-1}\theta(c)$.
Then $gc^{-1}\in G^\theta$.
Hence $g\in CG^\theta=G^\theta C$.
It follows that the kernel of $\overline{\tau}$ is $G^\theta C$.

We recall by \cite[8.1]{rich2} that $K^*=F^*\cdot K$.
Hence there is a surjective map $F^*\rightarrow\Gamma$, $a\mapsto aG^\theta\cdot e$.
Moreover, since $F\subset G^\theta$ and $az\cdot e = a\cdot e$ for any $a\in F^*,z\in (Z\cap A)$, this map factors through the cosets of $F(Z\cap A)$ in $F^*$.
This proves (c).
Finally, the homomorphism $F^*\rightarrow Z\cap A$, $a\mapsto a^2$ is surjective by the definition of $F^*$ and the fact that $A$ is a torus.
Suppose $a^2 =z^2$ for some $z\in Z\cap A$.
Then $(z^{-1}a)^2$ is the identity element.
Hence $z^{-1}a\in F\Rightarrow a\in F(Z\cap A)$.
This completes the proof.
%\qed
\end{proof}

An involution is {\it split} (or of {\it maximal rank}) if the maximal $\theta$-split torus $A$ is a maximal torus of $G$, and {\it quasi-split} if $Z_G(A)$ is a maximal torus of $G$.
Recall (see Sect. \ref{sec:2.2}) that, relative to a maximal torus $S$ containing $A$, there is a basis $\Delta_S$ for $\Phi_S$, a subset $I$ of $\Delta_S$, and a graph automorphism $\psi$ of $\Phi_S$ such that $\theta^*(\beta)=-w_I(\psi(\beta))$ for any $\beta\in\Phi_S$.
With this notation, $\theta$ is quasi-split if $I=\emptyset$, and is split if in addition the action of $\psi$ is trivial.

\begin{corollary}\label{splitcmpts}
Suppose $G$ is almost simple and simply-connected.

(a) Let $\theta$ be split.
The irreducible components of ${\cal N}$ are in one-to-one correspondence with the elements of $Z/Z^2$.
Hence ${\cal N}$ has 4 components if $G$ is of type $D_{2n}$, has 2 components if $G$ is of type $A_{2n-1},B_n,C_n,D_{2n+1},E_7$, and is irreducible if $G$ is of type $A_{2n},E_6,E_8,F_4$, or $G_2$.

(b) Let $\theta$ be quasi-split.
Then the irreducible components of ${\cal N}$ are in one-to-one correspondence with the elements of $(Z\cap A)/{\tau(Z)}$.

(c) Let $\theta$ be any involutive automorphism and let $G$ be of one of the following types: $A_{2n},E_6,E_8,F_4$, or $G_2$.
Then ${\cal N}$ is irreducible.
\end{corollary}

\begin{proof}
Since $G$ is semisimple and simply-connected, the isotropy subgroup $G^\theta$ is connected by \cite[8.1]{steinberg}.
Hence the irreducible components of ${\cal N}$ are in one-to-one correspondence with the elements of $Z\cap A/\tau(C)$ by Thm. \ref{gthetaorbs}.
If $\theta$ is split or quasi-split, then a regular nilpotent element of ${\mathfrak p}$ is also a regular element of ${\mathfrak g}$, hence $C=Z(G)$.
Thus $\tau(C)=\tau(Z)$.
If $\theta$ is split, then $A$ is a maximal torus of $G$, hence $Z\subset A$.
This proves (a) and (b).
For (c), the centre $Z$ of $G$ has odd order, hence so does $Z\cap A$.
Therefore $(Z\cap A)/{(Z\cap A)^2}$ is trivial.
But now by Thm. \ref{gthetaorbs}(d), ${\cal N}$ is irreducible.
%\qed
\end{proof}

Note that by Rk. \ref{nonsc}, the description of the number of irreducible components of ${\cal N}$ holds without the assumption of simply-connectedness.
Using the notation $({\mathfrak g},{\mathfrak k})$, the split
involutions are as follows:

 - Type $A_n$, $(\mathfrak{sl}(n+1),\mathfrak{so}(n+1))$ (or
   $(\mathfrak{gl}(n+1),\mathfrak{so}(n+1))$ if $p\, |\,(n+1)$),

 - Type $B_n$,
   $(\mathfrak{so}(2n+1),\mathfrak{so}(n)\oplus\mathfrak{so}(n+1))$,

 - Type $C_n$, $(\mathfrak{sp}(2n),\mathfrak{gl}(n))$,

 - Type $D_n$,
   $(\mathfrak{so}(2n),\mathfrak{so}(n)\oplus\mathfrak{so}(n))$,

 - Type $E_6$, $({\mathfrak e}_6,\mathfrak{sp}(8))$,

 - Type $E_7$, $({\mathfrak e}_7,\mathfrak{sl}(8))$,

 - Type $E_8$, $({\mathfrak e}_8,\mathfrak{so}(16))$,

 - Type $F_4$, $({\mathfrak f}_4,\mathfrak{sp}(6)\oplus\mathfrak{sl}(2))$,

 - Type $G_2$, $({\mathfrak g}_2,\mathfrak{sl}(2)\oplus\mathfrak{sl}(2))$.

Hence Cor. \ref{splitcmpts} confirms no. 2 of Table 1, and no.s 1,2,3,4,6 of Table 2, listed
in \cite[p. 161]{sek}.
In Sect. 6.3 we deal with the remaining cases.

\subsection{A $\theta$-equivariant Springer isomorphism}
\label{sec:5.5}

Assume once more that $G$ satisfies the conditions (A)-(C) of \S 3.
Let ${\cal U}(G)$ be the closed set of unipotent elements in $G$ and let ${\cal N}({\mathfrak g})$ be the nilpotent cone in ${\mathfrak g}$.
We let ${\cal U}=\{ u\in{\cal U}(G)\,|\,\theta(u)=u^{-1}\}$.
By \cite[6.1]{rich2}, ${\cal U}\subset P$, where $P=\{ g^{-1}\theta(g)\,|\,g\in G\}$.
It is well-known (see for example \cite{sands}) that if the characteristic of $k$ is good for $G$, then there exists a $G$-equivariant isomorphism of affine varieties $\psi:{\cal U}(G)\longrightarrow{\cal N}({\mathfrak g})$, sometimes known as the Springer map.
It was also stated without proof in \cite[\S 10]{barrich} that there is a $K$-equivariant isomorphism from ${\cal U}$ to ${\cal N}$.
We get the desired result in our case with the following proposition.
Part (c) is due to McNinch (\cite[Thm. 35]{mcninch}).

\begin{proposition}\label{iso}
There is a $G$-equivariant isomorphism of affine varieties $\Psi:{\cal U}(G)\longrightarrow{\cal N}({\mathfrak g})$ such that:

(a) $\Psi(u^{-1})=-\Psi(u)$ ($u\in {\cal U}(G)$),

(b) $\Psi(\theta(u))=d\theta(\Psi(u))$ ($u\in{\cal U}(G)$),

(c) $\Psi(u^p)=\Psi(u)^{[p]}$ ($u\in{\cal U}(G)$).

Moreover, if (i) $p>3$ or (ii) $G$ has no component of type $D_4$, then we may assume that (b) holds for all automorphisms of $G$.
\end{proposition}

\begin{proof}
As ${\cal U}(G)\subseteq G^{(1)}$ and ${\cal N}({\mathfrak g})\subseteq\Lie(G^{(1)})$ we may assume that $G$ is semisimple.
Let $G_1,G_2,\ldots,G_l$ be the minimal normal subgroups of $G$ and let ${\mathfrak g}_i=\Lie(G_i)$ for $1\leq i\leq l$.
Then $G=G_1\times G_2\times\ldots \times G_l$ and ${\mathfrak g}={\mathfrak g}_1\oplus{\mathfrak g}_2\oplus\ldots\oplus{\mathfrak g}_l$.
Let $H$ (resp. $L$) be the subgroup of $G$ generated by all $G_i$ isomorphic to $G_1$ (resp. all $G_i$ not isomorphic to $G_1$) and let ${\mathfrak h}=\Lie(H),{\mathfrak l}=\Lie(L)$.
Then $G=H\times L$ and ${\mathfrak g}={\mathfrak h}\oplus{\mathfrak l}$.
Moreover ${\cal U}(G)={\cal U}(H)\times{\cal U}(L),{\cal N}({\mathfrak g})={\cal N}({\mathfrak h})\oplus{\cal N}({\mathfrak l})$.
Any automorphism of $G$ stabilizes $H$ and $L$.
Hence we may assume that all minimal normal subgroups of $G$ are isomorphic to $G_1$.
Identify $G$ with the product $G_1\times G_1\times\ldots \times G_1$ ($l$ times).
Thus we write an element of $G$ as $(g_1,g_2,\ldots ,g_l)$, $g_i\in G_1$.
The symmetric group $S_l$ acts on $G$: $\tau(g_1,g_2,\ldots,g_l)=(g_{\tau(1)},g_{\tau(2)},\ldots ,g_{\tau(l)})$.
Furthermore, any automorphism of $G$ can be written in the form $\tau\circ(\theta_1,\theta_2,\ldots ,\theta_l)$, where $\theta_i\in\Aut(G_1)$, $(\theta_1,\theta_2,\ldots ,\theta_l)(g_1,g_2,\ldots,g_l)=(\theta_1(g_1),\theta_2(g_2),\ldots,\theta_l(g_l))$ and $\tau\in S_l$.
Thus it will suffice to prove the proposition in the case where $G$ is almost simple.
There are three cases: (i) $G$ is not of type $A_n$, (ii) $G=\SL(n,k)$ with $p\nmid n$, and (iii) $G=\SL(n,k)$ with $p\, |\, n$.
In case (iii) replace $G$ by $\GL(n,k)$.

In all three cases, it is well-known (see for example \cite[I.5]{sands}) that there exists a representation $\rho:G\longrightarrow\GL(V)$ such that:

(i) $d\rho:{\mathfrak g}\longrightarrow\mathfrak{gl}(V)$ is injective,

(ii) The associated trace form $\kappa_\rho:{\mathfrak g}\times{\mathfrak g}\longrightarrow k$, $(x,y)\mapsto\tr(d\rho(x),d\rho(y))$ is non-degenerate.

We construct a new representation $\sigma:G\longrightarrow\GL(V\oplus V)$ defined by $g\mapsto
\left(
\begin{array}{ll}
\rho(g) & 0 \\
0 & {^t}\rho(g)^{-1}
\end{array}
\right).
$

The associated trace form $\kappa_\sigma=2\kappa_\rho$.
Replacing $(\rho,V)$ by $(\sigma,V\oplus V)$, we may assume that $(\rho,V)$ satisfies the further properties:

(iii) $d\rho({\mathfrak g})\subseteq\mathfrak{sl}(V)$,

(iv) $\tr(\rho(g)d\rho(x))=-\tr(\rho(g^{-1})d\rho(x))$ for all $g\in G,x\in{\mathfrak g}$.

Finally, construct another representation $\sigma:{\mathfrak g}\longrightarrow\mathfrak{gl}(V\oplus V)$ defined by $g\mapsto
\left(
\begin{array}{ll}
\rho(g) & 0 \\
0 & \rho(\theta(g))
\end{array}
\right) \in\GL(V\oplus V)$.

By the $\theta$-invariance of the trace (see the proof of Thm. \ref{redthm}) $\kappa_\sigma=2\kappa_\rho$.
Moreover, it is easy to see that $\sigma$ satisfies (i)-(iv) and that:

(v) $\tr(\sigma(\theta(g))d\sigma(x))=\tr(\sigma(g)d\sigma(d\theta(x)))$ for all $g\in G,x\in{\mathfrak g}$.

Identify ${\mathfrak g}$ with its image $d\sigma({\mathfrak g})$ and let ${\mathfrak g}^\bot=\{ x\in\mathfrak{gl}(V)| \tr(xy)=0\,\forall y\in{\mathfrak g}\}$.
It follows from (ii) and (iii) that $\mathfrak{gl}(V)={\mathfrak g}\oplus{\mathfrak g}^\bot$ and that $I_V\in{\mathfrak g}^\bot$.
Let $\iota:\GL(V)\hookrightarrow\mathfrak{gl}(V)$ be the map embedding $\GL(V)$ as a Zariski open subset of $\mathfrak{gl}(V)$ and let $\pr_{\mathfrak g}:\mathfrak{gl}(V)\twoheadrightarrow{\mathfrak g}$ be the projection onto ${\mathfrak g}$ induced by the direct sum decomposition $\mathfrak{gl}(V)={\mathfrak g}\oplus{\mathfrak g}^\bot$.
Introduce the map $\eta=\pr_{\mathfrak g}\circ\iota\circ\sigma:G\longrightarrow{\mathfrak g}$.
It follows from \cite[Cor. 6.3]{barrich} that $\eta$ restricts to an isomorphism $\Psi:{\cal U}(G)\longrightarrow{\cal N}({\mathfrak g})$.

We claim that (iv) and (v) imply, respectively, (a) and (b) of the proposition.
Identify $\GL(V)$ with its image $\iota(\GL(V))$.
By (iv) we have $\kappa_\sigma(\eta(g),x)=-\kappa_\sigma(\eta(g^{-1}),x)$ for all $x\in{\mathfrak g}$.
It follows that $\eta(g^{-1})=-\eta(g)$.
This proves (a).
By (v), $\kappa_\sigma(\eta(\theta(g)),x)=\kappa_\sigma(\eta(g),d\theta(x))$ for all $x\in{\mathfrak g}$.
But $\kappa_\sigma(d\theta(\eta(g)),x)=\kappa_\sigma(\eta(g),d\theta(x))$ for all $x\in{\mathfrak g}$, hence $d\theta(\eta(g))=\eta(\theta(g))$ for any $g\in G$.
This proves (b).

The proof that $\eta(g^p)=\eta(g)^{[p]}$ is in \cite[Thm. 35]{mcninch}.
It can be applied perfectly well here without affecting the rest of the proof.

We have constructed the isomorphism $\Psi$ invariant with respect to a given involution $\theta$.
But $\Aut G$ is generated over $\Int G$ by the group $\Gamma$ of graph automorphisms (for $G=\GL(n,k)$ with $p\, |\, n$ and $n\neq 2$ this follows from Lemma \ref{GLautos}).
Moreover the group of graph automorphisms is either trivial, or cyclic of order 2 (for types $A_n\;(n\geq 2),D_n\;(n\geq 5)$, and $E_6$), or isomorphic to the symmetric group $S_3$ (for type $D_4$).

Choose a set of coset representatives $C$ for $\Gamma$.
If $p>3$ then we can easily adapt the proof above to make $\eta$ invariant with respect to every element of $C$.
If there is a component of type $D_4$, then we need the assumption $p>3$ for the trace form $\kappa_\sigma$ to be non-zero.
Hence it is straightforward with these restrictions to construct an isomorphism $\Psi$ satisfying (b) for every element of $C$.
But then $\Psi$ satisfies (b) for every element of $\Aut G$.
%\qed
\end{proof}

\begin{corollary}\label{isocor}
There is a $K^*$-equivariant isomorphism of affine varieties $\Psi:{\cal U}\longrightarrow{\cal N}$.
\end{corollary}

\section{A reductive subalgebra}
\label{sec:6}

\subsection{Preparation}
\label{sec:6.1}

Fix a cocharacter $\omega:k^\times\longrightarrow A$ as in Lemma \ref{regconj}, and let $Y_\omega=\{ x\in{\mathfrak g}(2;\omega)\,|\,\overline{Z_G(\omega)\cdot x}={\mathfrak g}(2;\omega)\}$, $Y_{-\omega}=\{ x\in{\mathfrak g}(-2;\omega)\,|$ $\overline{Z_G(\omega)\cdot x}={\mathfrak g}(-2;\omega)\}$.
Then $\omega$ (resp. $-\omega$) is an associated cocharacter for any $x\in Y_\omega$ (resp. $x\in Y_{-\omega}$).
Let $S$ be a maximal torus of $G$ containing $A$.
Recall (\cite{springer2} and \cite[1.3-4]{springer} - see also Sect. \ref{sec:2.2}) that there exists a basis $\Delta_S$ for $\Phi_S$, a subset $I$ of $\Delta_S$, and a graph automorphism $\psi:\Phi_S\rightarrow\Phi_S$ (stabilizing $\Delta_S$ and $I$) such that:

 - $\theta^*(\alpha)=-w_I(\psi(\alpha))$, $\alpha\in\Phi_S$,

 - $\theta^*(\alpha)=\alpha$, $\alpha\in I$,

 - $\alpha|_A = 1$ if $\alpha\in I$, and for $\alpha,\beta\in\Delta_S\setminus I$, $\alpha|_A=\beta|_A\;\Leftrightarrow\;beta\in\{\alpha,\psi(\alpha)\}$.

 - The set $\Pi=\{\alpha|_A\,:\,\alpha\in\Delta_S\setminus I\}$ is a basis for $\Phi_A$.

Fix $S,\Delta_S,I,\psi,\Pi$ as above.
Let $\Phi_A^*$ be the set of $\alpha\in\Phi_A$ such that $\alpha /2\notin\Phi_A$.
For $\alpha\in\Phi_A$, denote by $\Psi_\alpha$ the set of all $\beta\in\Phi_S$ such that $\beta|_A$ is an integer multiple of $\alpha$: $\Psi_\alpha$ is a closed symmetric subset of $\Phi_S$.
For $\beta\in\Phi_S$ let $U_\beta$ be the unique closed connected $S$-stable subgroup of $G$ such that $\Lie(U_\beta)={\mathfrak g}_\beta$.
Let $L_\alpha$ be the subgroup of $G$ generated by $S$ together with all subgroups $U_\beta,\beta\in\Psi_\alpha$.
Then $L_\alpha$ is a $\theta$-stable connected reductive subgroup of $G$ and $U_\beta\subset L_\alpha$ if and only if $\beta\in\Psi_\alpha$ (\cite[Pf. of 4.6]{rich2}).
In fact, we are only concerned here with the following case:

\begin{lemma}
Let $\alpha\in\Pi$.
Then $L_\alpha$ is a standard Levi subgroup of $G$ relative to $(S,\Delta_S)$.
\end{lemma}

\begin{proof}
Let $\beta\in\Delta_S$ be such that $\beta|_A=\alpha$.
Then $\theta^*(\beta)=-w_I(\psi(\beta))\in -(\psi(\beta)+{\mathbb Z}I)$.
Hence $\Psi_\alpha=\Phi_J$, where $J=I\cup\{\beta,\psi(\beta)\}$.
%\qed
\end{proof}

Recall (\cite[\S 1]{vust}) that $Z_G(\omega)=Z_G(A)=M\cdot A$ (almost direct product), where $M=Z_K(A)^\circ$.
It is clear from the definition that $Z_{L_\alpha}(\omega_\alpha)=Z_G(A)$.
Once more we denote by $\langle .\, ,.\rangle:X(A)\times Y(A)\longrightarrow {\mathbb Z}$ the natural pairing of abelian groups.

\begin{corollary}\label{omegaalphacor}
There exists a cocharacter $\omega_\alpha:k^\times\longrightarrow A\cap L_\alpha^{(1)}$ such that $\langle\alpha,\omega_\alpha\rangle=2$.
We have $\omega_\alpha=\omega+\mu_\alpha$ for some $\mu_\alpha\in Y((Z(L_\alpha)\cap A)^\circ)$.
\end{corollary}

\begin{proof}
All of our earlier results apply to the $\theta$-stable Levi subgroup $L_\alpha$ of $G$.
In particular, there exists a cocharacter $\omega_\alpha:k^\times\longrightarrow A\cap L_\alpha^{(1)}$ such that $\langle\alpha,\omega_\alpha\rangle=2$ by Lemma \ref{regconj}.
Now, clearly $(\omega_\alpha-\omega)\in Y(A)$.
But $\langle\alpha,\omega_\alpha-\omega\rangle=0$, hence $\omega_\alpha-\omega\in Y(Z(L_\alpha)^\circ)$.
%\qed
\end{proof}

Let $E=X(A)\otimes_{\mathbb Z}{\mathbb R}$ and let $(.\, ,.):E\times E\rightarrow {\mathbb R}$ be a $W_A$-equivariant inner product.
The set $\Phi_A^*$ is a root system in $E$ with Cartan integers $\langle\alpha,\beta\rangle=2(\alpha,\beta)/{(\beta,\beta)}$, $\alpha,\beta\in\Pi$ (\cite[\S 4]{rich2}).

\begin{lemma}\label{omegaalphacartan}
We have $\langle\beta,\omega_\alpha\rangle=\langle\beta,\alpha\rangle$ for all $\alpha,\beta\in\Pi$.
\end{lemma}

\begin{proof}
Let $E^*$ be the dual space to $E$, naturally identified with $Y(A)\otimes_{\mathbb Z}{\mathbb R}$.
The inner product $(.\, ,.)$ induces a $W_A$-equivariant isomorphism $E\rightarrow E^*$.
Note that for $x\in E$, $s_\alpha(x)=-x\;\Leftrightarrow\; x\in{\mathbb R}\alpha$.
Moreover, $E^*={\mathbb R}\omega_\alpha\oplus (Y((Z(L_\alpha)\cap A)^\circ)\otimes_{\mathbb Z}{\mathbb R})$.
Hence for $y\in E^*$, $s_\alpha(y)=-y\;\Leftrightarrow\; y\in{\mathbb R}\omega_\alpha$.
It follows that the isomorphism $E\rightarrow E^*$ sends $\alpha$ to $c\omega_\alpha$ for some $c\in{\mathbb R}^\times$.
Thus $(\beta,\alpha)=c\langle\beta,\omega_\alpha\rangle$ for all $\beta\in\Phi_A$.
But $\langle\alpha,\omega_\alpha\rangle=2$, hence $c=(\alpha,\alpha)/2$.
Therefore $\langle\beta,\omega_\alpha\rangle=2(\beta,\alpha)/(\alpha,\alpha)= \langle\beta,\alpha\rangle$ for all $\alpha,\beta\in\Pi$.
%\qed
\end{proof}

It follows from the construction of $\omega_\alpha$ that there is an open $Z_G(\omega)$-orbit on ${\mathfrak g}(\alpha;A)$, which we denote $Y_\alpha$.
Since $L_\alpha$ is a Levi subgroup of $G$, $\omega_\alpha$ is an associated cocharacter (in $G$) for any $x_\alpha\in Y_\alpha$.

\begin{lemma}\label{sl2s}
Let $E_\alpha\in Y_\alpha$.
Then $d\omega_\alpha(1)=\xi_\alpha [E_\alpha,d\theta(E_\alpha)]$ for some $\xi_\alpha\in k^\times$.
\end{lemma}

\begin{proof}
By properties of associated cocharacters, ${\mathfrak z}_{\mathfrak g}(E_\alpha)\cap {\mathfrak g}(-\alpha;A)=0$.
Hence $[E_\alpha,d\theta(E_\alpha)]\neq 0$.
But $\dim A\cap L_\alpha^{(1)}=1$, hence $\dim{\mathfrak a}\cap\Lie(L_\alpha^{(1)})=1$.
It follows that there exists $\xi_\alpha\in k^\times$ such that $d\omega_\alpha(1)=\xi_\alpha [E_\alpha,d\theta(E_\alpha)]$.
%\qed
\end{proof}

\begin{lemma}\label{diffindpt}
The differentials $d\alpha:{\mathfrak a}\longrightarrow{\mathfrak a}$, $\alpha\in\Pi$, are linearly independent.
\end{lemma}

\begin{proof}
It follows at once from the definitions that $\cap_{\alpha\in\Pi}\ker d\alpha={\mathfrak z}({\mathfrak g})\cap{\mathfrak a}$.
Moreover, ${\mathfrak z}({\mathfrak g})=\Lie(Z(G)^\circ)$ by \cite[2.3]{me}.
But $Z(G)^\circ$ is a $\theta$-stable torus, hence $Z(G)^\circ=(Z\cap K)^\circ\cdot (Z\cap A)^\circ$ by Lemma \ref{stabletori}.
Therefore ${\mathfrak z}({\mathfrak g})\cap{\mathfrak a}=\Lie((Z\cap A)^\circ)$.
But $\dim A-\dim (Z\cap A)^\circ=\rank\Phi_A$ (see for example \cite[Rk. 4.8]{rich2}).
This completes the proof.
%\qed
\end{proof}

\begin{corollary}\label{omegasindpt}
The toral elements $d\omega_\alpha(1)$ are linearly independent.
\end{corollary}

\begin{proof}
Let $E_\alpha\in Y_{\alpha}$ for each $\alpha\in\Pi$.
By Lemma \ref{sl2s} there exist $\xi_\alpha\in k^\times$ such that $d\omega_\alpha(1)=\xi_\alpha[E_\alpha,d\theta(E_\alpha)]$ for each $\alpha$.

Let $\kappa$ be a non-degenerate $(\theta,G)$-equivariant symmetric bilinear form on ${\mathfrak g}$, let $S$ be a maximal torus of $G$ containing $A$, and let ${\mathfrak s}=\Lie(S)$.
By $S$-equivariance, the restriction of $\kappa$ to ${\mathfrak s}$ is non-degenerate; by $\theta$-equivariance, the restriction to ${\mathfrak a}$ is also non-degenerate.
Let $a\in{\mathfrak a}$.
Then $\kappa(a,d\omega_\alpha(1))=\xi_\alpha d\alpha(a)\kappa(E_\alpha,d\theta(E_\alpha))$.
Since $\kappa|_{{\mathfrak a}\times{\mathfrak a}}$ is non-degenerate, $\kappa(E_\alpha,d\theta(E_\alpha))\neq 0$ and the isomorphism ${\mathfrak a}\rightarrow{\mathfrak a}^*$ induced by $\kappa$ sends $d\omega_\alpha(1)$ to a non-zero multiple of $d\alpha$.
By Lemma \ref{diffindpt}, the toral elements $d\omega_\alpha(1)$ are linearly independent.
%\qed
\end{proof}

\subsection{Optimal cocharacters and $Y_\omega$.}
\label{sec:6.2}

Let $H$ be a reductive algebraic group, and let $\rho:H\longrightarrow\GL(V)$ be a rational representation.
Recall that $v\in V$ is $H$-unstable if $0\in\overline{\rho(H)(v)}$: otherwise $v$ is {\it $H$-semistable}.
Note that the $H$-unstable elements are the points of $\pi_{V,H}^{-1}(\pi_{V,H}(0))$.
We have the Hilbert-Mumford criterion (see \cite{mum}, for example):

{\it - $v$ is $H$-unstable if and only if there exists a cocharacter $\lambda:k^\times\longrightarrow H$ such that $v$ is $\lambda(k^\times)$-unstable.}

Let $T$ be a maximal torus of $H$, and let $W_T=N_H(T)/T$.
Let $Y(T)$ be the lattice of cocharacters in $T$ and let $E^*=Y(T)\otimes_{\mathbb Z}{\mathbb R}$.
Let $(.\, ,.):Y(T)\times Y(T)\longrightarrow{\mathbb Z}$ be a $W_T$-equivariant, positive definite symmetric bilinear form, extended linearly to an inner product $(.\, ,.):E^*\times E^*\longrightarrow{\mathbb R}$.
There is a corresponding length function $||.||:E^*\longrightarrow{\mathbb R}^{\geq 0}$, $\lambda\mapsto(\lambda,\lambda)^{1/2}$.
Any cocharacter $\lambda:k^\times\longrightarrow H$ is $H$-conjugate to an element of $Y(T)$, hence we can describe the set of cocharacters in $H$ as the union $Y(H)=\cup Y(hTh^{-1})$.
Moreover, if $\lambda,\mu\in Y(T)$, then $\lambda$ and $\mu$ are $H$-conjugate if and only if they are $W_T$-conjugate.
It follows that the length function can be extended to an $H$-equivariant function $||.||:Y(H)\longrightarrow{\mathbb R}^{\geq 0}$.

Let $\lambda\in Y(H)$ and let $h\in H$.
We say that the limit $\lim_{t\rightarrow 0}\lambda(t)h\lambda(t^{-1})$ exists if the morphism $k^\times\rightarrow H$, $t\mapsto\lambda(t)h\lambda(t^{-1})$ can be extended to a morphism $\eta:k\rightarrow H$.
If $\eta$ exists then it is unique: we write $\lim_{t\rightarrow 0}\lambda(t)h\lambda(t^{-1})$ for the image $\eta(0)$.
We associate to any cocharacter $\lambda$ the following subgroups of $H$:
$$P(\lambda):=\{ h\in H\,|\,\lim_{t\rightarrow 0}\lambda(t)h\lambda(t^{-1})\;\mbox{exists}\},$$
$$U(\lambda):=\{ h\in H\,|\,\lim_{t\rightarrow 0}\lambda(t)h\lambda(t^{-1})= I_H\}, Z(\lambda)=Z_H(\lambda).$$
(Here $I_H$ is the identity element of $H$.)
Then $P(\lambda)$ is a parabolic subgroup of $H$ with Levi decomposition $P(\lambda)=Z(\lambda)U(\lambda)$.

For $\lambda\in Y(H)$ and $i\in{\mathbb Z}$ set $V(i;\lambda)=\{ v\in V\,|\,\rho(\lambda(t))(v)=t^i v\;\forall t\in k^\times\}$: hence $V=\oplus_{i\in{\mathbb Z}} V(i;\lambda)$.
Let $v\in V$, $v=\sum_{i\in{\mathbb Z}} v_i$, $v_i\in V(i;\lambda)$.
We write $m(v,\lambda)$ for the mininum $i\in{\mathbb Z}$ such that $v_i\neq 0$.
The (non-trivial) cocharacter $\lambda$ is {\it optimal} for $v$ if $m(v,\lambda)/||\lambda||\geq m(v,\mu)/||\mu||$ for all $0\neq\mu\in Y(H)$.
A cocharacter $\lambda$ is {\it primitive} if $\lambda/m\in Y(H)\Rightarrow m=\pm 1$.

The main result of the Kempf-Rousseau theory is the following (\cite{kempf,rousseau}):

\begin{theorem}[Kempf, Rousseau]
Let $v$ be an $H$-unstable element of $V$.

(a) There exists at least one optimal cocharacter $\lambda\in Y(H)$ for $v$.

(b) There is a parabolic subgroup $P(v)$ of $G$ such that $P(v)=P(\lambda)$ for any optimal cocharacter $\lambda$ for $v$.
The centralizer $Z_H(v)\subset P(v)$.

(c) Let $\Lambda_v$ be the set of all cocharacters in $H$ which are primitive and optimal for $v$.
Any two elements of $\Lambda_v$ are conjugate by an element of $P(v)$.
Each maximal torus of $P(\lambda)$ contains a unique element of $\Lambda_v$.
\end{theorem}

Let $T$ be a maximal torus of $H$, and let $\lambda\in Y(T)$.
We denote by $T^\lambda$ the subtorus of $T$ generated by all cocharacters $\mu$ with $(\lambda,\mu)=0$, and by $Z^\bot(\lambda)$ the subgroup of $Z(\lambda)$ generated by $Z(\lambda)^{(1)}$ and $T^\lambda$.
Then $Z^\bot(\lambda)$ is a closed subgroup of $Z(\lambda)$ of codimension 1, and is independent of the choice of maximal torus $T$ containing $\lambda$.
We have the following criterion for optimality (Kirwan \cite{kir}, Ness \cite{ness}):

\begin{proposition}[Kirwan, Ness]
Let $i\geq 1$, and let $v\in V(i;\lambda)$.
Then $\lambda$ is optimal for $v$ if and only if $v$ is $Z^\bot(\lambda)$-semistable.
\end{proposition}

Consider the adjoint representation $\Ad:G\longrightarrow \GL({\mathfrak g})$.
Here $x\in {\mathfrak g}$ is $G$-unstable if and only if it is nilpotent.
In \cite{premnil}, Premet showed that every nilpotent element $x\in{\mathfrak g}$ has a cocharacter $\lambda$ which is both optimal for and associated to $x$.
(In general optimality depends on the choice of length function on $Y(G)$.)
Let $\lambda$ be any associated cocharacter for $x$.
Then $\lambda$ is optimal for $x$, and either $\lambda$ or $\lambda/2$ is primitive (\cite[Thm. 2.3, Thm. 2.7]{premnil}).
On the other hand, if $\lambda$ is optimal for $x$ and $x\in{\mathfrak g}(2;\lambda)$, then $\lambda$ is an associated cocharacter for $x$ (\cite[Thm. 14]{mcninch3}).

Let $S$ be a maximal torus of $G$ containing $A$, and let $E=X(S)\otimes_{\mathbb Z}{\mathbb R}$.
By \cite[2.6(iv)]{rich2}, $S$ is $\theta$-stable.
Let $W_S=N_G(S)/S$, let $\Gamma$ be the group of automorphisms of $S$ generated by $W_S$ and $\theta$, and let $(.\, ,.):E\times E\longrightarrow{\mathbb R}$ be a $\Gamma$-equivariant inner product such that $(\alpha,\beta)\in{\mathbb Z}$ for all $\alpha,\beta\in X(S)$.
The inner product induces a $\Gamma$-equivariant isomorphism $E\rightarrow E^*$.
Moreover, $E^*$ identifies with $Y(S)\otimes_{\mathbb Z}{\mathbb R}$.
Hence we write $(.\, ,.)$ also for the induced inner product on $E^*$.
Let $E_-$ (resp. $E^*_-$) denote the $(-1)$ eigenspace in $E$ (resp. $E^*$).
Then $E_-$ (resp. $E^*_-$) can be identified with $X(A)\otimes_{\mathbb Z}{\mathbb R}$ (resp. $Y(A)\otimes_{\mathbb Z}{\mathbb R}$).
The isomorphism $E\rightarrow E^*$ restricts to a $W_A$-equivariant isomorphism $E_-\rightarrow E^*_-$.
Recall (\cite[\S 1]{vust}) that $Z_G(A)=M\cdot A$ (almost direct product), where $M=Z_K(A)^\circ$.
Clearly $Z_G(A)^{(1)}\subseteq M$.
Since $\omega$ is regular in $A$, $Z_G(\omega)=Z_G(A)$.
Let $A^\omega$ denote the subtorus of $A$ generated by all $\mu(k^\times)$, with $\mu\in Y(A)$ such that $(\mu,\omega)=0$.

\begin{lemma}\label{zbot}
$Z^\bot(\omega)=M\cdot A^\omega$.
\end{lemma}

\begin{proof}
Let $S_0=(S\cap K)^\circ$.
By $\theta$-equivariance, $(\mu,\omega)=0$ for all $\mu\in Y(S_0)$.
Hence $Z^\bot(\omega)$ contains $S_0\cdot Z_G(A)^{(1)}=M$.
The lemma now follows at once.
%\qed
\end{proof}

Let $\alpha\in\Pi$ and let $L_\alpha$ be the (Levi) subgroup of $G$ introduced in Sect. \ref{sec:6.1}.
Note that $Z_{L_\alpha}(\omega_\alpha)=Z_G(\omega)=M\cdot A$.
Let $Z_{L_\alpha}^\bot(\omega_\alpha)$ be the subgroup of $Z_G(A)$ generated by $Z_G(A)^{(1)}$ and $S^{\omega_\alpha}$ (using similar notation to that used above).

\begin{lemma}\label{zbot2}
(i) $Z_{L_\alpha}^\bot(\omega_\alpha)=M\cdot (Z(L_\alpha)\cap A)^\circ$.

(ii) Let $x_\alpha\in{\mathfrak g}(\alpha;A)$.
Then $x_\alpha\in Y_\alpha$ if and only if $x_\alpha$ is $M$-semistable.
\end{lemma}

\begin{proof}
By Lemma \ref{zbot} applied to $L_\alpha$, $Z_{L_\alpha}^\bot(\omega_\alpha)=M\cdot A^{\omega_\alpha}$.
But $A=(Z(L_\alpha)\cap A)^\circ\cdot\omega_\alpha(k^\times)$, hence (i) follows.
Part (ii) now follows from the Kirwan-Ness criterion.
%\qed
\end{proof}

For ease of notation, let $\pi_\alpha=\pi_{{\mathfrak g}(\alpha;A),M}$.
We can choose homogeneous generators $f_1,f_2,\ldots f_l$ for $k[{\mathfrak g}(\alpha;A)]^M$.
Let the respective degrees be $d_1,d_2,\ldots d_l$.
Recall (Rk. \ref{geoquot}) that there is a natural action of $A$ on ${\mathfrak g}(\alpha;A)\quot M$, induced by the action on ${\mathfrak g}(\alpha;A)$.
Clearly $a\cdot f_i=\alpha(a)^{-d_i}f_i$ for any $a\in A$.
Let $U_\alpha$ be a vector space with basis $u_1,u_2,\ldots u_l$ and let $A$ act on $U_\alpha$ by: $a\cdot u_i=\alpha(a)^{d_i}u_i$, extending linearly to all of $U_\alpha$.
Hence the morphism ${\mathfrak g}(\alpha;A)\longrightarrow U_\alpha$, $x_\alpha\mapsto \sum f_i(x_\alpha)u_i$ induces an $A$-equivariant embedding $\iota_\alpha:{\mathfrak g}(\alpha;A)\quot M\hookrightarrow U_\alpha$.
Since the $f_i$ are homogeneous, $\iota_\alpha(\pi_\alpha(0))=0$.
Hence $x_\alpha\in Y_\alpha$ if and only if $\iota_\alpha(\pi_\alpha(x_\alpha))\neq 0$ (by Lemma \ref{zbot2}).

Let $r_0=\rank\Phi_A^*$.
Embed $A$ diagonally in the product $Z_G(A)^{r_0}$, and let $H=M^{r_0}\subset Z_G(A)^{r_0}$.
Clearly $H$ commutes with $A$.
Let the coordinates of $Z_G(A)^{r_0}$ be indexed by the elements of $\Pi$, and let $Z_G(A)^{r_0}$ act on ${\mathfrak g}(2;\omega)=\oplus_{\alpha\in\Pi}{\mathfrak g}(\alpha;A)$: $(g_\alpha)\cdot\sum y_\alpha=\sum(g_\alpha\cdot y_\alpha)$.
It is easy to see that the quotient ${\mathfrak g}(2;\omega)\quot H$ is naturally isomorphic to $\prod_{\alpha\in\Pi}{\mathfrak g}(\alpha;A)\quot M$.
Identify ${\mathfrak g}(2;\omega)\quot H$ with $\prod_{\alpha\in\Pi}{\mathfrak g}(\alpha;A)\quot M$, let $U=\oplus_{\alpha\in\Pi}U_\alpha$, and let $\iota=(\prod \iota_\alpha):{\mathfrak g}(2;\omega)\quot H\longrightarrow U$.
Then $\iota$ is an $A$-equivariant embedding.
Hence by Rk. \ref{geoquot} the following diagram is commutative:
\begin{diagram}
{\mathfrak g}(2;\omega) & \rTo & {\mathfrak g}(2;\omega)\quot H & \rTo & U \\
\dTo & & \dTo & & \dTo \\
{\mathfrak g}(2;\omega)\quot A^\omega & \rTo & {\mathfrak g}(2;\omega)\quot A^\omega H & \rTo & U\quot A^\omega
\end{diagram}

(Note that by construction $\iota(\pi_{{\mathfrak g}(2;\omega),H}(0))=0$.)

\begin{lemma}\label{alphaunstable}
(i) Let $u\in U$.
Then $u$ is $A^\omega$-unstable if and only if $u_\alpha=0$ for some $\alpha\in\Pi$.

(ii) Let $x=\sum_{\alpha\in\Pi} x_\alpha\in{\mathfrak g}(2;\omega)$.
Then $x$ is $A^\omega H$-semistable if and only if $x_\alpha\in Y_\alpha$ for all $\alpha\in\Pi$.
\end{lemma}

\begin{proof}
Since $A^\omega=(A^\omega\cap G^{(1)})^\circ\cdot (Z(G)\cap A)^\circ$ and $(Z(G)\cap A)$ acts trivially on $U$, we may clearly assume that $G$ is semisimple.
Suppose that $u\in U$ is $A^\omega$-unstable.
By the Hilbert-Mumford criterion, there exists $\mu\in Y(A^\omega)$ such that $u$ is $\mu(k^\times)$-unstable.
After replacing $\mu$ by $-\mu$, if necessary, we may assume that $u\in\sum_{i\geq 1}U(i;\mu)$.
Note that $U_\alpha\subset\sum_{i\geq 1}U(i;\mu)$ if and only if $\langle\alpha,\mu\rangle>0$.
Hence if $u_\alpha\neq 0$ for all $\alpha$, then $\langle\alpha,\mu\rangle>0$ for all $\alpha\in\Pi$.
But this implies that $\mu$ and $\omega$ are in the same Weyl chamber in $Y(A)$, which contradicts the assumption that $(\mu,\omega)=0$.

Suppose therefore that $u_\alpha=0$ for some $\alpha\in\Pi$.
Recall that $\omega=\omega_\alpha+\mu_\alpha$ for some $\mu_\alpha\in Y(Z(L_\alpha))$.
Hence $(\omega,\omega)=(\omega_\alpha,\omega_\alpha)+(\mu_\alpha,\mu_\alpha)$ and $(\omega_\alpha,\omega)=(\omega_\alpha,\omega_\alpha)$.
It follows that $c=(\omega_\alpha,\omega)/{(\omega,\omega)}<1$.
Let $m\in{\mathbb N}$ be such that $\nu=m(\omega_\alpha-c\omega)\in Y(A)$.
Then in fact $\nu\in Y(A^\omega)$.
Moreover, $\langle\alpha,\nu\rangle>0$ and $\langle\beta,\nu\rangle<0$ for all $\beta\in\Pi\setminus\{\alpha\}$.
Hence $u$ is $\nu(k^\times)$-semistable.
This proves (i).

For ease of notation, let $V={\mathfrak g}(2;\omega)$ and let $V_\alpha={\mathfrak g}(\alpha;A)$.
Suppose $x=\sum x_\alpha\in V$.
Recall (Rk. \ref{geoquot}) that $\pi_{V,A^\omega H}=\pi_{V,A^\omega H/H}\circ\pi_{V,H}$.
Moreover, $V\quot H$ embeds as an $A$-stable subset of $U$.
It follows that $x$ is an $A^\omega H$-unstable element of $V$ if and only if $\iota(\pi_{V,H}(x))$ is an $A^\omega$-unstable element of $U$.
But by (i), this holds if and only if $\iota_\alpha(\pi_\alpha(x_\alpha))=0$ for some $\alpha\in\Pi$.
Hence, by Lemma \ref{zbot2} $x$ is $A^\omega H$-semistable if and only if $x_\alpha\in Y_\alpha$ for all $\alpha$.
%\qed
\end{proof}

\begin{corollary}
Let $x\in {\mathfrak g}(2;\omega)$ be such that $x_\alpha\in Y_\alpha$ for all $\alpha\in\Pi$.
Then $x\in Y_\omega$.
\end{corollary}

\begin{proof}
By the Kirwan-Ness criterion, $x\in Y_\omega$ if and only if $x$ is $Z^\bot(\omega)$-semistable.
If $x$ is $Z^\bot(\omega)$-unstable, then it is clearly also $A^\omega H$-unstable.
But then $x_\alpha\notin Y_\alpha$ for some $\alpha\in\Pi$ by Lemma \ref{alphaunstable}.
%\qed
\end{proof}

Hence we have the following equivalent conditions:

\begin{proposition}\label{yomega}
Let $x=\sum_{\alpha\in\Pi}x_\alpha\in{\mathfrak g}(2;\omega)$.
Then the following are equivalent:

(i) $x\in Y_\omega$,

(ii) $[{\mathfrak g}^\omega,x]={\mathfrak g}(2;\omega)$,

(iii) $x_\alpha\in Y_\alpha$ for each $\alpha\in\Pi$,

(iv) $[{\mathfrak g}^\omega,x_\alpha]={\mathfrak g}(\alpha;A)$ for each $\alpha\in\Pi$.
\end{proposition}

\begin{proof}
The equivalence of (i) and (ii) is an immediate consequence of the separability of orbits, (see Lemma \ref{globalinf}).
Hence (iii) and (iv) are also equivalent (since $L_\alpha$ is a Levi subgroup of $G$).
Suppose $[{\mathfrak g}^\omega,x]={\mathfrak g}(2;\omega)$.
Then $\oplus_{\alpha\in\Pi}[{\mathfrak g}^\omega,x_\alpha]={\mathfrak g}(2;\omega)$, hence $[{\mathfrak g}^\omega,x_\alpha]={\mathfrak g}(\alpha;A)$.
This shows that (ii) $\Rightarrow$ (iv).
But by Lemma \ref{alphaunstable}, (iv) $\Rightarrow$ (ii).
This completes the proof.
%\qed
\end{proof}

\begin{rk}
The above proposition differs slightly from \cite[Prop. 19]{kostrall}, which it seeks to imitate.
Kostant-Rallis' version considers only elements of ${\mathfrak g}(2;\omega)$ which are contained in the real form ${\mathfrak g}_{\mathbb R}$.
Then $x\in Y_\omega\cap{\mathfrak g}_{\mathbb R}$ if and only if $x_\alpha\neq 0$ for each $\alpha$.
\end{rk}

\subsection{Construction of ${\mathfrak g}^*$}
\label{sec:6.3}

In \cite{kostrall}, Kostant and Rallis constructed a reductive subalgebra ${\mathfrak g}^*$ of ${\mathfrak g}$ containing ${\mathfrak a}$ as a Cartan subalgebra.
We will now generalise this to positive characteristic.
Fix $E\in Y_\omega$ and let $d\omega_\alpha(1)=H_\alpha=\xi_\alpha [E_\alpha,d\theta(E_\alpha)]$.
Let $F_\alpha=\xi_\alpha d\theta(E_\alpha)$.
Hence $\{ H_\alpha,E_\alpha,F_\alpha\}$ is an $\mathfrak{sl}(2)$-triple for each $\alpha$.

\begin{lemma}\label{commrels}
We have the following relations:

(a) $[H_\alpha,H_\beta]=0$ $(\alpha,\beta\in\Pi)$,

(b) $[H_\alpha,E_\beta]=-\langle\beta,\alpha\rangle E_\beta$ $(\alpha,\beta\in\Pi)$,

(c) $[H_\alpha,F_\beta]=\langle\beta,\alpha\rangle F_\beta$ $(\alpha,\beta\in\Pi)$,

(d) $[E_\alpha,F_\beta]=0$ for $\alpha\neq\beta\in\Pi$,

(e) $(\ad E_\alpha)^{-\langle\beta,\alpha\rangle+1}(E_\beta)=(\ad F_\alpha)^{-\langle\beta,\alpha\rangle+1}(F_\beta)=0$ for $\alpha\neq\beta\in\Pi$,

(f) $E_\alpha^{[p]}=F_\alpha^{[p]}=0$, and $H_\alpha^{[p]}=H_\alpha$ for every $\alpha\in\Pi$.
\end{lemma}

\begin{proof}
(a) is immediate since $H_\alpha\in{\mathfrak a}$; (b) and (c) follow from Lemma \ref{omegaalphacartan}.
If $\alpha\neq\beta\in\Pi$, then $\alpha-\beta\notin\Phi_A$.
Hence (d) follows.
Clearly, $\beta+m\alpha\in\Phi_A^*$ $\Leftrightarrow$ $\beta+m\alpha\in\Phi_A$.
But the integers $\langle\beta,\alpha\rangle$ are the Cartan integers for $\Phi_A^*$.
Hence $\beta+(1-\langle\beta,\alpha\rangle)\notin\Phi_A^*$, which proves (e).
Finally, if $\alpha\in\Phi_A$ then $3\alpha\notin\Phi_A$ by Lemma \ref{pisgood}.
Hence $E_\alpha^{[p]}=F_\alpha^{[p]}=0$.
Since $H_\alpha=d\omega_\alpha(1)$, $H_\alpha$ is a toral element.
This proves (f).
%\qed
\end{proof}

\begin{proposition}\label{constrpre}
Let ${\mathfrak b}^*={\mathfrak b}^*(E)$ be the subalgebra of ${\mathfrak g}$ generated by the elements $E_\alpha,F_\alpha,H_\alpha$.
Then ${\mathfrak b}^*$ is a $d\theta$-stable restricted subalgebra of ${\mathfrak g}$, ${\mathfrak a}\cap\Lie(G^{(1)})$ is a Cartan subalgebra of ${\mathfrak b}^*$, $[{\mathfrak b}^*,{\mathfrak b}^*]={\mathfrak b}^*$, and ${\mathfrak b}^*$ is an almost classical Lie algebra of universal type with root system $\Phi_A^*$.
Hence there exists a simply-connected semisimple group $B^*$ such that $\Lie(B^*)={\mathfrak b}^*$.
\end{proposition}

\begin{proof}
Since the set $\{ H_\alpha,E_\alpha,F_\alpha\}$ is $d\theta$-stable, so is ${\mathfrak b}^*$.
Furthermore, $E_\alpha^{[p]}=F_\alpha^{[p]}=0$ and $H_\alpha^{[p]}=H_\alpha$ by Lemma \ref{commrels}(f).
It follows that ${\mathfrak b}^*$ is a restricted subalgebra of ${\mathfrak g}$.
Let $G_{(1)},G_{(2)},\ldots,G_{(l)}$ be the minimal $\theta$-stable normal subgroups of $G^{(1)}$ and let ${\mathfrak g}_{(1)}=\Lie(G_{(1)}),{\mathfrak g}_{(2)}=\Lie(G_{(2)}),\ldots,{\mathfrak g}_{(l)}=\Lie(G_{(l)})$.
Hence $\Lie(G^{(1)})={\mathfrak g}_{(1)}\oplus{\mathfrak g}_{(2)}\oplus\ldots\oplus{\mathfrak g}_{(l)}$.
Moreover $\Phi_A^*\cong \Phi_{(1)}^*\cup\Phi_{(2)}^*\cup\ldots\cup\Phi_{(l)}^*$ is the decomposition of the root system into simple components, where $\Phi_{(i)}^*=\Phi(G_{(i)},A\cap G_{(i)})^*$.
Thus ${\mathfrak b}^*={\mathfrak b}_{(1)}^*\oplus{\mathfrak b}_{(2)}^*\oplus\ldots\oplus{\mathfrak b}_{(l)}^*$, where ${\mathfrak b}^*_{(i)}={\mathfrak b}^*\cap{\mathfrak g}_{(i)}$.
But therefore we have only to prove the proposition in the case $G=G_{(1)}$.
Hence we may assume that $\Phi_A^*$ is irreducible.

Let $\{ H_\alpha^{\mathbb C},E_\beta^{\mathbb C},F_\beta^{\mathbb C}:\alpha\in\Pi,\beta\in(\Phi_A^*)^+\}$ be a Chevalley basis for a complex semisimple Lie algebra ${\mathfrak g}_{\mathbb C}$ with root system $\Phi_A^*$.
Let ${\mathfrak g}_{\mathbb Z}$ be the ${\mathbb Z}$-subalgebra spanned by the elements $H_\alpha^{\mathbb C},E_\beta^{\mathbb C},F_\beta^{\mathbb C}$.
The $k$-Lie algebra ${\mathfrak g}_{\mathbb Z}\otimes_{\mathbb Z} k$ is an almost classical Lie algebra of universal type, and it is generated by $\{ H_\alpha^{\mathbb C}\otimes 1,E_\alpha^{\mathbb C}\otimes 1,F_\alpha^{\mathbb C}\otimes 1:\alpha\in\Pi\}$.
Hence by Lemma \ref{commrels} there is a unique Lie algebra homomorphism $\phi:{\mathfrak g}_{\mathbb Z}\otimes k\longrightarrow{\mathfrak b}^*$ such that $H_\alpha^{\mathbb C}\otimes 1\mapsto H_\alpha,E_\alpha^{\mathbb C}\otimes 1\mapsto E_\alpha,F^{\mathbb C}_\alpha\otimes 1\mapsto F_\alpha$.
Since ${\mathfrak b}^*$ is generated by the elements $E_\alpha,F_\alpha,\alpha\in\Pi$, $\phi$ is surjective.
The ideals of ${\mathfrak g}_{\mathbb Z}\otimes k$ are given in \cite[p. 446-7]{hog}.
Since $p$ is good, there is only one case of a non-trivial ideal: when $\Phi_A^*$ is of type $A_n$ and $p|(n+1)$, the centre is of dimension 1.
But by Cor. \ref{omegasindpt} the elements $H_\alpha$, $\alpha\in\Pi$ are linearly independent.
Hence $\phi$ is injective in all cases.
Thus ${\mathfrak b}^*\cong {\mathfrak g}_{\mathbb Z}\otimes k$.
Since ${\mathfrak b}^*$ is of universal type, there exists a simply-connected semisimple group $B^*$ such that $\Lie(B^*)={\mathfrak b}^*$ (see the discussion in \cite[\S 1]{hog}).
It remains to show that ${\mathfrak a}\cap\Lie(G^{(1)})$ is a Cartan subalgebra of ${\mathfrak b}^*$.
But by Cor. \ref{omegasindpt}, ${\mathfrak a}\cap\Lie(G^{(1)})$ is spanned by $H_\alpha\, ,\alpha\in\Pi$.
%\qed
\end{proof}

\begin{lemma}\label{wa}
Let ${\mathfrak a}'=\Lie(A\cap G^{(1)})$ and let $W^*=N_{B^*}({\mathfrak a}')/Z_{B^*}({\mathfrak a}')$.
Then $W_A=N_G({\mathfrak a})/Z_G({\mathfrak a})$ is naturally isomorphic to $W^*$.
\end{lemma}

\begin{proof}
Since the root system of $B^*$ is identified with $\Phi_A^*$, $N_{B^*}({\mathfrak a})/Z_{B^*}({\mathfrak a})$ is generated by the reflections $s_\alpha,\alpha\in\Pi$.
But so is $W_A$ by \cite[4.5]{rich2}.
%\qed
\end{proof}

We are now ready to present the main theorem of this section:

\begin{theorem}\label{constr}
Let $E\in Y_\omega$ and let ${\mathfrak g}^*(E)$ be the Lie subalgebra of ${\mathfrak g}$ generated by $E,d\theta(E)$ and ${\mathfrak a}$.

(a) ${\mathfrak g}^*(E)$ is a $d\theta$-stable restricted subalgebra of ${\mathfrak g}$, $[{\mathfrak g}^*(E),{\mathfrak g}^*(E)]={\mathfrak b}^*(E)$, and ${\mathfrak a}$ is a maximal toral algebra in ${\mathfrak g}^*(E)$.

(b) There exists a reductive group $G^*$ satisfying the standard hypotheses (A)-(C) of \S 3, such that $\Lie(G^*)={\mathfrak g}^*(E)$.

(c) There is an involutive automorphism $\theta^*$ of $G^*$ such that $d\theta^*=d\theta|_{\mathfrak g}$.
\end{theorem}

\begin{proof}
Let ${\mathfrak g}^*={\mathfrak g}^*(E),{\mathfrak b}^*={\mathfrak b}^*(E)$.
Since $[{\mathfrak a},E]=\sum_{\alpha\in\Pi}kE_\alpha$ and $[{\mathfrak a},d\theta(E)]$ $=\sum_{\alpha\in\Pi}kd\theta(E_\alpha)$, ${\mathfrak g}^*$ contains ${\mathfrak b}^*$.
Moreover, $[{\mathfrak b}^*,{\mathfrak b}^*]={\mathfrak b}^*$ by Prop. \ref{constrpre} and ${\mathfrak a}$ normalizes ${\mathfrak b}^*$.
Hence ${\mathfrak b}^*=[{\mathfrak g}^*,{\mathfrak g}^*]$.
Clearly ${\mathfrak g}^*$ is generated by ${\mathfrak a}$ and ${\mathfrak b}^*$.
Therefore ${\mathfrak g}^*$ is $d\theta$-stable and closed under the $p$-operation.
This proves (a).

By Prop. \ref{constrpre}, ${\mathfrak b}^*=\Lie(B^*)$, where $B^*$ is a simply-connected semisimple group.
Let $S$ be a maximal torus of $G$ containing $A$ and let $\Delta_S$ be a basis for $\Phi(G,S)$ such that $\Pi$ can be obtained as $\{\beta|_A\,:\,\beta\in\Delta_S\}$.
Let $S'=S\cap G^{(1)},A'=A\cap G^{(1)},{\mathfrak a}'=\Lie(A')$.
Since $G^{(1)}$ is simply-connected, $Y(S')=\oplus_{\beta\in\Delta_S}\beta^\vee$, where $\beta^\vee$ denotes the coroot corresponding to $\beta$.
Let $\alpha\in\Pi$ and let $\beta\in\Delta_S$ be such that $\beta|_A=\alpha$.
There are three possibilities: (i) $\theta^*(\beta)=-\beta$, (ii), $-\theta^*(\beta)$ and $\beta$ are orthogonal, and (iii) $-\theta^*(\beta)$ and $\beta$ generate a root system of type $A_2$.
But now we can describe $\omega_\alpha$ explicitly: in case (i), $\omega_\alpha=\beta^\vee$; in (ii) $\omega_\alpha=\beta^\vee-\theta^*(\beta)^\vee$; and in case (iii), $\omega_\alpha=2(\beta^\vee-\theta^*(\beta)^\vee)$.
Let $c_\alpha=1$ if $\alpha$ is of type (i) or (ii), and $c_\alpha=2$ if $\alpha$ is of type (iii).
It follows from Lemma \ref{basis} that $\{ \omega_\alpha/c_\alpha\,:\,\alpha\in\Pi\}$ is a basis for $Y(A')$.

Let $A_B^*$ be the unique maximal torus of $B^*$ such that $\Lie(A_B^*)={\mathfrak a}'$ (Lemma \ref{maxsplittori}).
Then $Y(A_B^*)$ can be identified with $\oplus_{\alpha\in\Pi}{\mathbb Z}\omega_\alpha\subset Y(A')$.
Hence $Y(A_B^*)$ embeds as a sublattice of $Y(A')$ of index $2^i$, where $i$ is the number of roots in $\Pi$ which are of type (iii).
Let $\{\chi_\alpha\,:\alpha\in\Pi\}$ be the basis for $X(A')$ which is dual to the basis $\{\omega_\alpha/c_\alpha\,:\,\alpha\in\Pi\}$ for $Y(A')$.
Then we can identify $X(A_B^*)$ with $\oplus_{\alpha\in\Pi}{\mathbb Z}(\chi_\alpha/c_\alpha)\subset X(A')\otimes_{\mathbb Z}{\mathbb Q}$.
Clearly $X(A')$ is a sublattice of $X(A_B^*)$ of index $2^i$.
Now the basis $\{\chi_\alpha\}$ can be lifted to a basis $\{\hat{\chi}_\alpha,z_j\,:\,\alpha\in\Pi\, ,1\leq j\leq r-r_0\}$ for $X(A)$.
(Here $r=\dim A$ and $r_0=\rank\Phi_A^*$.)
Let $\Lambda_X=\oplus_{\alpha\in\Pi}{\mathbb Z}(\hat{\chi}_\alpha/c_\alpha)\oplus{\mathbb Z}z_1\oplus\ldots\oplus{\mathbb Z}z_{r-r_0}\subset X(A)\otimes_{\mathbb Z}{\mathbb Q}$.
Clearly $\Lambda_X$ contains $X(A)$ as a sublattice of index $2^i$.
The pairing $\langle .\, ,.\rangle:X(A)\times Y(A)\longrightarrow{\mathbb Z}$ can be extended to a ${\mathbb Z}$-bilinear map $\langle .\, ,.\rangle:\Lambda_X\times Y(A)\longrightarrow{\mathbb Q}$.
Let $\Lambda_Y=\{\lambda\in Y(A)\,|\,\langle\chi,\lambda\rangle\in{\mathbb Z}\;\forall\chi\in\Lambda_X\}$.
Then $\Lambda_Y$ is a sublattice of $Y(A)$ of index $2^i$.

Let $A^*$ be the torus with character lattice $\Lambda_X$, that is $A^*=\Spec (k\Lambda_X)$.
Then $A^*$ contains $A_B^*$.
Since $\Lambda_Y$ is of index $2^i$ in $Y(A)$, we can identify $\Lie(A^*)$ with ${\mathfrak a}$.
Set $G^*=(B^*\times A^*)/{\diag(A_B^*)}$.
It is easy to see that $G^*$ is reductive and that $\Lie(G^*)$ can be identified with ${\mathfrak g}^*$.
To prove (b) we therefore have only to show that the restriction to ${\mathfrak g}^*$ of the $d\theta$-equivariant trace form $\kappa$ (see Cor. \ref{trace}) is non-degenerate.

Let ${\mathfrak s}=\Lie(S)$.
Since $\kappa$ is non-degenerate its restriction to ${\mathfrak s}$ is non-degenerate.
But $\kappa$ is also $d\theta$-equivariant.
Hence $\kappa(s,a)=0$ for any $s\in{\mathfrak s}\cap{\mathfrak k}$ and any $a\in{\mathfrak a}$.
It follows that the restriction $\kappa|_{{\mathfrak a}\times{\mathfrak a}}$ is non-degenerate.
To show that $\kappa|_{{\mathfrak g}^*}$ is non-degenerate, it will therefore suffice to show that the restriction to ${\mathfrak g}^*_\alpha\times{\mathfrak g}^*_{-\alpha}$ is non-degenerate for every $\alpha\in\Phi^*_A$.
(Here ${\mathfrak g}^*_\alpha={\mathfrak g}(\alpha;A)\cap{\mathfrak g}^*$, a one-dimensional root subspace for each $\alpha\in\Phi_A^*$).
But the Weyl group of $G^*$ is isomorphic to $W_A$ by Lemma \ref{wa}.
Hence to see that the restriction of $\kappa$ to ${\mathfrak g}^*$ is non-degenerate, we require only that $\kappa(E_\alpha,F_\alpha)\neq 0$ for each $\alpha\in\Pi$.
Since $\kappa$ is non-degenerate on ${\mathfrak a}$, there exists $a\in{\mathfrak a}$ such that $\kappa(a,H_\alpha)\neq 0$.
But $\kappa(a,H_\alpha)=d\alpha(a)\kappa(E_\alpha,F_\alpha)\neq 0$.
Hence $\kappa|_{{\mathfrak g}^*\times{\mathfrak g}^*}$ is non-degenerate.

Since $B^*$ is simply-connected, there exists a unique automorphism $\theta_B^*$ of $B^*$ such that $d\theta_B^*=d\theta|_{{\mathfrak b}^*}$ by Lemma \ref{sccover}.
Hence the involutive automorphism of $B^*\times A^*$ given by $(g,a)\mapsto (\theta_B^*(g),a^{-1})$ induces an automorphism $\theta^*$ of $G^*=(B^*\times A^*)/\diag(A_B^*)$ satisfying $d\theta^*=d\theta|_{{\mathfrak g}^*}$.
%\qed
\end{proof}

As an immediate consequence of the theorem, all of our earlier results apply to the pair $(G^*,\theta^*)$.

\begin{rk}
It is possible to construct a group $G^*_0$ such that $\Lie(G^*_0)={\mathfrak g}^*$ and $A$ is a maximal torus of $G^*_0$.
It is clear from the proof of Thm. \ref{constr} that the universal covering of $(G^*_0)^{(1)}$ is isomorphic to $B^*$, and that $B^*\rightarrow (G^*_0)^{(1)}$ is separable, with kernel of order $2^i$.
Here $i$ is the number of roots $\alpha\in\Pi$ which are of type (iii) (that is, if $\beta\in\Delta_S$ satisfies $\beta|_A=\alpha$, then $\beta$ and $-\theta^*(\beta)$ generate a root system of type $A_2$).
It can be seen from the classification of involutions (proved in odd
characteristic by Springer \cite{springer}) that there is at most one root of type (iii) for each component of the root system of $G$.
Suppose $G$ is almost simple, hence so is $G^*(=B^*)$.
Since the universal covering $G^*\rightarrow G_0^*$ maps $Z(G^*)$ onto $Z(G)\cap A$, we can easily calculate the order of $Z(G)\cap A$ for an arbitrary involution.
\end{rk}

\begin{lemma}\label{cpts2}
Let $G$ be an almost simple (simply-connected) group.

(1) Suppose $\theta$ is quasi-split, but not split.

(a) ${\cal N}$ has two irreducible components if $G$ is of type $A_{2n+1}$ or $D_{2n+1}$.

(b) Otherwise ${\cal N}$ is irreducible (types $A_{2n},D_{2n},E_6$).

(2) Let $\theta$ be an involution which is neither split nor quasi-split.
If $G$ is of type $A,E_6,E_8,F_4$, or if $\theta$ is an outer involution in type $D$, then ${\cal N}$ is irreducible.
\end{lemma}

\begin{proof}
(1) Let $Z=Z(G)$.
Recall from Cor. \ref{splitcmpts} that the components of ${\cal N}$
are in one-to-one correspondence with the elements of $Z\cap
A/\tau(Z)$, where $\tau:G\rightarrow G$ is given by $g\mapsto g^{-1}\theta(g)$.
Suppose $G$ is of type $E_6$.
Then $Z$ is a cyclic group of order 3, hence $(Z\cap A)/(Z\cap A)^2$ is trivial.
By Thm. \ref{gthetaorbs}, ${\cal N}$ is irreducible.
Similarly, ${\cal N}$ is irreducible if $G$ is of type $A_{2n}$.
For $G$ of type $A_{2n+1}$ (resp. $D_{2n+1},D_{2n}$) we can see from \cite[pp. 664-665]{springer} that $\Phi_A^*$ is of type $C_{n+1}$ (resp. $B_{2n-1},B_{2n-1}$).
Hence $Z(G^*)$ is of order 2 in each case.
Unless $G$ is of type $D_{2n}$, $\theta$ is inner by \cite{springer}, hence $\theta(z)= z$ for any $z\in Z(G)$.
On the other hand, an outer automorphism acts non-trivially on the centre.
It follows that $\tau(Z)$ is trivial unless $G$ is of type $D_{2n}$, in which case it is of order 2.
This shows that ${\cal N}$ has the number of irreducible components indicated.

(2) If $G$ is of type $E_6,E_8$, or $F_4$, then $Z/Z^2$ is trivial, hence ${\cal N}$ is irreducible by Thm. \ref{gthetaorbs}.
For an inner automorphism in type $A$, $\Phi_A^*$ is of type $C$, hence $Z(G^*)$ is of order 2.
Moreover, there exists a root $\alpha\in\Pi$ of type (iii); hence
$Z\cap A$ is trivial.
It follows that ${\cal N}$ is irreducible.
Suppose $\theta$ is a non-split outer automorphism in type $A_{2n+1}$.
Then $\Phi_A^*$ is of type $A_n$, and there is no root of type (iii).
Therefore $Z\cap A$ is of order $(n+1)$.
But (since $\theta$ is outer) we have $z\mapsto z^{-1}$ for $z\in Z$.
Thus $\tau(Z)= Z^2$ is of order $(2n+2)/2=(n+1)$.
Therefore $Z\cap A=\tau(Z)$, which implies that ${\cal N}$ is irreducible.

Finally, suppose $\theta$ is an outer involution in type $D$.
Then $\Phi_A^*$ is of type $B$, hence $Z(G^*)$ is of order 2.
There is no root of type (iii), hence $Z\cap A$ is also of order 2.
But $\theta$ acts non-trivially on the centre, hence $\tau(Z)\neq 1$.
It follows that $Z\cap A/{\tau(Z)}$ is trivial.
%\qed
\end{proof}

Lemma \ref{cpts2} provides us with two more classes of involution for
which ${\cal N}$ has two irreducible components: the quasi-split
involutions in type $A_{2n+1}$ and $D_{2n+1}$ are, respectively
$(\mathfrak{gl}(2n+2),\mathfrak{gl}(n+1)\oplus\mathfrak{gl}(n+1))$
and $(\mathfrak{so}(4n+2),\mathfrak{so}(2n+2)\oplus\mathfrak{so}(2n))$.

We now check the remaining (non-quasi-split) cases.
The classification of involutions in \cite{springer} associates to
each class of involution a unique {\it Araki diagram}: the Araki diagram for $\theta$ is a copy of the Dynkin diagram on $\Delta_S$, with the action of $\psi$ indicated, and the vertices in $I$ (resp. $\Delta_S\setminus I$) coloured black (resp. white).
But then one can easily write down the weighted Dynkin diagram corresponding to $\omega$ (and hence to a regular nilpotent element of ${\mathfrak p}$): $h(\alpha)=2$ if $\alpha\in\Delta_S\setminus I$, and $h(\alpha)=0$ if $\alpha\in I$.
Lemma \ref{cpts2} and \cite{springer} reduce us to the following cases:

(i) Non-split involutions in type $B_n$.
Here there are $(n-1)$ classes of involution, with corresponding
weighted Dynking diagrams $$\begin{array}{llll}
2 & 0 & \cdots & 0
\end{array},\;\;\;\; \begin{array}{lllll}
2 & 2 & 0 & \cdots & 0 
\end{array},\;\;\;\;\ldots \;\;\; ,\;\;\;
\begin{array}{llll}
2 & \cdots & 2 & 0
\end{array}$$
In each case $\Phi_A^*$ is of type $B$, and there is no root
$\alpha\in\Pi$ of type (iii).
Hence $Z\cap A$ is of order 2.
For type $B$ it is easier to carry out the calculations in the adjoint
group $\SO(2n+1)$, which we embed in the standard way in $\SL(2n+1)$.
Let $e$ be a regular nilpotent element of ${\mathfrak p}$ and let $C$
be its `reductive part'.
The determination of the number of irreducible components of ${\cal
N}$ therefore comes down to the determination of whether $C$ is
contained in $K$ or not.
(Here $G^\theta/K$ is of order 2.)
The embedding of $G$ in $\SL(2n+1)$ allows us to classify the nilpotent orbits in ${\mathfrak g}$ by partitions of $(2n+1)$, see for example \cite[3.5]{som}.
(The only partitions which occur in type $B$ are those such that $i$ appears an even number of times if $i$ is even.)
The partitions of $(2n+1)$ corresponding to the above weighted Dynkin
diagrams are, respectively, $3^1.1^{2(n-1)},\; 5^1.1^{2(n-2)},\;\ldots ,(2n-1)^1.1^2$.

The pair corresponding to a weighted Dynkin diagram as above with $m$
2's is
$(\mathfrak{so}(2n+1),\mathfrak{so}(m)\oplus\mathfrak{so}(2n+1-m))$.
It follows that if $m$ is even and $e$ is a regular nilpotent element
of ${\mathfrak p}$, then $\theta$ is conjugate to
$\Ad\lambda(\sqrt{-1})$, where $\lambda$ is an associated cocharacter
for $e$.
But then $Z_G(\lambda)\subset K$, hence $C\subset K$.
It follows that in this case, ${\cal N}$ has two irreducible
components.

Suppose therefore that $m$ is odd.
It is easy to see that $K\cong\SO(m)\times\SO(2n+1-m)$, and that
$G^\theta\cong\{ (g,h)\in\O(m)\times\O(2n+1-m)\, |\,\det g=\det h\}$.
Here $C/C^\circ$ is of order 2 by Sommers' theorem.
In fact, we can see by direct calculation that $C\cong \O(2n+1-m)$,
and that $C/C^\circ$ is generated by an element of
$G^\theta\hookrightarrow\O(m)\times\O(2n+1-m)$ of the form $(-I,n)$,
where $\det n=-1$.
But therefore $CK=G^\theta$.
It follows that ${\cal N}$ is irreducible in this case.

(ii) Non-split involutions in type $C_n$.
We consider $G=\Sp(2n,k)$ as a subgroup of $\SL(2n,k)$ in the standard
way.
There are $[n/2]$ classes of non-split involution of $G$, with
corresponding weighted Dynkin diagrams 
$$\begin{array}{llllll}
2 & 0 & 0 & \cdots & 0 \end{array},\;\;\;
\begin{array}{llllll}
2 & 0 & 2 & 0 & \cdots & 0 \end{array},\;\;\;\ldots\;\;\;,\;\;\;
\left\{
\begin{array}{cc}
{\begin{array}{llllll}
2 & 0 & 2 & \cdots & 0 & 2 \end{array}} & \mbox{if $n$ is even,} \\
{\begin{array}{llllll}
2 & 0 & 2 & \cdots & 2 & 0 \end{array}} & \mbox{if $n$ is odd.}
\end{array}\right.$$

In each case, the roots $\Phi_A^*$ are of type $B$, and with the
exception of the case $\begin{array}{llllll}
2 & 0 & 2 & \cdots & 0 & 2
\end{array}$, there is a root $\alpha\in\Pi$ of type (iii).
This shows that $Z\cap A$ is trivial in each except this final case,
which is $(\mathfrak{sp}(4n),\mathfrak{sp}(2n)\oplus\mathfrak{sp}(2n))$.
Here a regular nilpotent element of ${\mathfrak p}$ is of partition
type $(2n)^2$.
Up to conjugacy, $\theta$ is equal to conjugation by $\Int\left(\begin{matrix}
A_0 & {} & 0 \\
{} & \ddots & {} \\
0 & {} & A_0
\end{matrix}\right)$, where $A_0=\left(\begin{matrix}
1 & {} & {} & 0 \\
{} & -1 & {} & {} \\
{} & {} & -1 & {} \\
0 & {} & {} & 1
\end{matrix}\right)$.

Then $e=e_{13}+e_{24}+\ldots-e_{2n-2,2n}\in{\mathfrak
g}$ is a regular nilpotent element of ${\mathfrak p}$.
Hence if $\lambda$ is the unique diagonal cocharacter which is
associated to $e$, then $c=\left(\begin{matrix}
n & {} & 0 \\
{} & \ddots & {} \\
0 & {} & n
\end{matrix}\right)\in Z_G(\lambda)\cap Z_G(e)$ and $c^{-1}\theta(c)=-1$, where
$n=\left(\begin{matrix}
0 & 1 \\
1 & 0
\end{matrix}\right)$.
It follows from Cor. \ref{splitcmpts} that ${\cal N}$ is irreducible in this case.

(iii) Inner involutions in type $D_{2n}$.
There are $(n+1)$ classes of involutions, producing $\Phi_A^*$ of types $B_2,B_4,\ldots ,B_{2n-2},C_n,C_n$.
The corresponding weighted Dynkin diagrams are:
$$\begin{array}{llllll}
{} & {} & {} & {} & {} & 0 \\
2 & 2 & 0 & \cdots & 0 & {} \\
{} & {} & {} & {} & {} & 0
\end{array} ,\;\;\;\;
\begin{array}{llllllll}
{} & {} & {} & {} & {} & {} & {} & 0 \\
2 & 2 & 2 & 2 & 0 & \cdots & 0 & {} \\
{} & {} & {} & {} & {} & {} & {} & 0 \\
\end{array}, \;\;\;\;
\ldots \;\;\;\; ,\;\;
\begin{array}{llllll}
{} & {} & {} & {} & {} & 0 \\
2 & 2 & 2 & \cdots & 2 & {} \\
{} & {} & {} & {} & {} & 0
\end{array},$$
$$\begin{array}{llllll}
{} & {} & {} & {} & {} & 2 \\
0 & 2 & 0 & \cdots & 2 & {} \\
{} & {} & {} & {} & {} & 0
\end{array},\;
\begin{array}{llllll}
{} & {} & {} & {} & {} & 0 \\
0 & 2 & 0 & \cdots & 2 & {} \\
{} & {} & {} & {} & {} & 2
\end{array}$$

Moreover, we have respectively:
${\mathfrak
k}=\mathfrak{so}(4n-2)\oplus\mathfrak{so}(2),\,\mathfrak{so}(4n-4)\oplus\mathfrak{so}(4),\,\ldots
\,
,\mathfrak{so}(2n+2)\oplus\mathfrak{so}(2n-2),\,\mathfrak{gl}(2n),\mathfrak{gl}(2n)$.
(The final two cases are conjugate by an outer involution of $G$.)
The nilpotent orbits in ${\mathfrak g}$ are classified in a standard way by partitions of $4n$, see for example \cite[3.5]{som}.
(The only partitions which occur in type $D$ are those such that $i$ appears an even number of times if $i$ is even.)
The partitions corresponding to the above weighted Dynkin diagrams are $3^1.1^{4n-3},\; 7^1.1^{4n-7},\ldots ,(4n-5)^1.1^5,(2n)^2,(2n)^2$.
Hence by Sommers' theorem \cite{som,premnil,mcninchsom}, in each of these cases the group $C=Z_G(\lambda)\cap Z_G(e)$ is connected modulo $Z(G)$.
(Here $e$ is a regular nilpotent element of ${\mathfrak p}$ and $\lambda$ is an associated cocharacter for $e$.)
Moreover, there is no root of type (iii).
Hence $Z\cap A/\tau(C) = Z\cap A\cong Z(G^*)$.
Thus ${\cal N}$ has two irreducible components.

(iv) Inner involutions in type $D_{2n+1}$.
There are $n$ classes of involutions, producing $\Phi_A^*$ of types $B_2,B_4,\ldots ,B_{2n-2},B_n$.
The corresponding weighted Dynkin diagrams are:
$$\begin{array}{llllll}
{} & {} & {} & {} & {} & 0 \\
2 & 2 & 0 & \cdots & 0 & {} \\
{} & {} & {} & {} & {} & 0
\end{array},\;\;\;
\begin{array}{llllllll}
{} & {} & {} & {} & {} & {} & {} & 0 \\
2 & 2 & 2 & 2 & 0 & \cdots & 0 & {} \\
{} & {} & {} & {} & {} & {} & {} & 0 \\
\end{array}, \;\;\;\;
\ldots\;\;\;\; ,\;\;
\begin{array}{llllll}
{} & {} & {} & {} & {} & 0 \\
2 & 2 & \cdots & 2 & 0 & {} \\
{} & {} & {} & {} & {} & 0
\end{array},$$
$$\mbox{and}\;\;
\begin{array}{llllll}
{} & {} & {} & {} & {} & 2 \\
0 & 2 & 0 & \cdots & 0 & {} \\
{} & {} & {} & {} & {} & 2
\end{array}$$

We have, respectively: ${\mathfrak
k}=\mathfrak{so}(4n)\oplus\mathfrak{so}(2),\mathfrak{so}(4n-2)\oplus\mathfrak{so}(4),\,\ldots\,
,\mathfrak{so}(2n+4)\oplus\mathfrak{so}(2n-2)$, and $\mathfrak{gl}(2n+1)$.
In the final case $\theta^*(\alpha_{2n})=-(\alpha_{2n-1}+\alpha_{2n+1})$.
Thus $\alpha_{2n}|_A=\alpha_{2n+1}|_A$ is of type (iii), $\Rightarrow\;A\cap Z(G)$ is trivial $\Rightarrow\;{\cal N}$ is irreducible.
For the first $(n-1)$ diagrams, the corresponding partitions of $(4n+2)$ are: $3^1.1^{4n-1},\; 7^1.1^{4n-5},\ldots ,$ $(4n-5)^1.1^7$.
By Sommers' theorem $Z_G(\lambda)\cap Z_G(e)$ is connected modulo $Z(G)$ in each case.
It follows that ${\cal N}$ has two irreducible components.

(v) (Inner) involutions in type $E_7$.
Here there are two classes of involutions, with weighted Dynkin diagrams:
$$\begin{array}{llllll}
2 & 2 & 2 & 0 & 2 & 0 \\
{} & {} & 0 & {} & {} & {}
\end{array}
\;\;\;\;\;\mbox{and}\;\;\;\;\;
\begin{array}{llllll}
2 & 0 & 0 & 0 & 2 & 2 \\
{} & {} & 0 {} & {} & {}
\end{array}.$$

For the first class, which is $({\mathfrak
e}_7,\mathfrak{so}(12)\oplus\mathfrak{sl}(2))$, $\Phi_A^*$ is of type $F_4$, hence ${\cal N}$ is irreducible (since the fundamental group of $F_4$ is trivial).
For the second, which is $({\mathfrak e}_7,{\mathfrak e}_6\oplus k)$, $\Phi_A^*$ is of type $C_3$ and there is no root of type (iii).
Hence $(Z\cap A)/\tau(Z)$ is of order 2.
Moreover, by Sommers' theorem (\cite[p. 558]{som} and \cite{premnil,mcninchsom}) $Z_G(\lambda)\cap Z_G(e)$ is connected modulo $Z(G)$.
Therefore ${\cal N}$ has two irreducible components.

This completes the process of computing the number of
irreducible components of ${\cal N}$.
The non-irreducible cases match those given by Sekiguchi in \cite{sek} for
$k={\mathbb C}$.

\begin{proposition}
The classes of involution for which ${\cal N}$ is non-irreducible are
as follows.

 - Type $A$: $(\mathfrak{gl}(n),\mathfrak{so}(n))$, $(\mathfrak{gl}(2n),\mathfrak{gl}(n)\oplus\mathfrak{gl}(n))$.

 - Type $B$:
   $(\mathfrak{so}(2n+1),\mathfrak{so}(2m)\oplus\mathfrak{so}(2(n-m)+1))$,
   {\bf only} if the even part $2m < 2(n-m)+1$,

 - Type $C$: $(\mathfrak{sp}(2n),\mathfrak{gl}(n))$,

 - Type $D$:
   $(\mathfrak{so}(2n),\mathfrak{so}(2m)\oplus\mathfrak{so}(2(n-m))$,
   $(\mathfrak{so}(4n),\mathfrak{gl}(2n))$,
   $(\mathfrak{so}(4n+2),\mathfrak{so}(2n+1)\oplus\mathfrak{so}(2n+1))$,

 - Type $E_7$:
   $(\mathfrak{e}_7,\mathfrak{sl}(8)),(\mathfrak{e}_7,\mathfrak{e}_6\oplus k)$.

In each of these cases ${\cal N}$ has two irreducible components, except for
$(\mathfrak{so}(4n),\mathfrak{so}(2n)\oplus\mathfrak{so}(2n))$, where
there are four components.
\end{proposition}

\subsection{Applications}
\label{sec:6.4}

We draw a number of conclusions from Theorem \ref{constr}.
Let $S$ be a maximal torus of $G$ containing $A$, and let $\Delta_S$ be a basis for $\Phi_S$ from which $\Pi$ is obtained (see Sect. \ref{sec:2.2}).
We can now show that each fibre of the quotient morphism $\pi_{\mathfrak p}:{\mathfrak p}\longrightarrow{\mathfrak p}\quot K$ has a dense open $K^*$-orbit.

\begin{lemma}\label{fibrelemma}
Let $s\in{\mathfrak a}$, let $L=Z_G(s)^\circ$, and let ${\mathfrak l}=\Lie(L)$.
There is a dense open $(K^*\cap L)$-orbit in ${\cal N}({\mathfrak l}\cap{\mathfrak p})$.
\end{lemma}

\begin{proof}
Since $s\in{\mathfrak a}$, $L$ is a $\theta$-stable Levi subgroup of $G$ containing $A$.
Let $F_L^*=\{ a\in A\,|\,a^2\in Z(L)\}$.
As there is a surjective map from $F_L^*/{F(Z(L)\cap A)}$ to the set of $L^\theta$-orbits in ${\cal N}({\mathfrak l}\cap{\mathfrak p})_{reg}$ (Thm. \ref{gthetaorbs}), it will suffice to show that the map $F^*/{(Z\cap A)}\rightarrow F_L^*/{(Z(L)\cap A)}$ induced by the embedding $F^*\hookrightarrow F_L^*$ is surjective.
Let $r_0=\rank (A\cap G^{(1)})$.
The basis $\Pi=\{\alpha_1,\alpha_2,\ldots,\alpha_{r_0}\}$ determines an isomorphism $(\alpha_1,\alpha_2,\ldots,\alpha_{r_0}):A/{(Z\cap A)}\longrightarrow (k^\times)^{r_0}$.
(Separability follows from Lemma \ref{diffindpt}.)
The subgroup $F/{(Z\cap A)}$ maps onto the set of $r_0$-tuples of the form $(\pm 1,\ldots,\pm 1)$.
Since any Levi subgroup of $G^*$ is conjugate to a standard Levi subgroup, there exists $w\in W(G^*,A^*)$ such that $w(Z_{G^*}(s))$ is standard.
But $W(G^*,A^*)=W_A$ by Lemma \ref{wa}.
Hence, after replacing $s$ by some $W_A$-conjugate, if necessary, there is a subset $J\subseteq \Pi$ such that ${\mathfrak l}$ is spanned by ${\mathfrak g}^A$ and the subspaces ${\mathfrak g}(\alpha;A)$ with $\alpha\in {\mathbb Z}J\cap \Phi_A$.
Then $J=\{\beta_1,\beta_2,\ldots,\beta_{r_1}\}$ determines an isomorphism $(\beta_1,\beta_2,\ldots,\beta_{r_1}):A/{(Z(L)\cap A)}\rightarrow (k^\times)^{r_1}$.
It is now easy to see that the projection onto the $\beta_i$-coordinates gives a surjective homomorphism $A/{(Z\cap A)}\rightarrow A/{(Z(L)\cap A)}$ which sends $F^*/{(Z\cap A)}$ onto $F_L^*/{(Z(L)\cap A)}$.
%\qed
\end{proof}

Hence:

\begin{theorem}\label{fibres}
Every fibre of $\pi_{\mathfrak p}$ contains a dense (open) $K^*$-orbit.
\end{theorem}

\begin{proof}
Let $\xi\in{\mathfrak p}\quot K$ and let $s$ be a semisimple element of $\pi_{\mathfrak p}^{-1}(\xi)$.
We may assume after conjugating by an element of $K$, if necessary, that $s\in{\mathfrak a}$.
Let $L=Z_G(s)=Z_G(s)^\circ,{\mathfrak l}=\Lie(L)$.
Thus $\pi^{-1}(\xi)=K\cdot \{s+{\cal N}({\mathfrak l}\cap{\mathfrak p})\}$.
By Lemma \ref{fibrelemma} there is an open $K^*\cap L$-orbit in ${\cal N}({\mathfrak l}\cap{\mathfrak p})$.
Hence there is a dense $K^*$-orbit in $\pi^{-1}(\xi)$.
%\qed
\end{proof}

\begin{rk}
Let $P=\{ g^{-1}\theta(g)\,|\,g\in G\}$.
Let $x\in G$ and let $x=su$ be the Jordan-Chevalley decomposition of $x$, where $s$ is the semisimple part and $u$ the unipotent part.
Then $x\in P$ if and only if $\theta(u)=u^{-1}$ and $s$ is contained in a maximal $\theta$-split torus of $G$ (\cite[6.1]{rich2}).
Let ${\cal U}$ denote the set of unipotent elements in $P$; recall (Cor. \ref{isocor}) that there is a $K^*$-equivariant isomorphism $\Psi:{\cal U}\longrightarrow{\cal N}$.
Fix a maximal $\theta$-split torus $A$ of $G$.
By \cite[11.3-4]{rich2} the action of $K^*$ on $P$ is well-defined and the embedding $A\hookrightarrow P$ induces an isomorphism $A/W_A\longrightarrow P\quot K\cong P\quot K^*$.
Hence each fibre of $\pi_P:P\longrightarrow P\quot K$ is $K^*$-stable and contains a unique closed (semisimple) $K$-orbit.
In \cite[Rk. 10.4]{rich2} Richardson conjectured that each fibre of $\pi_P:P\rightarrow P\quot K$ has a dense open $K^*$-orbit.
However, this is not true, as we now show.

It follows from the above that every fibre of $\pi_P$ can be written as $K\cdot a({\cal U}\cap Z_G(a))=K\cdot a({\cal U}\cap Z_G(a)^\circ)$ for some $a\in A$.
Let $a\in A$, let $L=Z_G(a)^\circ$ and let $V_1,V_2,\ldots,V_l$ be the irreducible components of ${\cal U}\cap L$.
By Cor. \ref{isocor} and Lemma \ref{nil1} the $V_i$ are of equal dimension, and each contains an open $(L^\theta)^\circ$-orbit which is just the intersection with the set of $\theta$-regular elements of $L$.
(An element $x\in P$ is $\theta$-regular if $\dim Z_G(x)=\dim Z_G(A)$.
Note that $v\in V_i$ is $\theta$-regular in $L$ if and only if $av$ is $\theta$-regular in $G$.)
It follows that each irreducible component of $\pi_P^{-1}(\pi_P(a))$ is of the form $K\cdot aV_i$ for some $i$, and $K\cdot aV_i$ is an irreducible component of $\pi_P^{-1}(\pi_P(a))$ for all $i$.

It is now easy to see that there is a dense open $K^*$-orbit in $\pi_P^{-1}(\pi_P(a))$ if and only if $Z_G(a)\cap K^*$ permutes the components $V_i$ transitively.
Let $G$ be almost simple, of type $E_8,F_4$, or $G_2$, and let $\theta$ be a split involution of $G$.
Since $G$ is both simply-connected and adjoint, $K^*=G^\theta=K$ and $L=Z_G(a)=Z_G(a)^\circ$.
It follows that there is a dense open $K^*$-orbit in $\pi_P^{-1}(\pi_P(a))$ if and only if $L^\theta$ permutes the components of ${\cal U}\cap L$ transitively.
Let $a$ be a non-regular element of order 2.
As $G$ is adjoint, $Z(L)/Z(L)^\circ$ is cyclic of order 2 (see \cite[Prop. 3.2]{premnil}).
Hence $Z(L)/Z(L)^2$ is cyclic of order 2.
By Cor. \ref{splitcmpts}, the $L^\theta$-orbits in ${\cal N}({\mathfrak l}\cap{\mathfrak p})$ are parametrised by the elements of $Z(L)/Z(L)^2$.
Hence by \ref{isocor} there are two regular $L^\theta$-orbits in ${\cal U}\cap L$.
It follows that there is more than one regular $K^*$-orbit in $\pi_P^{-1}(\pi_P(a))$.
\end{rk}

Let $x\in{\mathfrak g}$ be such that $x^{[p]}=0$.
McNinch has associated to $x$ a family of {\it optimal} homomorphisms $\rho:\SL(2)\longrightarrow G$.
These behave in a similar way to the $\mathfrak{sl}(2)$-triples in zero (or large) characteristic.
Let $\chi:k^\times\longrightarrow\SL(2)$, $\chi(t)=
\left(
\begin{array}{ll}
t & 0 \\
0 & t^{-1}
\end{array}\right)$, let $X=
\left(
\begin{array}{ll}
0 & 1 \\
0 & 0
\end{array}\right)$, and let $Y=
\left(
\begin{array}{ll}
0 & 0 \\
1 & 0
\end{array}\right)$.
A homomorphism $\rho:\SL(2)\longrightarrow G$ is optimal for $x$ if $d\rho(X)=x$, and $\rho\circ\chi$ is an associated cocharacter for $x$ in $G$.
We have the following facts:

{\it - Optimal homomorphisms exist: for any associated cocharacter $\lambda$ for $x$ there is a unique homomorphism $\rho:\SL(2)\longrightarrow G$ such that $d\rho(X)=x$ and $\rho\circ\chi=\lambda$ (\cite{mcninch} and \cite[Prop. 44]{mcninch2}).

 - Any two optimal $\SL(2)$-homomorphisms for $x$ are conjugate by an element of $Z_G(x)^\circ$. (\cite[Thm. 46]{mcninch2}),

 - $Z_G(x)\cap Z_G(\lambda)=Z_G(\rho(\SL(2))$ (\cite[Cor. 45]{mcninch2}).}

Recall that a homomorphism $\rho:\SL(2)\longrightarrow G$ is {\it good} (cf. Seitz \cite{seitz}) if all weights of $\rho\circ\chi$ on ${\mathfrak g}$ are less than or equal to $(2p-2)$.

{\it - A homomorphism $\rho:\SL(2)\longrightarrow G$ is optimal for some $x$ if and only if it is good (\cite[Prop. 55]{mcninch2}),

 - The representation $(\Ad\circ\rho,{\mathfrak g})$ is a tilting module for $\SL(2)$ (This follows from \cite[Prop. 4.2]{seitz}. See \cite[Prop. 36 and Pf. of Prop. 37]{mcninch2}).}

Let $E,\omega$ be as in Thm. \ref{constr} and let ${\mathfrak g}^*={\mathfrak g}^*(E)$.
Let $\alpha\in\Pi$: then $E_\alpha^{[p]}=0$ by Lemma \ref{pisgood}.
Moreover, $\omega_\alpha$ is an associated cocharacter for $E_\alpha$ in $L_\alpha$.
But $L_\alpha$ is a Levi subgroup of $G$, hence $\omega_\alpha$ is associated to $E_\alpha$ in $G$.
Let $L^*_\alpha$ be the (unique) Levi subgroup of $G^*$ such that $\Lie(L^*_\alpha)= {\mathfrak a}\oplus kE_\alpha\oplus kd\theta(E_\alpha)$.
Then $E_\alpha$ is distinguished in $\Lie(L^*_\alpha)$.
By our construction of $G^*$ (see the proof of Thm. \ref{constr}) $\omega_\alpha$ also defines a cocharacter in $A^*$.
Hence $\omega_\alpha(k^\times)\subset (L^*_\alpha)^{(1)}$ by the argument used in the proof of Lemma \ref{omegaalphacartan}.
It follows that there exist optimal homomorphisms $\rho_\alpha:\SL(2)\longrightarrow G$ and $\rho'_\alpha:\SL(2)\longrightarrow G^*$ for $E_\alpha$ such that $\rho_\alpha\circ\chi=\omega_\alpha=\rho'_\alpha\circ\chi$.
By uniqueness, $\rho_\alpha(\SL(2))\subset L_\alpha$ and $\rho'_\alpha(\SL(2))\subset L^*_\alpha$.
By Lemma \ref{sl2s}, $\xi_\alpha d\theta(E_\alpha)$ is the unique element $F_\alpha\in{\mathfrak g}(-\alpha;A)$ such that $[E_\alpha,F_\alpha]=d\omega_\alpha(1)$.
Therefore $d\rho_\alpha(Y)=d\rho'_\alpha(Y)=F_\alpha$.
It follows that $d\rho_\alpha(x)=d\rho'_\alpha(x)$ for all $x\in\mathfrak{sl}(2)$.
Hence we can show:

\begin{lemma}
(i) ${\mathfrak g}^*$ is normalized by $\rho_\alpha(\SL(2))$.

(ii) $\Ad\rho_\alpha(g)|_{{\mathfrak g}^*}=\Ad\rho'_\alpha(g)$ for all $g\in \SL(2)$.

(iii) Let $H$ be the minimal closed subgroup of $G$ containing the subgroups $\rho_\alpha(\SL(2))$, $\alpha\in\Pi$.
Then $H$ is contained in $N_G({\mathfrak g}^*)$.

(iv) $\Ad H|_{{\mathfrak g}^*}=\Ad G^*$.
\end{lemma}

\begin{proof}
Let $\beta\in\Phi_A^*$, $\beta\neq\pm\alpha$, let $\beta-i\alpha,\ldots,\beta+j\alpha$ be the $\alpha$-chain through $\beta$, let ${\mathfrak g}_{(\beta)}={\mathfrak g}(\beta-i\alpha;A)\oplus\ldots\oplus{\mathfrak g}(\beta+j\alpha;A)$ and let $U={\mathfrak g}^A\oplus\sum{\mathfrak g}(\gamma;A)$, the sum taken over all $\gamma\in\Phi_A\setminus\{\beta-i\alpha,\ldots,\beta+j\alpha\}$.
Hence ${\mathfrak g}={\mathfrak g}_{(\beta)}\oplus U$ and each summand is $L_\alpha$-stable, therefore $\rho_\alpha(\SL(2))$-stable.
Since any direct summand in a tilting module is a tilting module (\cite[Thm. 1.1]{donkin2}), ${\mathfrak g}_{(\beta)}$ is a direct sum of indecomposable tilting modules for $\rho_\alpha(\SL(2))$.
For each positive integer $c$ there is a unique tilting module $T(c)$ for $\SL(2)$ with highest weight $c$: $T(c)$ is simple if $c<p$ (see \cite[Lemma 1.3]{seitz}).
But now by our condition on $p$, ${\mathfrak g}_{(\beta)}$ is a direct sum of simple $\rho_\alpha(\SL(2))$-modules.
Moreover, each tilting summand is infinitesimally irreducible, hence ${\mathfrak g}_{(\beta)}$ is completely reducible as a $\rho_\alpha(\SL(2))$-module, and as an $\mathfrak{sl}(2)$-module (with $\mathfrak{sl}(2)$ acting via $\ad\circ (d\rho_\alpha)$.
It follows that every $\mathfrak{sl}(2)$-submodule of ${\mathfrak g}_{(\beta)}$ is $\rho_\alpha(\SL(2))$-stable.

For $\gamma\in\Phi_A^*$, let ${\mathfrak g}^*_\gamma={\mathfrak g}^*\cap {\mathfrak g}(\gamma;A)$ (a one-dimensional root subspace), and let ${\mathfrak g}^*_{(\beta)}={\mathfrak g}^*_{\beta-i\alpha}\oplus\ldots\oplus{\mathfrak g}^*_{\beta+j\alpha}$.
Then ${\mathfrak g}^*_{(\beta)}$ is a simple $d\rho_\alpha(\mathfrak{sl}(2))$-submodule of ${\mathfrak g}_{(\beta)}$, hence is $\rho_\alpha(\SL(2))$-stable.
(In fact ${\mathfrak g}^*_{(\beta)}$ is isomorphic to $T(\langle\beta+j\alpha,\alpha\rangle)$.)
Moreover, ${\mathfrak g}^*={\mathfrak g}^*_{-\alpha}\oplus{\mathfrak a}\oplus{\mathfrak g}^*_\alpha\oplus\sum{\mathfrak g}^*_{(\beta)}$, and ${\mathfrak g}^*_{-\alpha}\oplus{\mathfrak a}\oplus{\mathfrak g}^*_\alpha= d\rho_\alpha(\mathfrak{sl}(2))\oplus({\mathfrak z}({\mathfrak l}_\alpha)\cap{\mathfrak a})$.
It follows that ${\mathfrak g}^*$ is $\rho_\alpha(\SL(2))$-stable.
This proves (i).
But now (iii) follows immediately.

We have decomposed ${\mathfrak g}^*$ as $\oplus V_\gamma$, where each $V_\gamma$ is a simple $d\rho_\alpha(\mathfrak{sl}(2))$-module of dimension $\leq 4$ ($\leq 3$ if $p=3$).
Each summand is also a simple tilting module for $\rho_\alpha(\SL(2))$ (resp. $\rho'_\alpha(\SL(2))$).
But now, since $d\rho_\alpha(x)=d\rho'_\alpha(x)$ for all $x\in\mathfrak{sl}(2)$, we must have: $\Ad\rho_\alpha(g)(v_\gamma)=\Ad\rho'_\alpha(g)(v_\gamma)$ for all $g\in\SL(2)$.
This proves (ii).
But $\Ad G^*$ is generated by the subgroups $\Ad\rho'_\alpha(\SL(2))$.
Hence (iv) follows.
%\qed
\end{proof}

\begin{corollary}\label{conjugacy}
For elements of ${\mathfrak g}^*$, $G^*$-conjugacy implies $G$-conjugacy.
\end{corollary}

Let ${\mathfrak k}^*={\mathfrak k}\cap{\mathfrak g}^*$, ${\mathfrak p}^*={\mathfrak p}\cap{\mathfrak g}^*$.
Clearly ${\mathfrak g}^*={\mathfrak k}^*\oplus{\mathfrak p}^*$ is the symmetric space decomposition of ${\mathfrak g}^*$.

\begin{lemma}\label{regstar}
Let $x\in{\mathfrak p}^*$.
The following are equivalent:

(i) $x$ is a ($\theta$-)regular element of ${\mathfrak p}$,

(ii) $x$ is a regular element of ${\mathfrak g}^*$,

(iii) ${\mathfrak z}_{{\mathfrak k}^*}(x)=0$,

(iv) $\dim{\mathfrak z}_{{\mathfrak p}^*}(x)=r=\dim{\mathfrak a}$.
\end{lemma}

\begin{proof}
Since ${\mathfrak a}$ is a maximal toral algebra of ${\mathfrak g}^*$, the equivalence of (ii)-(iv) follows immediately from Lemma \ref{regs}.
Suppose $x\in{\mathfrak p}^*$, and $x$ is a regular element of ${\mathfrak p}$.
Then $\dim{\mathfrak z}_{{\mathfrak p}^*}(x)\leq r$, hence (iv) holds.

It remains to show that if $x$ is a regular element of ${\mathfrak p}^*$, then $x$ is regular in ${\mathfrak p}$.
Let $e$ be a regular nilpotent element of ${\mathfrak p}^*$.
Then $e$ is $G^*$-conjugate to $E$.
But therefore $e$ is $G$-conjugate to $E$ by Cor. \ref{conjugacy}, hence $\dim{\mathfrak z}_{\mathfrak g}(e)=\dim{\mathfrak g}^\omega$, that is, $e$ is regular in ${\mathfrak p}$.
Suppose therefore that $x$ is a non-nilpotent regular element of ${\mathfrak p}^*$, and that $x=x_s+x_n$ is the Jordan-Chevalley decomposition of $x$.
After replacing $x$ by a $(G^*)^{\theta^*}$-conjugate, if necessary, we may assume that $x_s\in{\mathfrak a}$.
Let $L=Z_G(x_s),L^*=Z_{G^*}(x_s),{\mathfrak l}=\Lie(L),{\mathfrak l}^*=\Lie(L^*)$.
Let $\Pi_L$ be a basis for $\Phi(L,A)$, and let $\omega_L:k^\times\longrightarrow A\cap L^{(1)}$ be the unique cocharacter such that $\langle\alpha,\omega_L\rangle=2$ for all $\alpha\in\Pi_L$ (Cor. \ref{regconj}).
There exists a unique cocharacter $\omega^*_L:k^\times\longrightarrow A^*\cap (L^*)^{(1)}$ satisfying the same conditions: hence $\omega_L^*$ can be identified with $\omega_L$ (the embedding $Y(A^*)\hookrightarrow Y(A)$ sends $\omega_L^*$ to $\omega_L$).
We can therefore choose a representative $E_L$ for the open $Z_L(\omega_L)$-orbit in ${\mathfrak l}(2;\omega_L)$ such that $E_L\in{\mathfrak l}^*$.
Clearly $E_L$ is a regular nilpotent element of ${\mathfrak l}^*$.
By the argument used for Thm. \ref{constr}, ${\mathfrak l}^*$ is the subalgebra of ${\mathfrak l}$ generated by ${\mathfrak a},E_L$, and $d\theta(E_L)$.
Hence $L^*$ and $L$ stand in the same relation as do $G^*$ and $G$.

Since $x$ is regular in ${\mathfrak g}^*$, $x_n$ is a regular nilpotent element of ${\mathfrak l}^*$.
But then $x_n$ is $L^*$-conjugate to $E_L$, hence $L$-conjugate to $E_L$.
It follows that $\dim({\mathfrak l}\cap{\mathfrak z}_{\mathfrak g}(x_n))=\dim Z_G(A)$.
Thus $x$ is regular in ${\mathfrak p}$.
This completes the proof.
%\qed
\end{proof}

\begin{lemma}\label{gstarconj}
For semisimple elements of ${\mathfrak p}^*$, $(G^*)$-conjugacy is equivalent to $K$-conjugacy.
\end{lemma}

\begin{proof}
Let $a,a'$ be semisimple elements of ${\mathfrak p}^*$.
Since any two maximal tori of ${\mathfrak p}^*$ are conjugate by an element of $G^*$ (resp. $K$), we may clearly assume that $a,a'\in{\mathfrak a}$.
But now $a,a'$ are $K$-conjugate if and only if they are $W_A$-conjugate, hence if and only if they are $G^*$-conjugate.
%\qed
\end{proof}

Let $e$ be a nilpotent element of ${\mathfrak p}^*$ satisfying the equivalent conditions of Lemma \ref{regstar}.
By Lemma \ref{assoc} there is an associated cocharacter $\lambda:k^\times\longrightarrow (G^*)^{\theta^*}$ for $e$.
As $e$ is regular, ${\mathfrak z}_{{\mathfrak k}^*}(e)$ is trivial.
Therefore $[{\mathfrak p}^*,e]={\mathfrak k}^*$ and $[{\mathfrak k}^*,e]$ is of codimension $r=\dim{\mathfrak a}$ in ${\mathfrak p}^*$.
Let ${\mathfrak v}$ be an $\Ad\lambda$-graded subspace of ${\mathfrak p}^*$ such that $[{\mathfrak k}^*,e]\oplus{\mathfrak v}={\mathfrak p}^*$.
We recall (by \cite[6.3-6.5]{veldkamp}, see also \cite[\S 3]{premtang} for the proof in good characteristic) that every element of $e+{\mathfrak v}$ is regular in ${\mathfrak g}^*$, that the embedding $e+{\mathfrak v}\hookrightarrow{\mathfrak g}^*$ induces an isomorphism $e+{\mathfrak v}\rightarrow{\mathfrak g}^*\quot {G^*}$, and that each regular orbit in ${\mathfrak g}^*$ intersects $e+{\mathfrak v}$ in exactly one point.

\begin{lemma}\label{slices}
Let $j$ be the composite of the isomorphisms $k[{\mathfrak p}]^K\longrightarrow k[{\mathfrak a}]^{W_A}\longrightarrow k[{\mathfrak g}^*]^{G^*}$ and let $f\in k[{\mathfrak p}]^K,g\in k[{\mathfrak g}^*]^{G^*}$.
Then $j(f)=g$ if and only if $f|_{e+{\mathfrak v}}=g|_{e+{\mathfrak v}}$.
Hence ${\mathfrak p}\quot K$ is isomorphic to $e+{\mathfrak v}$, and each regular $K^*$-orbit in ${\mathfrak p}$ intersects $e+{\mathfrak v}$ in exactly one point.
\end{lemma}

\begin{proof}
Clearly $j(f)=g\Leftrightarrow f|_{\mathfrak a}=g|_{\mathfrak a}$.
The set of regular elements in ${\mathfrak a}$ is a dense open subset.
Hence its image in ${\mathfrak a}\quot W_A={\mathfrak a}/W_A$ is dense.
It follows that the set $U$ of semisimple elements in $e+{\mathfrak v}$ is dense.
By Lemma \ref{gstarconj}, $f|_{\mathfrak a}=g|_{\mathfrak a}\Leftrightarrow f|_U=g|_U\Leftrightarrow f|_{e+{\mathfrak v}}=g|_{e+{\mathfrak v}}$.
Therefore the restriction $k[{\mathfrak p}]\rightarrow k[e+{\mathfrak v}]$ induces an isomorphism $k[{\mathfrak p}]^K\longrightarrow k[e+{\mathfrak v}]$.

Let $x\in{\mathfrak p}$ be regular.
Then any regular element of $\pi_{\mathfrak p}^{-1}(\pi_{\mathfrak p}(x))$ is $K^*$-conjugate to $x$ by Thm. \ref{fibres}.
There is a unique point $y\in e+{\mathfrak v}$ such that $\pi(y)=\pi(x)$.
Moreover, $y$ is regular by Lemma \ref{regstar}.
This completes the proof.
%\qed
\end{proof}

\begin{corollary}\label{regcond}
Let $k[{\mathfrak p}]^K=k[u_1,u_2,\ldots,u_r]$, where the $u_i$ are homogeneous polynomials, and let $x\in{\mathfrak p}$ be regular. Then the differentials $(du_i)_x,1\leq i\leq r$ are linearly independent.
\end{corollary}

\begin{proof}
Let $x$ be regular.
By Lemma \ref{slices} there is a unique $K^*$-conjugate $y$ of $x$ in $e+{\mathfrak v}$.
Therefore the differentials $(du_i)_x$ are linearly independent if and only if $(du_i)_y$ are linearly independent (since $k[{\mathfrak p}]^K=k[{\mathfrak p}]^{K^*}$).
But the restriction map $k[{\mathfrak p}]^K\rightarrow k[e+{\mathfrak v}]$ is an isomorphism.
The result follows immediately since $e+{\mathfrak v}$ is isomorphic to affine $r$-space.
%\qed
\end{proof}

\begin{lemma}\label{regcodim}
The set ${\mathfrak p}\setminus{\cal R}$ of non-regular elements in ${\mathfrak p}$ is of codimension $\geq 2$.
\end{lemma}

\begin{proof}
Let ${\mathfrak a}_{reg}$ denote the set of regular elements in ${\mathfrak a}$.
Since ${\mathfrak a}\setminus{\mathfrak a}_{reg}$ is a union of hyperplanes in ${\mathfrak a}$, $U=\pi_{\mathfrak a}({\mathfrak a}\setminus{\mathfrak a}_{reg})$ is of pure codimension 1 in ${\mathfrak a}/W_A\cong {\mathfrak p}\quot K$.
Let $V=\pi_{\mathfrak p}^{-1}(U)$, the complement of the set of regular semisimple elements in ${\mathfrak p}$.
For any $x\in U$, the irreducible components of $\pi^{-1}(x)$ are of dimension $\dim{\mathfrak p}-\dim{\mathfrak a}$, hence $V$ is a closed set in ${\mathfrak p}$ of codimension greater than or equal to 1.
It is easy to see that ${\mathfrak p}\setminus{\cal R} = Y=V\setminus({\cal R}\cap V)$.
But $\pi_{\mathfrak p}(Y)=U$ and each fibre of $\pi_{\mathfrak p}|_Y$ has dimension strictly less than $\dim{\mathfrak p}-\dim{\mathfrak a}$.
It follows that $Y$ is of codimension $\geq 2$ in ${\mathfrak p}$.
%\qed
\end{proof}

We can now apply Skryabin's theorem on infinitesimal invariants.
The action of $K$ on the polynomial ring $k[{\mathfrak p}]$ induces an action of the Lie algebra ${\mathfrak k}$ as homogeneous derivations of $k[{\mathfrak p}]$.
We denote by $k[{\mathfrak p}]^{\mathfrak k}=\{ f\in k[{\mathfrak p}]\,|\,(x\cdot f)=0\;\forall x\in{\mathfrak k}\}$.
It is easy to see that $k[{\mathfrak p}]^{\mathfrak k}$ contains the global invariants $k[{\mathfrak p}]^K$.
Moreover, the ring of $p$-th powers, $k[{\mathfrak p}]^{(p)}=\{ f^p\,|\, f\in k[{\mathfrak p}]\}$ is also contained in $k[{\mathfrak p}]^{\mathfrak k}$.

\begin{theorem}
(1) (a) $k[{\mathfrak p}]^{\mathfrak k}=k[{\mathfrak p}]^K\cdot k[{\mathfrak p}]^{(p)}$ and $k[{\mathfrak p}]^{\mathfrak k}$ is free of rank $p^r$ over $k[{\mathfrak p}]^{(p)}$.

(b) $k[{\mathfrak p}]^{\mathfrak k}$ is a locally complete intersection.

(c) If $\pi_{{\mathfrak p},{\mathfrak k}}:{\mathfrak p}\longrightarrow {\mathfrak p}\quot{\mathfrak k}=\Spec(k[{\mathfrak p}]^{\mathfrak k})$ is the canonical morphism then $\pi_{{\mathfrak p},{\mathfrak k}}({\cal R})$ is the set of all smooth rational points of ${\mathfrak p}\quot{\mathfrak k}$.

(2) Let $K_i$ denote the $i$-th Frobenius kernel of $K$ and let $k[{\mathfrak p}]^{(p^i)}$ denote the ring of all $p^i$-th powers of elements of $k[{\mathfrak p}]$.

(a) $k[{\mathfrak p}]^{K_i}=k[{\mathfrak p}]^K\cdot k[{\mathfrak p}]^{(p^i)}$ and $k[{\mathfrak p}]^{K_i}$ is free of rank $p^{ir}$ over $k[{\mathfrak p}]^{(p^i)}$.

(b) $k[{\mathfrak p}]^{K_i}$ is a locally complete intersection.

(c) Let $\pi_{{\mathfrak p},K_i}:{\mathfrak p}\longrightarrow{\mathfrak p}\quot {K_i}$ denote the quotient morphism.
Then $\pi_{{\mathfrak p},K_i}({\cal R})$ is the set of all smooth rational points of ${\mathfrak p}\quot{K_i}$.
\end{theorem}

\begin{proof}
This follows immediately from Cor. \ref{regcond}, Lemma \ref{regcodim} and \cite[Thm.s 5.4,5.5]{skry}.
%\qed
\end{proof}

\end{document}